\documentclass{article}
\usepackage{amssymb}
\usepackage{amsmath}
\usepackage{enumitem}
\usepackage{tikz-cd}
\usepackage{hyperref}

\usepackage[titletoc]{appendix}
\usepackage{amsthm}
\usepackage{parskip}
\usepackage{amsfonts}
\usepackage{caption}
\usepackage{mathrsfs}
\usepackage{afterpage}
\usepackage{float}
\usepackage{multirow}
\usepackage{graphicx}
\usepackage{color}
\usepackage[margin=1.0in]{geometry}
\newcommand*{\SHom}{\mathscr{H}\kern -.5pt om}
\newcommand*{\SExt}{\mathscr{E}\kern -.5pt xt}

\begin{document}
\title{Quantum Steenrod operations and Fukaya categories}
\date{}
\author{Zihong Chen}
\maketitle
\theoremstyle{definition}
\newtheorem{mydef}{Definition}[section]
\numberwithin{mydef}{section}
\newtheorem{rmk}[mydef]{Remark}
\newtheorem{conj}[mydef]{Conjecture}
\theoremstyle{plain}
\newtheorem{cor}[mydef]{Corollary}
\newtheorem{lemma}[mydef]{Lemma}
\newtheorem{thm}[mydef]{Theorem}
\newtheorem{prop}[mydef]{Proposition}

\begin{abstract}
This paper is concerned with quantum cohomology and Fukaya categories of a closed monotone symplectic manifold $X$, where we use coefficients in a field $\mathbf{k}$ of characteristic $p>0$. The main result of this paper is that the quantum Steenrod operations $Q\Sigma$ admit an interpretation in terms of certain operations on the (equivariant) Hochschild invariants of the Fukaya category of $X$, via suitable (equivariant) versions of the open-closed maps. As an application, we demonstrate how the categorical perspective provides new tools for computing $Q\Sigma$ beyond the reach of known technology. We also explore potential connections of our work to arithmetic homological mirror symmetry.
\end{abstract}
\tableofcontents

\makeatletter
\def\thm@space@setup{%
  \thm@preskip=\parskip \thm@postskip=0pt
}
\makeatother

\renewcommand{\theequation}{1.\arabic{equation}}
\setcounter{equation}{0}
\section{Introduction}
\subsection{Motivation}
This paper is concerned with quantum cohomology and Fukaya categories with coefficients in a field $\mathbf{k}$ of characteristic $p>0$.
For simplicity, we work in the context of a closed monotone symplectic manifold $(X,\omega)$, where the technical aspects of the relevant Floer-theoretic structures are straightforward. Why is there a particular interest in working over positive characteristics?
\begin{enumerate}[label=\arabic*)]
    \item As a most straightforward answer, even if one is only interested in characteristic $0$ phenomena, reducing coefficients mod $p$ can be a useful tool in studying them, an idea that is ubiquitous and fruitful in different branches of mathematics.\par\indent
    To give an example of this idea in symplectic topology, an upcoming work of the author \cite{Che2} deploys a reduction mod $p$ argument to prove \emph{the exponential type conjecture} for quantum connections, cf. \cite[Conjecture 3.4]{KKP} and \cite[Section 2.5]{GGI}, on all closed monotone symplectic manifolds. 
    \begin{thm}(\cite{Che2})
    The quantum connection $(QH^*(X,\mathbb{C})((t)),\nabla^{QH}_{\frac{d}{dt}})$, cf. (1.11), admits a finite direct sum decomposition
    \begin{equation}
    \nabla^{QH}_{\frac{d}{dt}}=\bigoplus_{\lambda}(\mathbb{C}((t)),\frac{d}{dt}-\frac{\lambda}{t^2})\otimes \nabla^{reg}_{\frac{d}{dt}},
    \end{equation}
    where $\lambda\in\mathbb{C}$ and $\nabla^{reg}_{\frac{d}{dt}}$ is gauge equivalent to a connection with simple poles at $t=0$ and monodromies given by roots of unity. In fancier language, (1.1) says that the quantum connection has \emph{unramified exponential type} and \emph{quasi-unipotent regularized monodromy}.
    \end{thm}
   There are two main tools in the proof of Theorem 1.1. One is a classical result of Katz on the local monodromy of nilpotent connections \cite{Ka}. The other is the collection of \emph{quantum Steenrod operations} in Gromov-Witten theory with mod $p$ coefficients, which is also the main object of study in this paper.

\item The effectiveness of characteristic $p$ methods often comes from the fact that there are structures that only exist in characteristic $p$, such as the Steenrod operations in classical topology, the Frobenius map and Cartier operators in algebraic geometry, and the $p$-curvature in differential equations.\par\indent
The same phenomenon appears in symplectic topology, and as mentioned in 1), the characteristic $p$ structures we are concerned with are the quantum Steenrod operations $Q\Sigma: QH^*(X,\mathbf{k})\otimes QH^*_{\mathbb{Z}/p}(X,\mathbf{k})\rightarrow QH^*_{\mathbb{Z}/p}(X,\mathbf{k})$; cf. section 3.3 for their definition. In fact, quantum Steenrod operations are (or are expected to be) related to all of the structures mentioned in the previous paragraph:
\begin{itemize}
\item By construction, $Q\Sigma$'s are deformations of the classical Steenrod operations using mod $p$ counts of equivariant moduli spaces of pseudo-holomorphic curves. 
\item $Q\Sigma_{c_1}$ is related to the $p$-curvature of the quantum connection, cf. \cite[Theorem 1.2]{Lee} and \cite[Proposition 3.3]{Che2}. This idea was crucial to the proof of the exponential type conjecture in \cite{Che2}.   
\item We also expect $Q\Sigma$ to be `mirror' to certain operations on differential forms in characteristic $p$ on the $B$-side, cf. section 6 for a discussion. 
\end{itemize}

\item It is well known that Fukaya categories and homological mirror symmetry have interesting arithmetic aspects \cite{LPe}, \cite{LPo}, \cite{GHHPS}. Even simplistic symplectic manifolds such as the cylinder, when one considers its Fukaya category over a general coefficient ring $R$, can have interesting homological mirrors that encode the arithmetic properties of $R$, cf. \cite{EL}.\par\indent
As a concrete example, take $X$ to be the intersection of two quadrics in $\mathbb{C}P^5$. Over $\mathbb{C}$, a celebrated result of Ivan Smith \cite[Theorem 1.1]{Smi} states that a component of the Fukaya category of $X$ is equivalent to the Fukaya category of a genus $2$ curve. In section $5$, we show that there is a simple obstruction for the prior mentioned equivalence to hold over any coefficient field not containing a square root of $-1$, cf. Proposition 5.2. We do not know whether this is the only obstruction, i.e. whether Smith's equivalence holds over any coefficient field that contains a $\sqrt{-1}$. However, what we \emph{can} show is that if Smith's equivalence were true over some field $\mathbf{k}$ (containing a $\sqrt{-1}$) of odd positive characteristic, it imposes strong constraints on the structure of the quantum cohomology of $X$ over $\mathbf{k}$.  In fact, having such an equivalence will determine the quantum Steenrod operations of $X$ over $\mathbf{k}$, cf. Proposition 5.3. \par\indent
The idea behind Proposition 5.3, and the main result of this paper, is that quantum Steenrod operations admit an interpretation in terms of the Fukaya category. More precisely, for any $A_{\infty}$-category, there is an action of its Hochschild cohomology on its $\mathbb{Z}/p$-equivariant Hochschild homology, which we call the \emph{$\mathbb{Z}/p$-equivariant cap product}, cf. Definition 2.12. When the $A_{\infty}$-category in question is the Fukaya category of $X$, this action recovers $Q\Sigma$ via suitable versions of the open-closed maps. 
\end{enumerate}
With these motivations in mind, let us now state the main results of the paper.

\subsection{Quantum Steenrod operations}
Quantum Steenrod operations are a collection of operations on equivariant quantum cohomology that come from genus $0$ Gromov-Witten theory with mod $p$ coefficients. First introduced by Fukaya \cite{Fuk}, and systematically developed by Wilkins \cite{Wil}, they have seen various links and applications to Hamiltonian dynamics \cite{Shel}, arithmetic mirror symmetry \cite{Sei3} and representation theory \cite{Lee}. They have the property of being covariantly constant with respect to the quantum connection in characteristic $p$ \cite{SW}, which gives effective computation in low degrees. Nonetheless, our understanding of quantum Steenrod operations in the general case is still limited:
\begin{itemize}
    \item Relatively few computations of quantum Steenrod operations have been made beyond low degrees;
    \item The relationship between quantum Steenrod operations and Gromov-Witten theory in characteristic $0$ remains mysterious.
\end{itemize}
This paper aims to advance our understanding of both questions above by systematically studying quantum Steenrod operations on a closed monotone symplectic manifold from the perspective of its Fukaya category.  \par\indent
We now describe the setup and state the main results. Throughout this paper, $(X,\omega)$ will be a monotone symplectic manifold of dimension $2n$. Let $A$ be some base ring. Let $QH^*(X,A)$ denote the quantum cohomology of $X$ with coefficient in $A$; for $\lambda\in A$, let $\mathrm{Fuk}(X,A)_{\lambda}$ denote the monotone Fukaya category with disk potential $\lambda\in A$, see section 3 for a review. In this paper, the coefficient ring $A$ will come in two flavors:
\begin{enumerate}[label=\arabic*)]
    \item  Characteristic $0$: we consider $R\subset K\subset \overline{\mathbb{Q}}$ the ring of integers of some number field $K$, up to inverting finitely many elements, e.g. $\mathbb{Z}[\frac{1}{2}], \mathbb{Z}[\frac{1}{2}, \sqrt{5}, i]$. \emph{Throughout the paper, we always assume that $2$ is inverted in $R$}.
    \item Positive Characteristics: we consider, for an \emph{odd} prime $p$, fields $\mathbf{k}$ of characteristic $p$.
\end{enumerate}
Let us start with the story in characteristic $p$. In section 2, we will introduce a \emph{$p$-fold Hochschild chain complex} ${}_pCC(\mathcal{A})$ for any $A_{\infty}$-category $\mathcal{A}$ over $\mathbf{k}$. This chain complex is quasi-isomorphic to the usual Hochschild chain complex $CC(\mathcal{A})$, but has a chain level $\mathbb{Z}/p$-action, which can be thought of as induced from the $S^1$-action on $CC(\mathcal{A})$ via the inclusion of the $p$-th roots of unity $\mathbb{Z}/p\subset S^1$. We denote its associated negative $\mathbb{Z}/p$-equivariant complex (or $\mathbb{Z}/p$-homotopy fixed point) as $CC^{\mathbb{Z}/p}(\mathcal{A})$. There is an action
\begin{equation}
\bigcap^{\mathbb{Z}/p}: HH^*(\mathcal{A})\times HH^{\mathbb{Z}/p}_*(\mathcal{A})\rightarrow  HH^{\mathbb{Z}/p}_*(\mathcal{A}),
\end{equation}
which we call the \emph{$\mathbb{Z}/p$-equivariant cap product}. It has the following properties:
\begin{enumerate}[label=C\arabic*)]
    \item  It is a graded multiplicative action: $\bigcap^{\mathbb{Z}/p}(\phi\cup\varphi,a)=(-1)^{|\phi||\varphi|\frac{p(p-1)}{2}}\bigcap^{\mathbb{Z}/p}(\phi,\bigcap^{\mathbb{Z}/p}(\varphi,a))$, where $\cup$ denotes the cup product on Hochschild cohomology. 
    \item It is additive in the second variable, and becomes additive in the first variable after multiplying by $t$ (the formal $S^1$-equivariant variable of degree $2$), the latter meaning $t\bigcap^{\mathbb{Z}/p}(\phi+\varphi,-)=t\bigcap^{\mathbb{Z}/p}(\phi,-)+t\bigcap^{\mathbb{Z}/p}(\varphi,-)$.
    \item If $A$ is cohomologically unital, in which case $HH^*(A)$ is unital, then $\bigcap^{\mathbb{Z}/p}$ is a unital action.
    \item It is Frobenius $p$-linear: for $a\in \mathbf{k}$, $\bigcap^{\mathbb{Z}/p}(a\varphi,-)=a^p\bigcap^{\mathbb{Z}/p}(\varphi,-)$.
\end{enumerate}
In section 4, we define a \emph{$\mathbb{Z}/p$-equivariant open-closed map} (when context is clear, we omit $\lambda$ from notation)
\begin{equation}
OC^{\mathbb{Z}/p}_{\lambda}: HH^{\mathbb{Z}/p}_*(\mathrm{Fuk}(X,\mathbf{k})_{\lambda})\rightarrow QH^{*+n}_{\mathbb{Z}/p}(X,\mathbf{k}),
\end{equation}
where the $\mathbb{Z}/p$-action on $QH^*$ is trivial. The right hand side of (1.3) can be more explicitly written as
\begin{equation}
QH^*(X,\mathbf{k})\otimes H^*(B\mathbb{Z}/p,\mathbf{k})=QH^*(X,\mathbf{k})[[t,\theta]],
\end{equation}
where $|t|=2,|\theta|=1,\theta^2=0$. (1.3) is a $\mathbb{Z}/p$-analogue of the Ganatra's cyclic open-closed map, and it will be a crucial ingredient in proving the categorical formula for quantum Steenrod operations.   \par\indent
Following \cite{SW}, quantum Steenrod operations can be viewed as an action of $QH^*(X,\mathbf{k})$ on $QH^*(X,\mathbf{k})^{\mathbb{Z}/p}$. In particular, fixing a cohomology class $b\in QH^*(X,\mathbf{k})$ there is a quantum Steenrod action associated to $b$
\begin{equation}
Q\Sigma_b: QH^*(X,\mathbf{k})[[t,\theta]]\rightarrow QH^{*+p|b|}(X,\mathbf{k})[[t,\theta]].
\end{equation}
We review the definition and some of its properties in more detail in section 3. \par\indent
Our main result is that the $\mathbb{Z}/p$-equivariant open-closed map intertwines the quantum Steenrod operations (1.5) with the categorical action (1.2). 
\begin{thm}
In the above setting, for all $b\in QH^*(X,\mathbf{k})$ and $\lambda\in \mathbf{k}$, the diagram
\begin{equation}
\begin{tikzcd}[row sep=1.2cm, column sep=0.8cm]
HH^{\mathbb{Z}/p}_*(\mathrm{Fuk}(X,\mathbf{k})_{\lambda})\arrow[rrrr,"{OC^{\mathbb{Z}/p}}"]\arrow[d,"CO(b)\bigcap^{\mathbb{Z}/p}-"]& & & &QH^{*+n}(X,\mathbf{k})[[t,\theta]]\arrow[d,"Q\Sigma_b"]  \\
HH^{\mathbb{Z}/p}_{*+p|b|}(\mathrm{Fuk}(X,\mathbf{k})_{\lambda})\arrow[rrrr,"{OC^{\mathbb{Z}/p}}"]& & & &QH^{*+p|b|+n}(X,\mathbf{k})[[t,\theta]]
\end{tikzcd}
\end{equation}
commutes, where $CO$ denotes the closed-open map, cf. section 4.3. 
\end{thm}
\subsection{Relation to cyclic homology and the cyclic open-closed map}
Theorem 1.2 transforms the computation of $Q\Sigma$ into that of the $\mathbb{Z}/p$-equivariant cap product, which is purely algebraic, and of the $\mathbb{Z}/p$-equivariant open-closed map. Luckily, the latter is not mysterious and is in fact closely related to the familiar cyclic open-closed map, cf. \cite{Gan2}, which we now briefly recall.\par\indent
A cohomologically unital $A_{\infty}$-category $\mathcal{A}$ has a non-unital Hochschild chain complex $(CC^{nu}_*(\mathcal{A}),b^{nu})$, which is a chain complex quasi-isomorphic to $(CC_*(\mathcal{A}),b)$ and
\begin{equation}
CC^{nu}_*(\mathcal{A})=CC_*(\mathcal{A})\oplus CC_*(\mathcal{A})[1]
\end{equation}
as vector spaces \cite{Gan2}. $CC^{nu}_*(\mathcal{A})$ has an $k[\epsilon]/\epsilon^2\simeq C_*(S^1)$-action given by the non-unital Connes' operator $B^{nu}$, cf. (2.28). The $t$-complex for negative cyclic homology $CC^{S^1}(\mathcal{A})$ is defined as $CC^{nu}_*(\mathcal{A})[[t]]$, with $|t|=2$ and differential $b^{nu}+tB^{nu}$, whose cohomology we denote by $HH^{S^1}_*(\mathcal{A})$. \cite{Gan2} defined the (negative) cyclic open-closed map
\begin{equation}
OC^{S^1}: HH^{S^1}_*(\mathrm{Fuk}(X,\mathbf{k})_{\lambda})\rightarrow QH^{*+n}(X,\mathbf{k})[[t]].
\end{equation}
The following theorem compares the cyclic and $\mathbb{Z}/p$-equivariant open-closed maps, and is the key to applying characteristic zero techniques in computing quantum Steenrod operations. 
\begin{thm}(\cite[Theorem 1.6]{Che1})
1) For any cohomologically unital $A_{\infty}$-category $\mathcal{A}$, there exists a quasi-isomorphism $\Phi_p: CC^{S^1}(\mathcal{A})\oplus CC^{S^1}(\mathcal{A})\theta\simeq CC^{\mathbb{Z}/p}(\mathcal{A})$, where
$\theta$ is a formal variable of degree $1$.\par\indent
2) The following diagram is homotopy commutative:
\begin{equation}
\begin{tikzcd}[row sep=1.2cm, column sep=0.8cm]
CC^{\mathbb{Z}/p}(\mathrm{Fuk}(X,\mathbf{k})_{\lambda})\arrow[rrr,"{OC^{\mathbb{Z}/p}}"]& & &QH(X,\mathbf{k})[[t,\theta]]  \\
CC^{S^1}(\mathrm{Fuk}(X,\mathbf{k})_{\lambda})\oplus CC^{S^1}(\mathrm{Fuk}(X,\mathbf{k})_{\lambda})\theta \arrow[u,"{\Phi_p}"]\arrow[urrr,"{OC^{S^1}\oplus OC^{S^1}\theta}"]& & &
\end{tikzcd}
\end{equation}\qed
\end{thm}
The proof of Theorem 1.3 is rather technical, and was obtained in a prequel to this paper \cite{Che1}.
\begin{rmk}
From an abstract point of view, Theorem 1.3. 1) is an analogue of the classical $\mathbb{Z}/p$-Gysin sequence in topology, which relates the homology of the $S^1$-homotopy fixed points of an $S^1$-space with that of the induced $\mathbb{Z}/p\subset S^1$-homotopy fixed points.
\end{rmk}
Given Theorem 1.3, once reduces the computation of $OC^{\mathbb{Z}/p}$ to that of $OC^{S^1}$, in characteristic $p$. One strategy for computing the latter is to first do the computation in characteristic $0$ and reduce coefficients mod $p$. In other words, we take a subring $R\subset \overline{\mathbb{Q}}$ satisfying certain properties (cf. section 1.5) and then
\begin{enumerate}[label=\arabic*)]
    \item  first, we compute the ordinary open-closed map over $R$, which can often be obtained via a direct geometric argument. 
    \item then, over $R\subset \overline{\mathbb{Q}}$ (here being in characteristic $0$ is crucial), one can sometimes derive the cyclic open-closed from its non-equivariant term by its property that it intertwines the Getzler-Gauss-Manin connection and the quantum connection, cf. Theorem 1.5. The upshot is that in characteristic $0$, the constraints imposed by covariant constancy with respect to formal connections are very strong; in contrast, there are usually lots of covariantly constant endomorphisms in characteristic $p$. Finally, one finds a suitable ring map $\pi: R\rightarrow \mathbf{k}$ to some field $\mathbf{k}$ of characteristic $p$ and reduces coefficient along $\pi$. 
\end{enumerate}
This strategy is implicitly implemented in Section 5, where we compute certain quantum Steenrod operations on the intersection of two quadrics in $\mathbb{C}P^5$.
\subsection{Cyclic open-closed map and the $t$-connection}
In this subsection, we review the classical definition of $t$-connections in symplectic topology and properties of the cyclic open-closed map over $\overline{\mathbb{Q}}$. Recall that the negative cyclic homology $HH^{S^1}(\mathcal{A})$ of an $\mathcal{A}_{\infty}$-category $A$ is equipped with \emph{the Getzler-Gauss-Manin $t$-connection} given on the chain level\footnote{We remark that to be precise, in order for this formula to work for $\mathcal{A}=\mathrm{Fuk}(X)_{\lambda}$ (which in our convention is an uncurved $A_{\infty}$-category, i.e. $m_0=0$, cf. section 3.2), we need to pass to a quasi-equivalent $A_{\infty}$-category $\mathrm{Fuk}(X)_{\lambda}\rightarrow \mathrm{Fuk}(X)^+_{\lambda}$ which contains a strict unit $e^+$, and set $m_0=\lambda\cdot e^+$ (which makes $\mathrm{Fuk}(X)_{\lambda}^+$ a \emph{weakly curved} $A_{\infty}$-category). This can be achieved, for instance, by using Fukaya's homotopy unit construction, cf. \cite[section 10]{Gan1}. An alternative approach is to stick to $\mathrm{Fuk}(X)_{\lambda}$, but use the connection obtained by tensoring (1.10) with the one-dimensional connection $\mathcal{E}^{-\frac{\lambda}{t^2}}$, cf. (6.4).} by
\begin{equation}
\nabla^{GGM}_{d/dt}:=\frac{d}{dt}+\frac{\Gamma}{2t}-\frac{i\{\sum_{k\geq 0}(2-k)m_k\}}{2t^2},
\end{equation}
where $\Gamma$ is the length operator on Hochschild chains and $i$ is the contraction of a cyclic chain by a Hochschild cochain, cf. \cite[section 3.2]{Hug}, \cite[section 5]{OS}. We remark that in the literature, the $S^1$-equivariant variable is often named $u$ (and (1.10) is named the `$u$-connection'), whereas we use $t$. \par\indent
There is an analogous $t$-connection on $S^1$-equivariant quantum cohomology, called the \emph{quantum $t$-connection} given by
\begin{equation}
\nabla^{QH}_{d/dt}:=\frac{d}{dt}+\frac{\mu}{t}-\frac{c_1\star}{t^2},
\end{equation}
where $\mu$ is the grading operator on cohomology, i.e. $\mu(x)=\frac{k-n}{2}x$ for $x\in H^k(X)$.
\begin{thm}
On the level of cohomology,
\begin{equation}
OC^{S^1}\circ \nabla^{GGM}_{d/dt}=\nabla^{QH}_{d/dt}\circ OC^{S^1}.
\end{equation}
\end{thm}
This was originally conjectured by \cite{GPS}. It was first proved by \cite[Theorem 1.7]{Hug} in a slightly simplified technical setting. The proof that is closer to our setting is \cite[Theorem 6.3.5]{PS}, which however considers the so-called `$q$-connection'; we explain how Theorem 1.5 can be derived from the results in loc.cit. in Appendix A. \par\indent
We recall the following general fact about formal $t$-connections.
\begin{lemma} \cite[Lemma 2.13]{Hug}
Write $\{e_j\}$ for a generalized eigenbasis of $c_1\star$, and write $\nabla=\frac{d}{dt}+\frac{A_0}{t^2}+\frac{A_1}{t}$ in this basis. In particular, $A_0$ is in Jordan normal form with diagonal entries the eigenvalues of $c_1\star$. Then, there exists a basis $\{v_j\}$ of $QH^*(X)[[t]]$ such that 
\begin{enumerate}[label=\arabic*)]
    \item  $v_j|_{t=0}=e_j$.
    \item Write $\nabla=\frac{d}{dt}+\frac{1}{t^2}\sum_{i=0}^{\infty} \tilde{A}_i t^i$ in the basis $\{v_j\}$, then $\tilde{A}_0=A_0$ and all $\tilde{A}_i$'s respect the generalized eigen-decomposition of $A_0$.
    \item The diagonal blocks of $\tilde{A}_1$ and $A_1$ agree.
\end{enumerate}\qed
\end{lemma}
Lemma 1.6 implies that there exists a decomposition of the quantum $t$-connection
\begin{equation}
QH^*(X,\overline{\mathbb{Q}})[[t]]=\bigoplus_{\lambda\in\mathrm{spec}(c_1\star)}QH^*(X,\overline{\mathbb{Q}})[[t]]_{\lambda}    
\end{equation}
such that $QH^*(X,\overline{\mathbb{Q}})[[t]]_{\lambda}|_{t=0}=QH^*(X,\overline{\mathbb{Q}})_{\lambda}$.
We call the decomposition (1.13) the \emph{elementary Hukuhara-Levelt-Turittin (HLT) decomposition}. On the categorical side, as $\lambda$ ranges over the eigenvalues of $c_1\star$, we consider the big direct sum $\bigoplus_{\lambda\in \mathrm{spec}(c_1\star)} CC^{S^1}(\mathrm{Fuk}(X,\overline{\mathbb{Q}})_{\lambda})$, where each summand is equipped with its own Getzler-Gauss-Manin $t$-connection. A consequence of Theorem 1.5 is that $OC^{S^1}$ must also intertwine the two decompositions (cf. \cite[Cor 6.5.]{Hug}), i.e.
\begin{equation}
OC^{S^1}(HH^{S^1}_*(\mathrm{Fuk}(X,\overline{\mathbb{Q}})_{\lambda}))\subset QH^{*+n}(X,\overline{\mathbb{Q}})[[t]]_{\lambda}.
\end{equation}
We note that even though the decomposition (1.13) induced by Lemma 1.6 was stated over the complex numbers, it in fact holds over a finitely generated extension $R\subset\overline{\mathbb{Q}}$ of $\mathbb{Z}$ (which allows us to reduce mod $p$, for almost all primes $p$). The precise statement is the following.
\begin{lemma}(\cite[Corollary 2.4]{Che2})
Suppose $R$ is an integral domain such that $QH^*(X,R)$ has a decomposition $QH^*(X,R)=\bigoplus_{\lambda\in\mathrm{spec}(c_1\star)}QH^*(X,R)_{\lambda}$ into generalized eigenspaces of $c_1\star$, and that the difference between two distinct eigenvalues of $c_1\star$ is invertible in $R$. Then this generalized eigen-decomposition extends to a unique decomposition of the $t$-connection (1.13) over $R$. 
\end{lemma} \qed
\subsection{Change of base ring}
Given we are ultimately interested in characteristic $p$ enumerative invariants and their interaction with characteristic $0$ phenomena, the goal of this subsection is to explain how to `extend scalars' from characteristic zero to suitable positive characteristics. \par\indent
Let $R\subset K\subset \overline{\mathbb{Q}}$ be the ring of integers of some number field, up to inverting finitely many numbers, and assume that $2$ is inverted in $R$. We consider the following conditions:
\begin{enumerate}[label=A\arabic*)]
    \item  $c_1\star\in\mathrm{End}(QH^*(X,R))$ admits a generalized eigenspace decomposition. Meaning, there is a finite direct sum decomposition $QH^*(X,R)=\bigoplus_{\lambda\in \mathrm{spec}(c_1\star)} QH^*(X,R)_{\lambda}$, where each summand $QH^*(X,R)_{\lambda}$ is a finite free $R$-module with a basis with respect to which $c_1\star$ is in canonical Jordan form with $\lambda$'s on the diagonal. 
    \item The difference between any two distinct eigenvalues of $c_1\star$ is invertible in $R$.
   \item The total open-closed map 
\begin{equation}
OC: \bigoplus_{\lambda\in\mathrm{spec}(c_1)}HH_*(\mathrm{Fuk}(X,R)_{\lambda})\rightarrow QH^{*+n}(X,R)
\end{equation} 
is an isomorphism.
\end{enumerate}
We now say a few words about conditions A1)-A3). Condition A2) guarantees that the following (cf. section 2.9 \cite{Sh1} for the result over $\mathbb{C}$) holds over $R$: for $\lambda, \lambda'\in \mathrm{spec}(c_1\star)$, the image of $QH^*(X,R)_{\lambda'}$ under
\begin{equation}
CO: QH^*(X,R)\rightarrow HH^*(\mathrm{Fuk}(X,R)_{\lambda})\quad\textrm{is $0$ if $\lambda\neq \lambda'$}.
\end{equation}
Assuming conditions A1) and A2), Lemma 1.8 implies that there is a unique elementary HLT decomposition over $R$, i.e. a decomposition 
\begin{equation}
QH^*(X,R)[[t]]=\bigoplus_{\lambda\in\mathrm{spec}(c_1\star)} QH^*(X,R)[[t]]_{\lambda}
\end{equation}
of the connection $\nabla^{QH}_{\frac{d}{dt}}$, with ${QH^*(X,R)[[t]]_{\lambda}}|_{t=0}=QH^*(X,R)_{\lambda}$. \par\indent
Assumption A3) is related to Abouzaid's generation criterion, by work of \cite{Gan1}. In particular, results in loc.cit. implies that if the total open-closed map (1.15) hits the unit in $QH^*(X,R)$, then it is an isomorphism. Since $R\subset \overline{\mathbb{Q}}$, (1.14) continues to hold over $R$, i.e. we have
\begin{equation}
OC^{S^1}(HH^{S^1}_*(\mathrm{Fuk}(X,R)_{\lambda}))\subset QH^{*+n}(X,R)[[t]]_{\lambda}.
\end{equation}
By a standard spectral sequence argument with respect to the $t$-filtration, A3) implies that the cyclic open-closed map
\begin{equation}
OC^{S^1}: HH^{S^1}_*(\mathrm{Fuk}(X,R)_{\lambda}))\rightarrow QH^{*+n}(X,R)[[t]]_{\lambda}
\end{equation}
is an isomorphism for each $\lambda\in\mathrm{spec}(c_1\star)$.\par\indent
Now suppose we have a homomorphism of rings $\pi: R\rightarrow \mathbf{k}$, where $\mathbf{k}$ is a field $\mathbf{k}$ of odd characteristic $p$. For a fixed $R$, such homomorphism (usually not unique) exists for almost all $p$. \par\indent 
Consider the diagram of abelian groups induced by the change of base rings $\pi$: 
\begin{equation}
\begin{tikzcd}[row sep=1.2cm, column sep=0.8cm]
HH^{S^1}_*(\mathrm{Fuk}(X,\mathbf{k})_{\pi(\lambda}))\arrow[rrr,"OC^{S^1}/{\mathbf{k}}"] & & & QH^{*+n}(X,\mathbf{k})[[t]]\\
HH^{S^1}_*(\mathrm{Fuk}(X,R)_{\lambda})\arrow[u,"\pi_*"]\arrow[rrr,"OC^{S^1}/R"] & & &QH^{*+n}(X,R)[[t]]\arrow[u,"\pi_*"]
\end{tikzcd}
\end{equation}
The $c_1\star$-eigenvalues over $\mathbf{k}$ are the images under $\pi$ of the $c_1\star$-eigenvalues over $R$ and the generalized $c_1\star$-eigenspaces satisfy (note by condition A2) $\pi$ is injective on $\mathrm{spec}(c_1\star)$)
\begin{equation}
QH^*(X,\mathbf{k})_{\pi(\lambda)}=QH^*(X,R)_{\lambda}\otimes_R\mathbf{k}.
\end{equation} 
We can then base-change the decomposition $QH^*(X,R)[[t]]=\bigoplus_{\lambda\in \mathrm{spec}(c_1\star)} QH^*(X,R)[[t]]_{\lambda}$ along $\pi$ and define the \emph{$\pi$-induced HLT decomposition} as
\begin{equation}
QH^*(X,\mathbf{k})[[t]]=\bigoplus_{\lambda\in \mathrm{spec}(c_1\star)} QH^*(X,\mathbf{k})[[t]]_{\pi(\lambda)},
\end{equation}
where $QH^*(X,\mathbf{k})[[t]]_{\pi(\lambda)}:=QH^*(X,R)[[t]]_{\lambda}\otimes_{R[[t]]}\mathbf{k}[[t]]$.\par\indent 
One can consider the non-equivariant analogue (i.e. set $t=0$) of the commutative diagram (1.20), which makes it obvious that condition A3) implies that (the total) $OC/\mathbf{k}$ also hits the unit, and thus is an isomorphism by \cite{Gan1}. Thus the total cyclic open-closed map over $\mathbf{k}$ is also an isomorphism on cohomology. Then, commutativity of (1.20) implies that
\begin{equation}
(OC^{S^1}/{\mathbf{k}})(HH^{S^1}_*(\mathrm{Fuk}(X,\mathbf{k})_{\pi(\lambda)}))\subset QH^{*+n}(X,\mathbf{k})[[t]]_{\pi(\lambda)}.
\end{equation}
Thus Theorem 1.3 implies the following. 
\begin{cor}
\begin{equation}
(OC^{\mathbb{Z}/p}/{\mathbf{k}}(HH_*^{\mathbb{Z}/p}(\mathrm{Fuk}(X,\mathbf{k})_{\pi(\lambda)}))\subset QH^{*+n}(X,\mathbf{k})[[t,\theta]]_{\pi(\lambda)},
\end{equation}
where
\begin{equation}
QH^*(X,\mathbf{k})[[t,\theta]]_{\pi(\lambda)}:=QH^*(X,\mathbf{k})[[t]]_{\pi(\lambda)}\otimes_{\mathbf{k}[[t]]}\mathbf{k}[[t,\theta]].
\end{equation}\qed
\end{cor}
Combining Theorem 1.3 with Corollary 1.8 we recover the following compatibility between quantum Steenrod operations and the splitting (1.13), cf. \cite[Corollary 4.4]{Che2} for a more general version.
\begin{cor}
In the same setting as above, let $b\in QH^*(X,\mathbf{k})_{\pi(\lambda)}$ and $x\in QH^*(X,\mathbf{k})[[t,\theta]]_{\pi(\lambda')}$.
\begin{enumerate}[label=\arabic*)]
    \item If $\lambda\neq \lambda'$,  then $Q\Sigma_b(x)=0$.
    \item If $\lambda=\lambda'$, then $Q\Sigma_b(x)\in QH^*(X,\mathbf{k})[[t,\theta]]_{\pi(\lambda)}$.
\end{enumerate}
\end{cor}
\noindent\emph{Proof}. By assumption A3), we know that the total $CO$ and $OC$ (and hence $OC^{\mathbb{Z}/p}$) are isomorphisms on cohomology over $\mathbf{k}$. Thus, there exists unique $\sigma\in HH^{\mathbb{Z}/p}_*(\mathrm{Fuk}(X,\mathbf{k})_{\pi(\lambda')})$ such that $OC^{\mathbb{Z}/p}(\sigma)=x$. By Theorem 1.2 applied to the monotone Fukaya category associated to $\pi(\lambda')$,
\begin{equation}
Q\Sigma_b(x)=OC^{\mathbb{Z}/p}_{\pi(\lambda')}(\bigcap^{\mathbb{Z}/p}(CO_{\pi(\lambda')}(b), \sigma)).
\end{equation}
This always lies in $QH^*(X,\mathbf{k})[[t,\theta]]_{\pi(\lambda)}$ by Corollary 1.8. This proves part 2). If $\lambda\neq \lambda'$, then $CO_{\pi(\lambda')}(b)=0$, which proves part 1).  \qed\par\indent
\subsection*{Organization}
The organization of this paper is as follows. Section 2 will be devoted to concepts in homological algebra. We start by reviewing the definitions of $A_{\infty}$-categories, bimodules, and their Hochschild (co)homology. Then, we introduce the $p$-fold Hochschild complex and its $\mathbb{Z}/p$-equivariant complex, and define the main algebraic operation of this paper: the $\mathbb{Z}/p$-equivariant cap product. In section 3, we review the definition of the monotone Fukaya category and the construction of the quantum Steenrod operations. In section 4, we define the $\mathbb{Z}/p$-equivariant open-closed map and prove the main result Theorem 1.2. In section 5, we apply our main results to computing the quantum Steenrod operations in the example of intersection of two quadrics in $\mathbb{C}P^5$, conditional on Smith's categorical equivalence over a positive characteristic field $\mathbf{k}$. In section 6, we discuss a conjectural $B$-side mirror of quantum Steenrod operations. 

\subsection*{Acknowledgements}
First and foremost, I would like to thank my advisor Paul Seidel for suggesting this line of inquiry and for his patient and helpful guidance. I would also like to thank Denis Auroux, Sheel Ganatra, Jae Hee Lee, Ivan Smith, and Nicholas Wilkins for helpful discussions at various points. This research was partially supported by the Simons Foundation, through a Simons Investigator grant (256290).

\renewcommand{\theequation}{2.\arabic{equation}}
\setcounter{equation}{0}
\section{Hochschild (co)homology of $A_{\infty}$-categories}
In section 2.1, we review the definition of $A_{\infty}$-categories and bimodules. In section 2.2 and 2.3, we recall the various chain models for computing Hochschild homology and cohomology in the literature. In section 2.4, we introduce the $\mathbb{Z}/p$-equivariant Hochschild complex and in section 2.5 we define the $\mathbb{Z}/p$-equivariant cap product, which will be the categorical analogue of quantum Steenrod operations. In this section, we work over an arbitrary coefficient ring, until otherwise specified.
\subsection{$A_{\infty}$-categories and bimodules}
We review some basic definitions, following \cite[section 2]{Gan1}.
\begin{mydef}
An (non-unital) $A_{\infty}$-\emph{category} is the data of a set of objects $\mathrm{Ob}\mathcal{A}$, a graded $k$-vector space $\mathrm{hom}_{\mathcal{A}}(X_0,X_1)$ for any pair of objects, and compositions for each $d\geq 1$,
\begin{equation}
\mu_{\mathcal{A}}^d: \mathrm{hom}_{\mathcal{A}}(X_0,X_1)\otimes\cdots\otimes \mathrm{hom}_{\mathcal{A}}(X_{d-1},X_d)\rightarrow \mathrm{hom}_{\mathcal{A}}(X_0,X_d)[2-d]    
\end{equation}
satisfying the $A_{\infty}$-relations
\begin{equation}
\sum_{m,n}(-1)^{\maltese_{n+m+1}^d}\mu_{\mathcal{A}}^{d-m+1}(a_1,\cdots,a_n,\mu_{\mathcal{A}}^m(a_{n+1},\cdots,a_{n+m}),a_{n+m+1},\cdots,a_d)=0,
\end{equation}
where the sign is determined by $\maltese_n=|a_{n+m+1}|+\cdots+|a_d|-(d-n-m)$.
\end{mydef}
Let $\|a_i\|:=|a_i|-1$ denote the reduced degree, then $\mu$ has degree $1$ with respect to $\|\cdot\|$.

\begin{mydef}
Let $\mathcal{C},\mathcal{D}$ be $A_{\infty}$-categories. An $A_{\infty}$ $\mathcal{C}-\mathcal{D}$ \emph{bimodule} $\mathcal{M}$ is the data of:
\begin{itemize}
    \item  for each $V\in \mathcal{C}, V'\in\mathcal{D}$, a graded vector space $\mathcal{M}(V,V')$
    \item  for $r,s\geq 0$ and $V_0,\cdots,V_r\in \mathcal{C}$, $W_0,\cdots,W_s\in \mathcal{D}$, structure maps
\begin{equation}
\mu^{r|1|s}_{\mathcal{M}}: \mathcal{C}(V_r,\cdots,V_0)\otimes \mathcal{M}(V_0,W_0)\otimes \mathcal{D}(W_0,\cdots,W_s)\rightarrow \mathcal{M}(V_r,W_s)
\end{equation}
of degree $1-r-s$, where we denote
$$\mathcal{C}(V_r,\cdots,V_0)=\mathrm{hom}_{\mathcal{C}}(V_{r-1},V_r)\otimes \cdots\otimes\mathrm{hom}_{\mathcal{C}}(V_0,V_1).$$
\end{itemize}
They are required to satisfy the following equation for all $r,s\geq 0$:
$$ \sum (-1)^{\maltese_{-s}^{j+1}} \mu_{\mathcal{M}}^{r-i|1|s-j}(v_r,\cdots,v_{i+1},\mu_{\mathcal{M}}^{i|1|j}(v_i,\cdots,v_1,\mathbf{m},w_1,\cdots,w_j),w_{j+1},\cdots,w_s)$$
$$+ \sum (-1)^{\maltese_{-s}^{k}} \mu_{\mathcal{M}}^{r-i+1|1|s}(v_r,\cdots,v_{k+i+1},\mu_{\mathcal{C}}^i(v_{k+i},\cdots,v_{k+1}),v_k,\cdots,v_1,\mathbf{m},w_1,\cdots,w_s)$$
\begin{equation}
+ \sum (-1)^{\maltese_{-s}^{-(l+j+1)}} \mu_{\mathcal{M}}^{r|1|s-j+1}(v_r,\cdots,v_1,\mathbf{m},w_1,\cdots,w_l,\mu_{\mathcal{D}}^j(w_{l+1},\cdots,w_{l+j}),w_{l+j+1}endts,w_s)=0.
\end{equation}
The signs are given by
\begin{equation}
\maltese_{-s}^{-(j+1)}:=\sum_{i=j+1}^s\|w_i\|\quad,\quad \maltese_{-s}^k=\sum_{i=1}^s\|w_i\|+|\mathbf{m}|+\sum_{j=1}^k\|v_j\|.
\end{equation}
\end{mydef}

\begin{mydef}
Let $\mathcal{A}$ be an $A_{\infty}$-category. The \emph{diagonal bimodule} $\mathcal{A}_{\Delta}$ is the $\mathcal{A}-\mathcal{A}$ bimodule defined by $\mathcal{A}_{\Delta}(X,Y)=\mathrm{hom}_{\mathcal{A}}(Y,X)$, and
\begin{equation}
\mu_{\mathcal{A}_{\Delta}}^{r|1|s}(x_r,\cdots,x_1,\mathbf{a},y_1,\cdots,y_s)=(-1)^{\maltese_{-s}^{-1}+1}\mu^{r+1+s}_{\mathcal{A}}(x_r,\cdots,x_1,\mathbf{a},y_1,\cdots,y_s).
\end{equation}
\end{mydef}
If $\phi, \varphi$ are multilinear maps, then we define
\begin{equation}
\phi\circ \varphi(x_n,\cdots,x_1):=\sum (-1)^{|\varphi|\cdot\maltese_s} \phi(x_n,\cdots,x_{r+1},\varphi(x_r,\cdots,x_{s+1}),x_s,\cdots,x_1).
\end{equation}
That is, we sum over all possible ways to insert $\varphi$ into $\phi$, and as before $\maltese_s$ denotes the sum of degrees of $x_s,\cdots,x_1$, using reduced degree unless it is a bimodule entry. \par\indent
When one of $x_n,\cdots,x_1$ is a $\mathcal{C}-\mathcal{D}$ bimodule entry and $\varphi$ is the bimodule structure map, we use $\phi\circ (\mu_{\mathcal{C}},\varphi,\mu_{\mathcal{D}})$ to denote the above expression plus summing over first applying $\mu_{\mathcal{C}}$ or $\mu_{\mathcal{D}}$ and then applying $\phi$, with the appropriate Koszul signs.
For instance, the $A_{\infty}$-category structural identity and the $A_{\infty}$-bimodule structural identity can be expressed as
$\mu_{\mathcal{A}}\circ \mu_{\mathcal{A}}=0$ and $\mu_{\mathcal{M}}\circ (\mu_{\mathcal{C}},\mu_{\mathcal{M}},\mu_{\mathcal{D}})=0$, respectively. 

\begin{mydef}
A \emph{pre-morphism} of $\mathcal{C}-\mathcal{D}$ bimodules of degree $k$ $\mathcal{F}: \mathcal{M}\rightarrow \mathcal{M}'$ is the data of:
\begin{equation}
\mathcal{F}^{r|1|s}: \mathcal{C}(V_r,\cdots,V_0)\otimes \mathcal{M}(V_0,W_0)\otimes \mathcal{D}(W_0,\cdots,W_s)\rightarrow \mathcal{M'}(V_r,W_s)\,,r,s\geq 0
\end{equation}
of degree $k-r-s$.
\end{mydef}
The pre-morphisms between $\mathcal{M},\mathcal{M}'$ form a chain complex $\mathrm{hom}_{\mathcal{C}-\mathcal{D}}(\mathcal{M},\mathcal{M}')$ where the differential is given by
\begin{equation}
\delta(\mathcal{F})=\mu_{\mathcal{M}'}\circ \mathcal{F}-(-1)^{|\mathcal{F}|}\mathcal{F}\circ (\mu_{\mathcal{C}},\mu_{\mathcal{M}},\mu_{\mathcal{D}}).
\end{equation}
A pre-morphism $\mathcal{F}$ is \emph{closed} if $\delta(\mathcal{F})=0$ (also called a \emph{morphism} of bimodules). The operation $\circ$ in (2.7) defines a composition \begin{equation}
\mathrm{hom}_{\mathcal{C}-\mathcal{D}}(\mathcal{M},\mathcal{M}')\times \mathrm{hom}_{\mathcal{C}-\mathcal{D}}(\mathcal{M}',\mathcal{M}")\rightarrow \mathrm{hom}_{\mathcal{C}-\mathcal{D}}(\mathcal{M},\mathcal{M}")
\end{equation}
making $\mathcal{C}-\mathcal{D}$-bimodules into a dg category. \par\indent
Finally, we recall the notion of tensor product of bimodules.
\begin{mydef}
Given a $\mathcal{C}-\mathcal{D}$ bimodule $\mathcal{M}$ and a $\mathcal{D}-\mathcal{E}$ bimodule $\mathcal{N}$, their \emph{tensor product over} $\mathcal{D}$
\begin{equation}
\mathcal{M}\otimes_{\mathcal{D}}\mathcal{N}
\end{equation}
is a $\mathcal{C}-\mathcal{E}$ bimodule whose underlying graded vector space is $\mathcal{M}\otimes T\mathcal{D}[1]\otimes \mathcal{N}$, where $T$ denotes the tensor algebra (where we consider composable morphisms), and $\mathcal{D}[1]$ denotes using the reduced degree on $\mathrm{hom}_{\mathcal{D}}$; the differential $\mu_{\mathcal{M}\otimes_{\mathcal{D}}\mathcal{N}}^{0|1|0}$ is given by
$$\mu_{\mathcal{M}\otimes_{\mathcal{D}}\mathcal{N}}^{0|1|0}(\mathbf{m},d_1,\cdots,d_k,\mathbf{n})=$$
$$\sum (-1)^{\maltese_{-(k+1)}^{-(t+1)}}\mu_{\mathcal{M}}^{0|1|t}(\mathbf{m},d_1,\cdots,d_t)\otimes d_{t+1}\otimes\cdots\otimes d_k\otimes \mathbf{n}+\sum \mathbf{m}\otimes d_1\otimes\cdots\otimes d_{k-s}\otimes \mu_{\mathcal{N}}^{s|1|0}(d_{k-s+1},\cdots,d_k,\mathbf{n})$$
\begin{equation}
+\sum (-1)^{\maltese_{-(k+1)}^{-(j+i+1)}}\mathbf{m}\otimes d_1\otimes\cdots\otimes d_{k-s}\otimes \mu_{\mathcal{D}}^i(d_{j+1},\cdots,d_{j+i})\otimes d_{j+i+1}\otimes\cdots\otimes d_k\otimes \mathbf{n}; 
\end{equation}
for $r>0$ or $s>0$, structure maps
\begin{equation}
\mu_{\mathcal{M}\otimes_{\mathcal{D}}\mathcal{N}}^{r|1|0}(c_1,\cdots,c_r,\mathbf{m},d_1,\cdots,d_k,\mathbf{n})=\sum (-1)^{\maltese_{-(k+1)}^{-(t+1)}}\mu_{\mathcal{M}}^{r|1|t}(c_1,\cdots,c_r,\mathbf{m},d_1,\cdots,d_t)\otimes d_{t+1}\otimes\cdots d_k\otimes \mathbf{n},
\end{equation}
\begin{equation}
\mu_{\mathcal{M}\otimes_{\mathcal{D}}\mathcal{N}}^{0|1|s}(\mathbf{m},d_1,\cdots,d_k,\mathbf{n},e_1,\cdots,e_s)=\sum \mathbf{m}\otimes d_1\otimes\cdots d_{k-j}\otimes \mu_{\mathcal{N}}^{j|1|s}(d_{k-j+1},\cdots d_k,\mathbf{n},e_1,\cdots,e_s)
\end{equation}
and $\mu_{\mathcal{M}\otimes_{\mathcal{D}}\mathcal{N}}^{r|1|s}=0$ if $r>0,s>0$.
\end{mydef}

\subsection{Hochschild cohomology}
The \emph{Hochschild cochain complex} of an $A_{\infty}$-category $\mathcal{A}$ with coefficients in an $\mathcal{A}$-bimodule $\mathcal{M}$ is
\begin{equation}
CC^r(\mathcal{A},\mathcal{M}):=\prod_{X_0,\cdots,X_k}\mathrm{Hom}_{\mathbf{k}}(\mathcal{A}(X_0,\cdots,X_k),\mathcal{M}(X_0,X_k))[r],
\end{equation}
where we use reduced degree for hom spaces of $\mathcal{A}$ and usual degree for $\mathcal{M}$. The differential is given by
\begin{equation}
d_{CC^*}(\varphi)=\mu_{\mathcal{M}}\circ \varphi-(-1)^{\|\varphi\|}\varphi\circ \mu_{\mathcal{A}}
\end{equation}
When $\mathcal{M}=\mathcal{A}_{\Delta}$, we denote $CC^*(\mathcal{A}):=CC^*(\mathcal{A},\mathcal{A}_{\Delta})$. 
The \emph{two-pointed Hochschild cochain complex} is by definition 
\begin{equation}
{}_2CC^*(\mathcal{A},\mathcal{M}):=\mathrm{hom}_{\mathcal{A}-\mathcal{A}}(\mathcal{A}_{\Delta},\mathcal{M}),
\end{equation}
equipped with the differential (2.9).
There is a chain map
\begin{equation}
\Psi:CC^*(\mathcal{A},\mathcal{M})\rightarrow {}_2CC^*(\mathcal{A},\mathcal{M})
\end{equation}
defined by
\begin{equation}
\Psi(\phi)(x_1,\cdots,x_k,\mathbf{a},y_1,\cdots,y_l):=\sum (-1)^{\dag}\mu_{\mathcal{M}}(x_1,\cdots,x_k,\mathbf{a},y_1,\cdots,y_i,\phi(y_{i+1},\cdots,y_{l-s}),y_{l-s+1},\cdots,y_l),
\end{equation}
where $\dag=|\phi|\cdot(\sum_{j=l-s+1}^l\|y_j\|)$. Moreover, $\Psi$ is a quasi-isomorphism whenever $\mathcal{A}$ is cohomologically unital, cf.\cite[Prop 2.5.]{Gan1}.\par\indent
There is a product structure on the Hochschild cochain complex $CC^*(\mathcal{A},\mathcal{A})$ called the \emph{cup product}. It is given by
\begin{equation}
\phi\cup\psi(x_1,\cdots,x_k):=\sum(-1)^{\clubsuit}\mu^k(x_1,\cdots,x_i,\phi(x_{i+1},\cdots,x_{i+r}),\cdots,\psi(x_{j+1},\cdots,x_{j+l}),\cdots,x_k),
\end{equation}
where $\clubsuit=|\phi|\cdot(\sum_{s=j+l+1}^{k}\|x_s\|+\sum_{t=i+r+1}^j\|x_t\|)+|\psi|\cdot(\sum_{s=j+k+1}^k\|x_s\|)$. The cup product defines an algebra structure on $HH^*(\mathcal{A})$, which is unital when $\mathcal{A}$ is cohomologically unital. \par\indent
There is also an algebra structure on ${}_2HH^*(\mathcal{A})$ whose product is given by composition of bimodule morphisms. On cohomology, $\Psi$ of (2.18) is a map of algebras and is unital when $\mathcal{A}$ is cohomologically unital.

\subsection{Hochschild homology and cyclic homology}
Let $\mathcal{A}$ be an $A_{\infty}$-category and $\mathcal{M}$ be an $\mathcal{A}$-bimodule. The \emph{Hochschild chain complex} of $\mathcal{A}$ with coefficients in $\mathcal{M}$ is
\begin{equation}
CC_*(\mathcal{A},\mathcal{M}):=\bigoplus_{X_0,X_1,\cdots,X_k} \mathcal{M}(X_k,X_0)\otimes \mathcal{A}(X_0,\cdots,X_k),
\end{equation}
with grading given by $\deg(\mathbf{m}\otimes x_1\otimes\cdots\otimes x_k)=|\mathbf{m}|+\sum_{i=1}^k\|x_i\|$ and differential given by
$$ b(\mathbf{m}\otimes x_1\otimes\cdots\otimes x_k)=\sum (-1)^{\sharp_j^i}\mu_{\mathcal{M}}(x_{k-j+1},\cdots,x_k,\mathbf{m},x_1,\cdots,x_i)\otimes x_{i+1}\otimes\cdots\otimes x_{k-j}$$
\begin{equation}
+\sum (-1)^{\maltese_{-k}^{-(s+j+1)}}\mathbf{m}\otimes x_1\otimes \cdots\otimes x_s\otimes \mu_{\mathcal{A}}^j(x_{s+1},\cdots,x_{s+j})\otimes x_{s+j+1}\otimes\cdots\otimes x_k,
\end{equation}
where the new symbol $\sharp$ denotes
\begin{equation}
\sharp_j^i=\big(\sum_{s=k-j+1}^k\|x_s\|\big)\cdot\big(|\mathbf{m}|+\sum_{t=1}^{k-j}\|x_t\|\big)+\maltese_{-(k-j)}^{-(i+1)}.
\end{equation}
\begin{rmk}
In words, the sign is given by the parity of: if a cyclic permutation is involved, we sum the degrees of elements that get moved in front, and multiply with the sum of degrees of all other elements; if furthermore a structure map (e.g. $\mu_{\mathcal{A}}$ or $\mu_{\mathcal{M}}$) is applied, we sum the degrees of elements to the right of $\mu$, and add to the previous sum. Reduced degree is used except for a bimodule entry.
\end{rmk}
The cohomology of this complex is called \emph{Hochschild homology}, and the above chain model is also called the cyclic bar complex. We denote $CC_*(\mathcal{A}):=CC_*(\mathcal{A},\mathcal{A}_{\Delta})$. \par\indent
Next, we recall from \cite{Gan2} that the \emph{non-unital Hochschild complex} is defined as
\begin{equation}
CC^{nu}_*(\mathcal{A}):=CC_*(\mathcal{A})\oplus CC_*(\mathcal{A})[1].
\end{equation}
To define the differential, one considers the following two operations.
$$b'(x_d\otimes\cdots\otimes x_1):=\sum (-1)^{\maltese_1^s} x_d\otimes\cdots\otimes x_{s+j+1}\otimes \mu(x_{s+j},\cdots,x_{s+1})\otimes x_s\otimes\cdots\otimes x_1+$$
\begin{equation}
\sum (-1)^{\maltese_1^{d-j}}\mu(x_d,\cdots,x_{d-j+1})\otimes x_{d-j}\otimes\cdots\otimes x_1,
\end{equation}
and
\begin{equation}
d_{\wedge\vee}(x_d\otimes\cdots\otimes x_1):=(-1)^{\maltese_2^d+\|x_1\|\cdot\maltese_2^d+1}x_1\otimes x_d\otimes\cdots\otimes x_2+(-1)^{\maltese^{d-1}_1}x_d\otimes\cdots\otimes x_1.
\end{equation}
Then, the differential $b^{nu}$ is defined as
\begin{equation}
b^{nu}:=
\begin{pmatrix}
b&d_{\wedge\vee}\\
0& b'
\end{pmatrix}.
\end{equation}
When $\mathcal{A}$ is cohomologically unital, the natural inclusion $CC_*(\mathcal{A})\rightarrow CC^{nu}_*(\mathcal{A})$ is a quasi-isomorphism, see \cite[Prop. 2.2]{Gan1}. \par\indent
The crucial property of the non-unital Hochschild complex is that there is a chain level $S^1$-action on $CC^{nu}_*(\mathcal{A})$, for a not necessarily strictly unital $A_{\infty}$-category $\mathcal{A}$. The $S^1$-action takes the form of the (non-unital) \emph{Connes operator}, given by
\begin{equation}
B^{nu}(x_k\otimes \cdots\otimes x_1,y_l\otimes\cdots\otimes y_1):=\sum_i (-1)^{\maltese_1^i\maltese_{i+1}^k+\|x_k\|+\maltese_1^k+1}(0,x_i\otimes \cdots\otimes x_1\otimes x_k\otimes \cdots\otimes x_{i+1}).
\end{equation}
The Connes operator satisfies $(B^{nu})^2=0$ and $b^{nu}B^{nu}+B^{nu}b^{nu}=0$, and therefore can be viewed as a unital dg action of $C_*(S^1)$ on $CC^{nu}_*(\mathcal{A},\mathcal{A})$.\par\indent
Finally, one defines the \emph{negative cyclic chain complex} as the homotopy fixed point of this $S^1$-action, that is
\begin{equation}
CC^{S^1}_*(\mathcal{A}):=CC^{nu}_*(\mathcal{A})[t]],
\end{equation}
with differential given by $b^{nu}+tB^{nu}$. The homology of this complex is called the \emph{negative cyclic homology} of $\mathcal{A}$.

\subsection{The $N$-fold Hochschild complex}
Let $N$ be any positive integer. As seen in the previous subsection, the Hochschild homology of an $A_{\infty}$-category is equipped with an $S^1$-action. It is natural to ask the following question: consider the inclusion $\mathbb{Z}/N\subset S^1$ as $N$-th roots of unities, what is the induced $\mathbb{Z}/N$-action on Hochschild homology? To answer this question at the chain level, we consider the following variant of the Hochschild chain complex.
\begin{mydef}
The \emph{$N$-fold Hochschild chain complex} ${}_NCC_*(\mathcal{A})$ is defined as
\begin{equation}
{}_NCC_*(\mathcal{A}):=CC_*(\mathcal{A},\mathcal{A}_{\Delta}\otimes_{\mathcal{A}}\cdots\otimes_{\mathcal{A}}\mathcal{A}_{\Delta}),
\end{equation}
where the tensor product is $N$-fold. For $N=2$, this was considered by e.g. \cite[Definition 2.30]{Gan1}, \cite[section 3.7.2]{LT}.
\end{mydef}
More explicitly, the underlying graded vector space of ${}_NCC_*(\mathcal{A})$ is given by
\begin{equation}
(\mathcal{A}\otimes T\mathcal{A}[1]\otimes \mathcal{A}\otimes T\mathcal{A}[1]\otimes\cdots\otimes \mathcal{A}\otimes T\mathcal{A}[1])^{\mathrm{diag}},
\end{equation}
where $\mathrm{diag}$ means we consider cyclically composable sequence of morphisms. The differential can be described schematically by the following picture, where the boldface marked points correspond to distinguished (bimodule) entries.
\begin{figure}[H]
 \centering
 \includegraphics[width=0.7\textwidth]{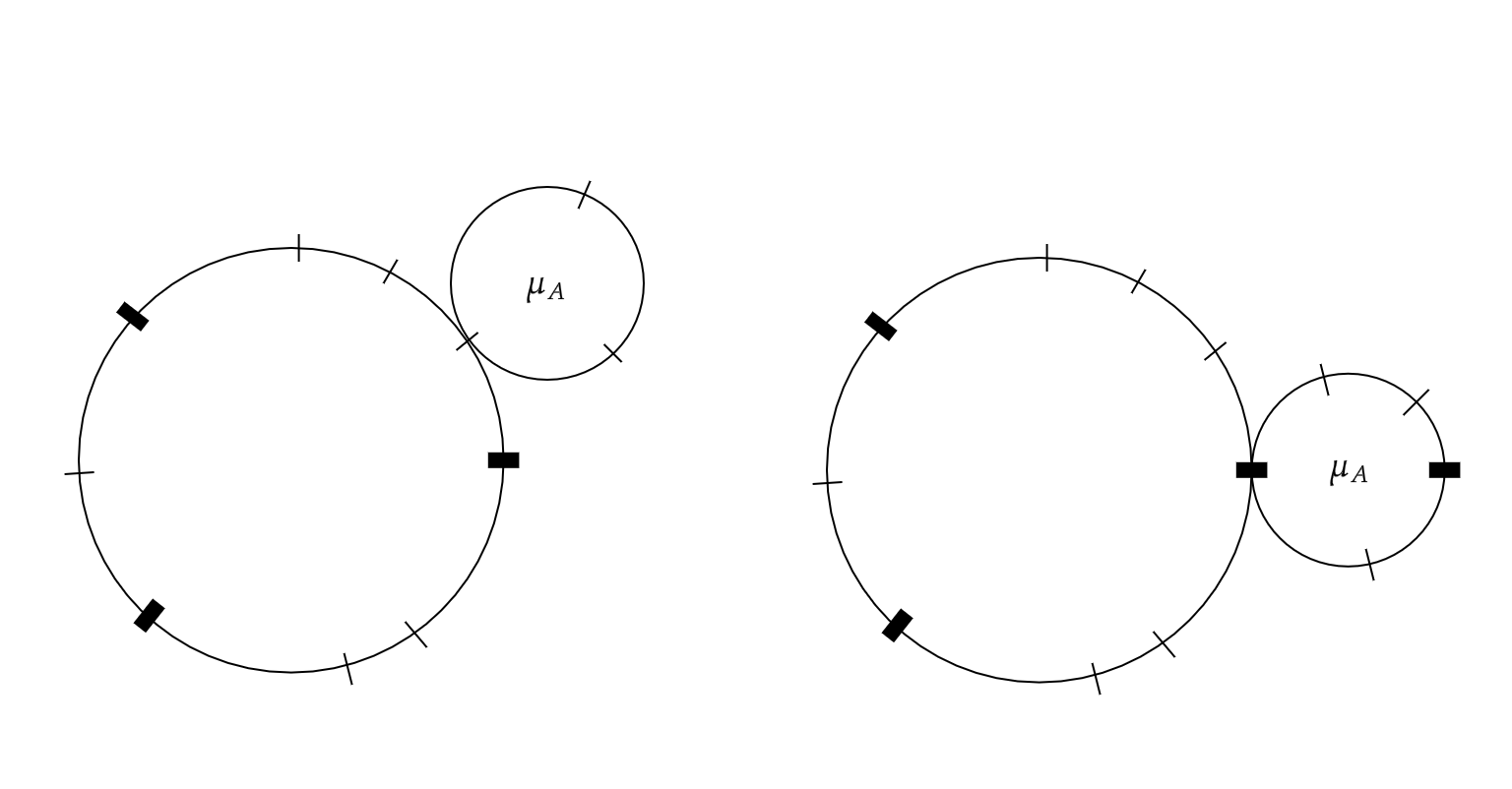}
 \caption{Two types of differentials in $d_{{}_NCC}$}
\end{figure}
When $\mathcal{A}$ is homologically unital, there is a canonical quasi-isomorphism of bimodules $\mathcal{A}_{\Delta}\otimes_{\mathcal{A}}\mathcal{A}_{\Delta}\simeq \mathcal{A}_{\Delta}$. Therefore, the $N$-fold Hochschild complex ${}_NCC_*(\mathcal{A})$ is quasi-isomorphic to ordinary Hochschild complex $CC_*(\mathcal{A})$. One can explicitly describe this quasi-isomorphism as follows. For $N\geq 2$, there is a map
\begin{equation}
\epsilon^0_{N-1,N}: {}_{N}CC_*(\mathcal{A})\rightarrow {}_{N-1}CC_*(\mathcal{A})
\end{equation}
defined by
$$
\epsilon^0_{N-1,N}(\mathbf{x}^1\otimes x_1^1\otimes \cdots\otimes x^1_{k_1}\otimes\mathbf{x}^2\otimes x_1^2\otimes \cdots\otimes x^2_{k_2}\otimes \cdots\otimes \mathbf{x}^N\otimes x_1^N\otimes \cdots\otimes x^N_{k_N}):=
$$
\begin{equation}
\sum (-1)^{(\sum_{j=i+1}^{k_N}\|x^N_{j}\|)(1+\sum_{l=1}^p(|\mathbf{x}^l|+\sum_{j=1}^{k_l}\|x^l_j\|))} \mu(x^N_{i+1},\cdots,x^N_{k_N},\mathbf{x}^1, x_1^1, \cdots, x^1_{k_1},\mathbf{x}^2,x^2_1,\cdots,x^2_j)\otimes x^2_{j+1}\otimes\cdots\otimes \mathbf{x}^N\otimes x_1^N\otimes \cdots\otimes x^N_{i},
\end{equation}
where the first term $\mu(\cdots)$ in (2.33), together with $\mathbf{x}^3,\cdots,\mathbf{x}^N$, becomes the new distinguished bimodule entries. When $\mathcal{A}$ is homologically unital, $\epsilon^0_{N,N-1}$ is a quasi-isomorphism for all $N\geq 2$, cf. the proof of \cite[Proposition 2.2]{Gan1}. Taking the composition of $\epsilon^0_{k,k-1}$ for $k=2,3,\cdots,N$, one obtains a quasi-isomorphism
\begin{equation}
\Phi^0_N:\epsilon^0_{N-1,N}\circ\cdots\circ \epsilon^0_{1,2}: {}_NCC_*(\mathcal{A})\rightarrow CC_*(\mathcal{A}).
\end{equation}
There is a chain level $\mathbb{Z}/N$-action on ${}_NCC_*(\mathcal{A})$, where the generator $\tau\in\mathbb{Z}/N$ acts by cyclically permuting the $N$ `blocks':
$$\tau: \mathbf{x}^1\otimes x_1^1\otimes \cdots\otimes x^1_{k_1}\otimes\mathbf{x}^2\otimes x_1^2\otimes \cdots\otimes x^2_{k_2}\otimes \cdots\otimes \mathbf{x}^N\otimes x_1^N\otimes \cdots\otimes x^N_{k_N}\mapsto $$
\begin{equation}
(-1)^{\dag}\mathbf{x}^N\otimes x_1^N\otimes \cdots\otimes x^N_{k_N}\otimes \mathbf{x}^1\otimes x_1^1\otimes \cdots\otimes x^1_{k_1}\otimes\cdots\otimes \mathbf{x}^{N-1}\otimes x_1^{N-1}\otimes \cdots\otimes x^{N-1}_{k_{N-1}},    
\end{equation}
where
\begin{equation}
\dag=\big(|\mathbf{x}^N|+\sum_{i=1}^{k_N} \|x^N_i\|\big)\cdot\big(\sum_{j=1}^{p-1} |\mathbf{x}^j|+\sum_{j=1}^{N-1}\sum_{i=1}^{k_j}\|x^j_{i}\|\big)    
\end{equation}
is the Koszul sign. One easily verifies that $\tau\circ d_{{}_NCC}=d_{{}_NCC}\circ\tau$.
\begin{mydef}
Consider the action 
\begin{equation}
\mathrm{hom}_{\mathcal{A}-\mathcal{A}}(\mathcal{M},\mathcal{M}')\otimes CC_*(\mathcal{A},\mathcal{M})\rightarrow CC_*(\mathcal{A},\mathcal{M}')    
\end{equation}
given by
\begin{equation}(\mathcal{F},\mathbf{m}\otimes x_n\otimes \cdots\otimes x_1)\mapsto (-1)^{\dag}\mathcal{F}(x_r,\cdots,x_1,\mathbf{m},x_n,\cdots,x_s)\otimes x_{s-1}\otimes\cdots\otimes x_{r+1},
\end{equation}
where $\dag=(\sum_{i=1}^r\|x_i\|)\cdot(\sum_{i=r+1}^n \|x_i\|+|\mathbf{m}|)+|\mathcal{F}|\maltese_{r+1}^{s-1}$. When $\mathcal{M}=\mathcal{M}'=\mathcal{A}_{\Delta}^{\otimes N}:=\mathcal{A}_{\Delta}\otimes_{\mathcal{A}}\cdots\otimes_{\mathcal{A}}\mathcal{A}_{\Delta}$ ($n$ times), precomposing (2.37) with the canonical chain map ${}_2CC^*(\mathcal{A})^{\otimes N}=\mathrm{End}_{\mathcal{A}-\mathcal{A}}(\mathcal{A}_{\Delta})^{\otimes N}\rightarrow \mathrm{End}_{\mathcal{A}-\mathcal{A}}(\mathcal{A}_{\Delta}^{\otimes N})$, we get a chain map
\begin{equation}
{}_N\prod: {}_2CC^*(\mathcal{A})^{\otimes N}\otimes {}_NCC_*(\mathcal{A})\rightarrow {}_NCC_*(\mathcal{A}),
\end{equation}
which can be schematically represented as in figure 2. This action descends to cohomology, which we call the \emph{$N$-fold $2$-pointed cap product}. Precompose with the ring map $\Psi:HH^*(\mathcal{A})\rightarrow {}_2HH^*(\mathcal{A})$, we obtain an action
\begin{equation}
{}_N\bigcap:HH^*(\mathcal{A})^{\otimes N}\otimes {}_NHH_*(\mathcal{A})\rightarrow {}_NHH_*(\mathcal{A}),
\end{equation}
which we call the \emph{$N$-fold cap product}.
\begin{figure}[H]
 \centering
 \includegraphics[width=0.7\textwidth]{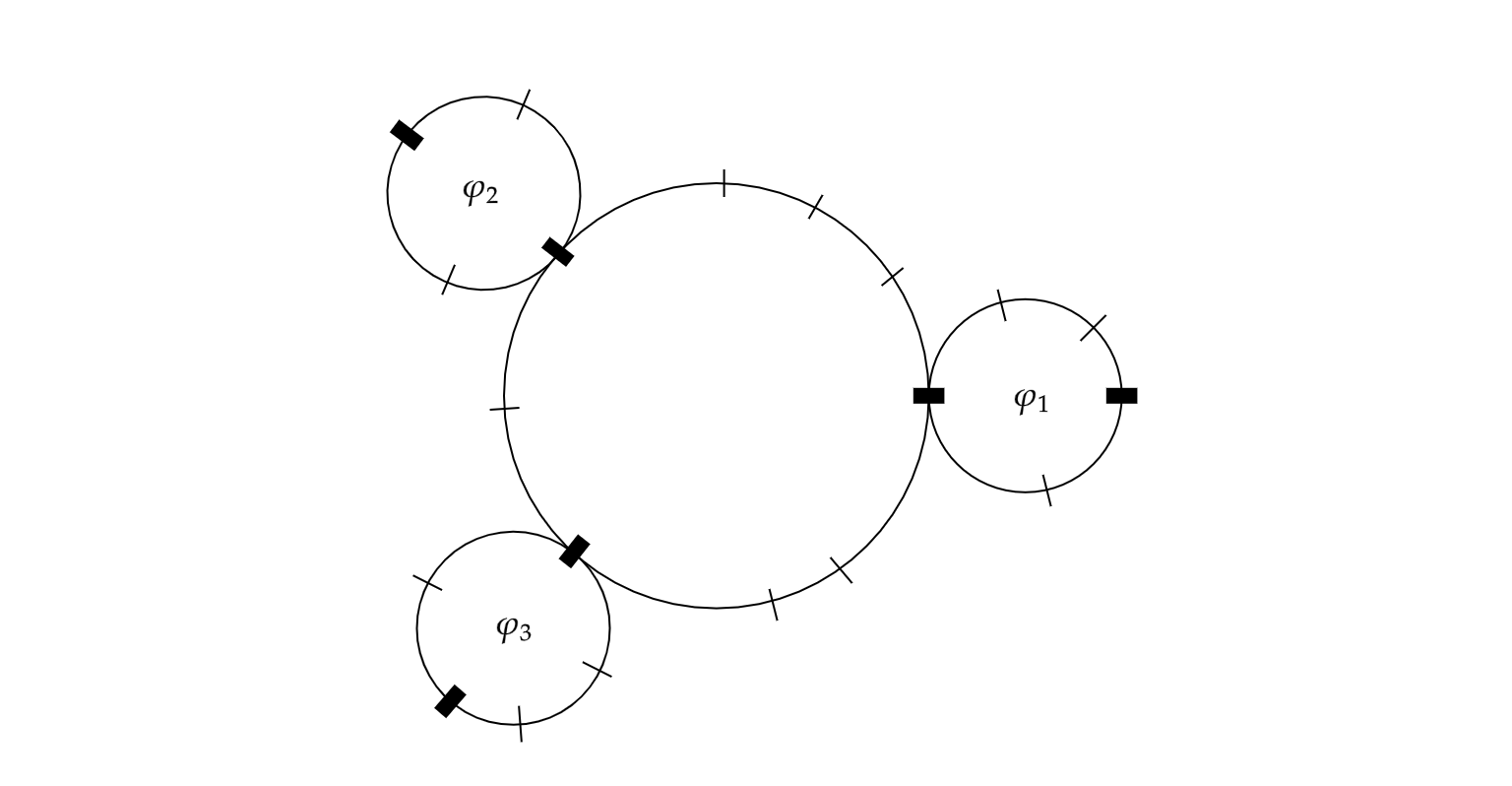}
 \caption{The action of $ {}_2CC^*(\mathcal{A},\mathcal{A})^{\otimes N}$ on ${}_NCC_*(\mathcal{A},\mathcal{A})$}
\end{figure}
\end{mydef}
It is straightforward from the definition that the $N$-fold $2$-pointed cap product (and hence the $N$-fold cap product) gives a unital graded associative algebra action on the chain level, i.e.
\begin{equation}
{}_N\prod\big((\phi_1\circ\phi_1',\cdots,\phi_N\circ\phi_N'),x\big)=(-1)^{\dagger}{}_N\prod\big((\phi_1,\cdots,\phi_N),{}_N\prod((\phi_1',\cdots,\phi_N'),x)\big)
\end{equation}
\begin{equation}
{}_N\prod\big((\mathrm{id},\mathrm{id},\cdots,\mathrm{id}),x\big)=x,  
\end{equation}
where $\dagger=\sum_{i=1}^N |\phi'_i|\cdot(|\phi_{i+1}|+\cdots+|\phi_N|)$ is the Koszul sign of reordering $(\phi_1,\cdots,\phi_N,\phi_1',\cdots,\phi_N')$ into $(\phi_1,\phi_1',\cdots,\phi_N,\phi_N')$. 
In addition, we make the following observation: let $\mathbb{Z}/N$ act on ${}_2CC^*(\mathcal{A},\mathcal{A})^{\otimes N}$ by cyclically permuting the tensor product (with appropriate Koszul signs) and take the diagonal $\mathbb{Z}/N$-action on ${}_2CC^*(\mathcal{A})^{\otimes N}\otimes {}_NCC_*(\mathcal{A})$, then ${}_N\prod$ is $\mathbb{Z}/N$-equivariant.

\subsection{The $\mathbb{Z}/p\mathbb{Z}$-equivariant cap product}
In this subsection, we let $p$ be an odd prime number and work over a field $\mathbf{k}$ of characteristic $p$. We first collect some standard facts about cochain complexes with a $\mathbb{Z}/p$-action and their homotopy fixed points. Then we define the $\mathbb{Z}/p$-equivariant cap product of $HH^*$ on $HH^{\mathbb{Z}/p}_*$.\par\indent
Let $X$ be a cohomologically graded complex with a $\mathbb{Z}/p$-actions. The associated \emph{negative $\mathbb{Z}/p$-equivariant complex}, or \emph{$\mathbb{Z}/p$ homotopy fixed point}, is defined as $X^{\mathbb{Z}/p}:=\mathrm{Rhom}_{k[\mathbb{Z}/p]}(k,X)$. An explicit chain complex computing $X^{\mathbb{Z}/p}$ is given by $X[[t,\theta]]$, where $|t|=2,|\theta|=1,\theta^2=0$ 
Let $\tau\in\mathbb{Z}/p$ be the standard generator, then the differential is given by
\begin{equation}
\begin{cases}
d_{eq}(x)=dx+(-1)^{|x|}(\tau x-x),\\
d_{eq}(x\theta)=dx\,\theta+(-1)^{|x|}(x+\tau x+\cdots+\tau^{N-1}x)t.     
\end{cases}  
\end{equation} 
and extended $t$-linearly.
\begin{lemma}
On $X^{\mathbb{Z}/p}$, the action by $\tau$ is homotopic to the identity.
\end{lemma}
\noindent\emph{Proof}. The homotopy is given by $h(x\,t^k)=0, h(x\,t^k\theta)=(-1)^{|x|}x\,t^k$.\qed\\\\
Suppose we have two $\mathbb{Z}/p$-complexes $X,Y$, then $X\otimes Y$ is also a $\mathbb{Z}/p$-complex equipped with the diagonal action. There is a canonical map
\begin{equation}
L_{XY}: X^{\mathbb{Z}/p}\otimes Y^{\mathbb{Z}/p}\rightarrow (X\otimes Y)^{\mathbb{Z}/p}
\end{equation}
given explicitly by (writing $X^{\mathbb{Z}/p}\otimes Y^{\mathbb{Z}/p}$ explicitly as $X\otimes Y[[t_1,t_2,\theta_1,\theta_2]]$)
\begin{equation}
\begin{cases}
x\otimes y\,t_1^{k_1}t_2^{k_2}\mapsto x\otimes y\,t^{k_1+k_2}\\
x\otimes y\,t_1^{k_1}\theta_1t_2^{k_2}\mapsto (-1)^{|y|}x\otimes \tau y\,t^{k_1+k_2}\theta\\
x\otimes y\,t_1^{k_1}t_2^{k_2}\theta_2\mapsto x\otimes y\,t^{k_1+k_2}\theta\\
x\otimes y\,t_1^{k_1}\theta_1 t_2^{k_2}\theta_2\mapsto (-1)^{|y|}\sum_{0\leq i<j\leq p-1}\tau^ix\otimes\tau^jy\,t^{k_1+k_2+1}
\end{cases}.
\end{equation}
Let $X={}_pCC(\mathcal{A})$ for some $A_{\infty}$-category $\mathcal{A}$. We shorthand ${}_pCC^{\mathbb{Z}/p}(\mathcal{A})$ as just $CC^{\mathbb{Z}/p}(\mathcal{A})$, and denote its homology by $HH^{\mathbb{Z}/p}(\mathcal{A})$. Combining (2.44) with (2.39), we obtain a chain map
\begin{equation}
\prod_{eq}: \big({}_2CC^*(\mathcal{A})^{\otimes p}\big)^{\mathbb{Z}/p}\otimes CC^{\mathbb{Z}/p}_*(\mathcal{A})\rightarrow CC^{\mathbb{Z}/p}_*(\mathcal{A}),
\end{equation}
which induces a cohomological level map
\begin{equation}
H^*_{\mathbb{Z}/p}\big({}_2CC^*(\mathcal{A})^{\otimes p}\big)\otimes HH^{\mathbb{Z}/p}_*(\mathcal{A})\rightarrow HH^{\mathbb{Z}/p}_*(\mathcal{A}),
\end{equation}
We view (2.47) as an action of the $\mathbf{k}$-algebra $H^*_{\mathbb{Z}/p}\big({}_2CC^*(\mathcal{A})^{\otimes p}\big)$ on $HH^{\mathbb{Z}/p}_*(\mathcal{A})$. Then this action is
\begin{itemize}
    \item graded multiplicative by (2.41) and the fact that the diagram
\begin{equation}
\begin{tikzcd}[row sep=1.2cm, column sep=0.8cm]
X^{\mathbb{Z}/p}\otimes Y^{\mathbb{Z}/p}\otimes Z^{\mathbb{Z}/p}\arrow[r,"\mathrm{id}\otimes L_{YZ}"]\arrow[d,"L_{XY}\otimes \mathrm{id}"] & X^{\mathbb{Z}/p}\otimes(Y\otimes Z)^{\mathbb{Z}/p}\arrow[d,"L_{X(Y\otimes Z)}"]\\
(X\otimes Y)^{\mathbb{Z}/p}\otimes Z^{\mathbb{Z}/p}\arrow[r,"L_{(X\otimes Y)Z}"]& (X\otimes Y\otimes Z)^{\mathbb{Z}/p}
\end{tikzcd}
\end{equation}
is homotopy commutative,
    \item unital by (2.42) and that the unit in $H^*_{\mathbb{Z}/p}\big({}_2HH^*(\mathcal{A})^{\otimes p}\big)$ has a chain representative $\mathrm{id}\otimes\mathrm{id}\otimes\cdots\otimes\mathrm{id}$.
\end{itemize}
\begin{lemma} \cite[Lemma 2.5]{SW}
Let $X$ be a cochain complex. Taking a cocycle $x\in X$ to its $p$-th power $x^{\otimes p}\in (X^{\otimes p})^{\mathbb{Z}/p}$ gives a well defined map
\begin{equation}
H^*(X)\rightarrow H^*_{\mathbb{Z}/p}(X^{\otimes p}),
\end{equation}
which becomes additive after multiplying by $t$.\qed
\end{lemma}
Observe that when $X$ is a dg algebra, the $p$-th power map (2.49) is also multiplicative; if furthermore $X$ is a unital, then the $p$-th power map is unital as well.\par\indent
Let $\mathcal{A}$ be a cohomologically unital $A_{\infty}$-category. Combining (2.47) and Lemma 2.10, one obtains a unital graded multiplicative action
\begin{equation}
\prod^{\mathbb{Z}/p}(-,-):=\prod_{eq}((-)^{\otimes p},-): {}_2HH^*(\mathcal{A})\otimes  HH^{\mathbb{Z}/p}_*(\mathcal{A})\rightarrow HH^{\mathbb{Z}/p}_*(\mathcal{A}),
\end{equation}
which is Frobenius $p$-linear and additive in the first entry after multiplying by $t$.
\begin{mydef}
$\prod^{\mathbb{Z}/p}$ is called the $\mathbb{Z}/p$-\emph{equivariant $2$-pointed cap product}.
\end{mydef}
Precomposing with the unital algebra map $\Psi: HH^*(\mathcal{A})\rightarrow {}_2HH^*(\mathcal{A})$, we obtain a unital multiplicative Frobenius linear action
\begin{equation}
\bigcap^{\mathbb{Z}/N}: HH^*(\mathcal{A})\otimes  HH^{\mathbb{Z}/p}_*(\mathcal{A})\rightarrow HH^{\mathbb{Z}/p}_*(\mathcal{A}),
\end{equation}
which is Frobenius $p$-linear and additive in the first entry after multiplying by $t$. 
\begin{mydef}
$\bigcap^{\mathbb{Z}/p}$ is called the $\mathbb{Z}/p$-\emph{equivariant cap product}.    
\end{mydef}
The above arguments show that $\bigcap^{\mathbb{Z}/p}$ satisfies the conditions C1)-C4) on page 2.

\renewcommand{\theequation}{3.\arabic{equation}}
\setcounter{equation}{0}
\section{The monotone Fukaya category}
In this section, we review the basic definitions of the monotone Fukaya category, which is the context we work in throughout the rest of the paper. Our exposition largely follows \cite{Sh1} and \cite{Oh}, to which we direct the readers for more details. In section 4.3, we review the definition of quantum Steenrod operations on a closed monotone symplectic manifold following \cite{SW}.
\subsection{Monotonicity}
\begin{mydef}
A symplectic manifold $(X,\omega)$ is \emph{monotone} if
$$\omega=2\tau c_1,$$
for some constant $\tau>0$. 
\end{mydef}
Furthermore, we only consider \emph{monotone} Lagrangians submanifolds $L$, meaning that
$$[\omega]=\lambda\mu: H_2(X,L)\rightarrow \mathbb{Z}$$
for some nonzero constant $\lambda$ (automatically $\lambda=\tau$, see \cite[Rmk 2.3.]{Oh}), where $\mu$ is the Maslov class, and that either $H^1(X)=0$ or the image of $\pi_1(L)$ in $\pi_1(X)$ is trivial. One important consequence of monotonicity is that for holomorphic disks/spheres of a fixed index, then there is a uniform bound on their energy, cf. \cite[Proposition 2.7]{Oh}. This allows us to appeal to Gromov compactness without using Novikov coefficients.

\subsection{The monotone Fukaya category}
In this subsection, we review the definition of the monotone Fukaya category associated with a closed monotone symplectic manifold $(X,\omega)$, over an arbitrary ring $R$. More precisely, to each $\lambda\in R$, one can associate a $k$-linear $\mathbb{Z}/2$-graded $A_{\infty}$-category $\mathrm{Fuk}(X,R)_{\lambda}$. \par\indent
Let $L$ be an oriented spin monotone Lagrangian submanifold $L\subset X$ equipped with a $R^*$-local system. Recall that monotonicity means that $\mu(L)=[\omega]$ when considered as classes in $H^2(X,L)$, where $\mu$ denotes the Maslov class. Orientability implies that the minimal Maslov number is $\geq 2$. By an abuse of notation,
we denote this datum simply by its underlying Lagrangian $L$. Let $\mathcal{J}$ denote the space of compatible almost complex structures and $\mathcal{H}:=C^{\infty}(X,\mathbb{R})$ the space of Hamiltonians. For each $L$, we fix $J_L\in \mathcal{J}$. For each pair $(L_0,L_1)$, we fix $J_t\in C^{\infty}([0,1],\mathcal{J})$ and $H_t\in C^{\infty}([0,1],\mathcal{H})$ such that $J_t=J_{L_t}$ when $t=0,1$.
If the $R^*$ local systems on both Lagrangians are trivial, the morphism space $CF^*(L_0,L_1)$ is the free $R$-module generated by time-$1$ Hamiltonian chords of $H_t$ from $L_0$ to $L_1$; in general, it is the direct sum of hom spaces between the fibers of the local systems at the startpoint and endpoint of the chord. \par\indent
Fix Lagrangians $L_0,L_1$. By standard transversality arguments, for generic almost complex structure $J_{L_0}, J_{L_1}$ and one parameter family $(H_t,J_t), t\in [0,1]$ such that $J_0=J_{L_0}, J_1=J_{L_1}$:
\begin{enumerate}[label=R\arabic*)]
    \item The moduli space $\mathcal{M}_1(L_0)$ of Maslov index $2$ $J$-holomorphic disks with one boundary marked point and boundary on $L_0$ is regular.
    \item The moduli space $\mathcal{M}_1(J_t)$ of pairs $(t,u)$, where $t\in [0,1]$ and $t$ is a Chern number $1$ $J_t$-holomorphic sphere with one marked point, is regular.
    \item For any time $1$ Hamiltonian chord $\gamma: [0,1]\rightarrow X$ starting on $L_0$ and ending on $L_1$, the map
\begin{equation}
(\gamma\circ t, \mathrm{ev}): \mathcal{M}_1(J_t)\rightarrow X\times X
\end{equation}
avoids the diagonal. In other words, all $J_t$ holomorphic spheres avoid $\gamma$.
    \item For Hamiltonian chords $x,y$ from $L_0$ to $L_1$, the moduli space $\mathcal{M}(x,y,J_t,H_t)$ of strips satisfying Floer's equation
\begin{equation}
(du-X_{H_t}\otimes dt)^{0,1}_{J_t}=0
\end{equation}
is regular.
\end{enumerate}
For each $L$,  since $L$ has minimal Maslov number $\geq 2$, the only possible nodal configuration in the Gromov compactification of $\mathcal{M}_1(L)$ is a $J_L$ holomorphic sphere of Chern number $1$ attached to a constant disk on $L$, and the moduli space of those has codimension $2$. In particular, $\mathcal{M}_1(L)$ has a well-defined pseudocycle fundamental class.  
If $L$ is equipped with the trivial local system, we define $w(L)\in \mathbb{Z}$ by
\begin{equation}
\mathrm{ev}_*[\mathcal{M}_1(L)]=w(L)[L]\in H_n(L,\mathbb{Z}).
\end{equation}
More generally, $\mathrm{ev}_*$ is weighted by the monodromy of the local system around the boundary of the disc, and $w(L)$ defines an element of $k$. \par\indent
The Floer differential $\mu^1:CF^*(L_0,L_1)\rightarrow CF^*(L_0,L_1)$ is defined as follows. Let $x_-,x_+\in CF^*(L_0,L_1)$,
then the coefficient of $x_-$ in $\mu^1(x_+)$ is the signed count (weighted by monodromy) of rigid elements of $\mathcal{M}(x_-,x_+)/\mathbb{R}$,
where $\mathcal{M}(x_-,x_+)$ is the moduli space of $u:\mathbb{R}\times [0,1]\rightarrow X$ such that
\begin{equation}
\begin{cases}
\partial_su+J_t(\partial_tu-X_{H_t})=0\\
u(s,0)\in L_0,\;\; u(s,1)\in L_1\\
\lim_{s\rightarrow \pm\infty}u(s,\cdot)=x_{\pm}.
\end{cases}
\end{equation}
By Gromov compactness and monotonicity, when $\mathcal{M}(x_-,x_+)$ is one dimensional (i.e. the Maslov index of $t$ is $1$), the space $\mathcal{M}(x_-,x_+)/\mathbb{R}$ is compact. When $\mathcal{M}(x_-,x_+)$ is $2$-dimensional, its Gromov compactification consists of broken strips $u_1,u_2$, each with Maslov index $1$, as well as a Maslov index $2$ disk bubbling off a Maslov index $0$ (hence constant in $s$) strip. For generic $J_t$, sphere bubbling cannot occur by regularity assumption (R3). Therefore, we have
\begin{equation}
\mu^1(\mu^1(x))=(w(L_0)-w(L_1))x.
\end{equation}
Hence, if $w(L_0)=w(L_1)$, then $(\mu^1)^2=0$.  \par\indent
We define the objects of $\mathrm{Fuk}(X,R)_{\lambda}$ to be oriented spin monotone Lagrangian submanifolds $L\subset X$ equipped with a $R^*$-local system, such that $w(L)=\lambda$. The morphism chain complexes
are defined as $(CF^*(L_0,L_1),\mu^1)$. \par\indent
We now describe the higher $A_{\infty}$ operations in $\mathrm{Fuk}(X,R)_{\lambda}$. Let $S$ be a surface with boundary and interior marked points (the boundary marked points are thought of as punctures).
Given a Lagrangian labeling $\mathbf{L}$ of the boundary components of $S$, a \emph{labeled Floer datum} for $S$ consists of the following data:
\begin{itemize}
    \item for each boundary marked point $\zeta$, a strip-like end $\epsilon_{\zeta}: \mathbb{R}^{\pm}\times [0,1]\rightarrow S$ at $\zeta$ (the strip-like end being positive or negative depending on whether $\zeta$ is an input or an output);
    \item a choice of $K\in \Omega^1(S,\mathcal{H})$ and $J\in C^{\infty}(S,\mathcal{J})$ such that $K(\xi)|_{L_C}=0$ for all $\xi\in TC$, where $C$ is a boundary component and $L_C$ is the corresponding Lagrangian label. Moreover, $K,J$ are compatible with strip-like ends in the sense that
\begin{equation}
\epsilon_{\zeta}^*K=H_{\zeta}(t)dt,\quad J(\epsilon_{\zeta}(s,t))=J_{\zeta}(t),
\end{equation}
where $H_{\zeta}, J_{\zeta}$ are the chosen Hamiltonian and almost complex structure for the pair of Lagrangians meeting at $\zeta$. We also require $J=J_L$ when restricted to a boundary component labeled $L$. The pair $(K,J)$ is called a \emph{perturbation datum}.
\end{itemize} 
The higher $A_{\infty}$-operations of $\mathrm{Fuk}(X,R)_{\lambda}$ are governed by the Delign-Mumford moduli space of disks with boundary marked points. Let $\mathcal{R}^{d+1}$ be the moduli space of disks with one boundary output and $d$ boundary inputs. It admits a compactification to a manifold with corners $\overline{\mathcal{R}}^{d+1}$ given by
\begin{equation}
\overline{\mathcal{R}}^{d+1}=\coprod_T \mathcal{R}^T,
\end{equation}
where $T$ ranges over all planar stable $d$-leafed trees and $\mathcal{R}^T:=\prod_{v\in \mathrm{Ve}(T)}\mathcal{R}^{|v|}$.\par\indent
We make a \emph{consistent choice of labeled Floer data} for $\mathcal{R}^{d+1}, d\geq 2$, meaning it is compatible with the product of Floer data of lower dimensional $\mathcal{R}^{d'}$'s near a boundary stratum, see \cite[section (9g),(9i)]{Sei1}.\par\indent
The higher operations $\mu^d, d\geq 2$ are then defined by counting isolated elements of the parametrized moduli space $\mathcal{M}(y_1,\cdots,y_d;y_-)$, which is the space of $(r,u)$, $r\in \mathcal{R}^{d+1}, u:\mathcal{S}_r\rightarrow X$ satisfying
\begin{equation}
(du-Y_K)^{0,1}_J=0,
\end{equation}
with appropriate Lagrangian boundary and asymptotic conditions, where $Y_K$ is the one-form on $S$ with value the Hamiltonian vector field associated to $K$. For a generic choice of Floer data, this moduli space is regular (\cite[section (9k)]{Sei1}).

\subsection{Quantum Steenrod operations}
In this subsection, we review the definition of quantum Steenrod operations given in \cite{SW} using the Morse chain model of quantum cohomology, and study some of its properties. \par\indent
\textbf{Quantum Cohomology}. Fix a ground ring $R$, a Morse function $f$ and metric $g$ such that the associated gradient flow is Morse-Smale. The Morse complex $CM^*(f)$ is generated by critical points of $f$ over $R$ with the $\mathbb{Z}$-grading given by the Morse index. The differential is given by counting gradient flow lines between two critical points whose indices differ by $1$. We use the cohomological convention, i.e. $|x|=\dim W^s(x), 2n-|x|=\dim W^u(x)$, where $W^s,W^u$ are the stable and unstable manifold of $x$, respectively. \par\indent
The chain level quantum cup product can be defined as follows, cf. \cite[section 3]{SW}. Given a homology class $A\in H_2(X,\mathbb{Z})$, inputs $x_0,x_1\in CM^*(f)$ and output $x_{\infty}\in CM^*(f)$,  let $\mathcal{M}_A(C,x_0,x_1,x_{\infty})$ be the moduli space of maps $u:C=\mathbb{C}P^1\rightarrow X$ with three marked points $z_0,z_1,z_{\infty}\in C$ such that $t$ lies in class $A$, and $z_0,z_1,z_{\infty}$ are constrained on $W^u(x_0), W^u(x_1)$ and $W^s(x_{\infty})$, respectively. Equivalently, it is the moduli space of $J$-holomorphic $u: \mathbb{C}P^1\rightarrow X$ together with gradient half-flowlines \begin{equation}
y_0,y_1:(-\infty,0]\rightarrow X,\quad y_{\infty}:[0,\infty)\rightarrow X    
\end{equation}
satisfying
$$y_k'=\nabla f(y_k), y_k(0)=u(z_k), \lim_{s\rightarrow -\infty}y_k(s)=x_k,$$
\begin{equation}
y_{\infty}'=\nabla f(y_{\infty}), y_{\infty}(0)=u(z_{\infty}), \lim_{s\rightarrow \infty}y_{\infty}(s)=x_{\infty}.
\end{equation}
For a generic choice of $J$, the moduli space $\mathcal{M}_A(C,x_0,x_1,x_{\infty})$ is regular of dimension $2c_1(A)+|x_{\infty}|-|x_0|-|x_1|$.
We define the coefficient of $x_{\infty}$ in the quantum cup product of $x_0$ with $x_1$ by counting isolated elements of $\mathcal{M}_A(C,x_0,x_1,x_{\infty})$, over all $A\in H_2(X,\mathbb{Z})$. For a fixed $c_1(A)$, by monotonicity and Gromov compactness, there are only finitely many homology class $A$ admitting $J$-holomorphic curves. Thus there is a well defined map $\star: CM^*(f)\otimes CM^*(f)\rightarrow CM^*(f)$, which is easily seen to descend to cohomology, and will be called the \emph{quantum cup product}. \par\indent
\textbf{Quantum Steenrod operations}. Fix a field $\mathbf{k}$ of odd characteristic. To define the quantum Steenrod operations, following \cite[section 4a]{SW}, we consider a moduli problem with fixed domain but parametrized (equivariant) Floer data. The relevant Floer data will be parametrized by
\begin{equation}
S^{\infty}:=\{w=(w_0,w_1,\cdots)\in\mathbb{C}^{\infty}: w_k=0\;\textrm{for}\;k\gg 0, \|w\|^2=1\}.
\end{equation}
For a prime $p$, there is a $\mathbb{Z}/p$-action on $S^{\infty}$ where the standard generator $\tau\in\mathbb{Z}/p$ acts by
\begin{equation}
\tau(w_0,w_1,\cdots)=(e^{2\pi i/p} w_0,e^{2\pi i/p}w_1,\cdots).    
\end{equation}
Consider the cells
\begin{equation}
\Delta_{2k}=\{w\in S^{\infty}: w_k\geq 0, w_{k+1}=w_{k+2}=\cdots=0\},
\end{equation}
\begin{equation}
\Delta_{2k+1}=\{w\in S^{\infty}: e^{i\theta}w_k\geq 0\;\textrm{for some}\;\theta\in[0,2\pi/p], w_{k+1}=w_{k+2}=\cdots=0\}.
\end{equation}
We identify the tangent space of $\Delta_{2k}$ at the point $w_{k}=1$ (and hence all other coordinates are zero) with $\mathbb{C}^k$, via the projection onto the first $k$ coordinates; we use the induced complex orientation on $\Delta_{2k}$. The tangent space of $\Delta_{2k+1}$ at the same point is canonically identified with $\mathbb{C}^k\times i\mathbb{R}$, and we use the complex orientation on the first factor and the positive vertical orientation on the second factor. With these chosen orientations, one has
\begin{equation}
\partial\Delta_{2k}=\Delta_{2k-1}+\tau\Delta_{2k-1}+\cdots+\tau^{p-1}\Delta_{2k-1},
\end{equation}
\begin{equation}
\partial\Delta_{2k+1}=\tau\Delta_{2k}-\Delta_{2k}.
\end{equation}
Let $C=\mathbb{C}P^1$, equipped with the $\mathbb{Z}/p$ action given by rotating by $e^{2\pi i/p}$. Denote the action of the generator by $\sigma$. Fix $p+2$ points $z_0=0, z_k=e^{2\pi ik/p}, k=1,2,\cdots,p, z_{\infty}=\infty$ on $C$. Choose perturbation data $K_w$ on $C$ parametrized by $w\in S^{\infty}$ such that
\begin{equation}
\sigma^*K_w=K_{\tau(w)}.
\end{equation}
Let $\mathcal{M}_A(\Delta_i\times C,x_0,x_1,\cdots,x_p,x_{\infty})$ denote the moduli space of pairs $(w,u)$, where $w\in \Delta_i$ and $u: C\rightarrow X$ in class $A$ satisfying
\begin{equation}
(du-Y_{K_w})^{0,1}_{J}=0
\end{equation}
and incidence conditions to $W^u(x_0),\cdots,W^u(x_p),W^s(x_{\infty})$ as before. We have
\begin{equation}
\dim\mathcal{M}_A(\Delta_i\times C,x_0,x_1,\cdots,x_p,x_{\infty})=i+2c_1(A)+|x_{\infty}|-|x_0|-|x_1|-\cdots-|x_{p}|.
\end{equation}
Moreover, because of the condition imposed by (3.17), we have an identification
\begin{equation}
\mathcal{M}_A(\tau^j(\Delta_i)\times C,x_0,x_1,\cdots,x_p,x_{\infty})\cong \mathcal{M}_A(\Delta_i\times C,x_{p-j+1},\cdots,x_p,x_1,\cdots,x_{p-j},x_{\infty})
\end{equation}
given by $(w,u)\mapsto (\tau^{-j}(w),u\circ \sigma^{-j})$. This defines maps, for $i\geq 0$,
\begin{equation}
\Sigma_A^i:=\Sigma_A(\Delta_i,\cdots):CM^*(f)\otimes CM^*(f)^{\otimes p}\rightarrow CM^{*-i-2c_1(A)}(f).
\end{equation}
Fix a Morse cocycle $b\in CM^*(f)$, and let $CM(f)[[t,\theta]]$ be the $\mathbb{Z}/p$-equivariant complex of $CM^*(f)$ (where $\mathbb{Z}/p$ acts trivially). One can combine the maps in (3.21) into a chain map
\begin{equation}
\Sigma_{A,b}: CM^*(f)\rightarrow (CM^*(f)[[t,\theta]])^{*+p|b|-2c_1(A)}
\end{equation}
given by
\begin{equation}
x\mapsto (-1)^{|b||x|}\sum_k\big(\Sigma_A(\Delta_{2k},x,b,\cdots,b)+(-1)^{|b|+|x|}\Sigma_A(\Delta_{2k+1},x,b,\cdots,b)\theta\Big)t^k.
\end{equation}
Up to homotopy, (3.23) only depends on the cohomology class $[b]\in QH^*(X,\mathbf{k})$, cf. \cite[Lemma 4.4.]{SW}. Finally, summing over $A$ (which is well defined by monotonicity) and extending $(t,\theta)$-linearly, we obtain a chain map 
\begin{equation}
\Sigma_b: CM^*(f)[[t,\theta]]\rightarrow CM^*(f)[[t,\theta]]    
\end{equation} 
of degree $p|b|$. We denote the cohomology level map of $\Sigma_b$ as $Q\Sigma_b$. Then, $Q\Sigma$ defines a Frobenius $p$-linear action of $QH^*(X,\mathbf{k})$ on $QH^*_{\mathbb{Z}/p}(X,\mathbf{k})$ such that
\begin{equation}
Q\Sigma_b\circ Q\Sigma_{b'}=(-1)^{\frac{p(p-1)}{2}|b||b'|} Q\Sigma_{b\star b'},
\end{equation}
cf. \cite[Proposition 4.8]{SW}. (3.25) is sometimes called the \emph{Quantum Cartan relation}.

\renewcommand{\theequation}{4.\arabic{equation}}
\setcounter{equation}{0}
\section{The $\mathbb{Z}/p\mathbb{Z}$-equivariant open-closed map}
In this section, we work over a field $\mathbf{k}$ of odd characteristic $p$. In subsection 4.1 and 4.2, we define the $p$-fold open-closed map and its equivariant enhancement. In subsection 4.3, we review the definition of two versions of the closed-open map following \cite{Gan1}. In subsection 4.4, we prove our main result that the $\mathbb{Z}/p$-equivariant open-closed map intertwines the $\mathbb{Z}/p$-equivariant cap product with quantum Steenrod operations.

\subsection{The $p$-fold open-closed map}
Fix $\lambda\in\mathbf{k}$ and let $\mathcal{F}=\mathrm{Fuk}(X)_{\lambda}$ be the monotone Fukaya category associated to $\lambda$. In this subsection, we construct a chain map
\begin{equation}
{}_pOC:{}_pCC_*(\mathcal{F})\rightarrow CM^*(f)
\end{equation}
of degree $n$. We remark that the case $p=2$ is constructed in \cite[section 5.6]{Gan2}.\par\indent
For a $p$-tuple of non-negative integers $k_1,\cdots,k_p$. Let
\begin{equation}
\mathcal{R}^1_{k_1,\cdots,k_p}
\end{equation}
be the moduli space of disks with one interior output marked point $y_{out}$ and $k_1+\cdots+k_p+p$ boundary input marked points $z^1,z^1_1,\cdots,z^1_{k_1},z^2,z^2_1,\cdots,z^2_{k_2},\cdots,z^p,z^p_1,\cdots,z^p_{k_p}$ in counterclockwise order such that up to automorphisms of the disk, $y_{out},z^1,z^2,\cdots,z^p$ lie at $0,\zeta,\zeta^2,\cdots,\zeta^p$, where $\zeta=e^{2\pi i/p}$. A representative of an element in $\mathcal{R}^1_{k_1,\cdots,k_p}$ is called \emph{standard} when $y_{out},z^1,\cdots,z^p$ satisfy the previous constraints. $z^1,\cdots,z^p$ are called the \emph{distinguished inputs}. 
We fix the orientation on (4.2) induced by $-dz^1_1\wedge\cdots\wedge dz^1_{k_1}\wedge\cdots\wedge dz^p_1\wedge\cdots\wedge dz^p_{k_p}$ on a standard representative. 
\begin{figure}[H]
 \centering
 \includegraphics[width=1.0\textwidth]{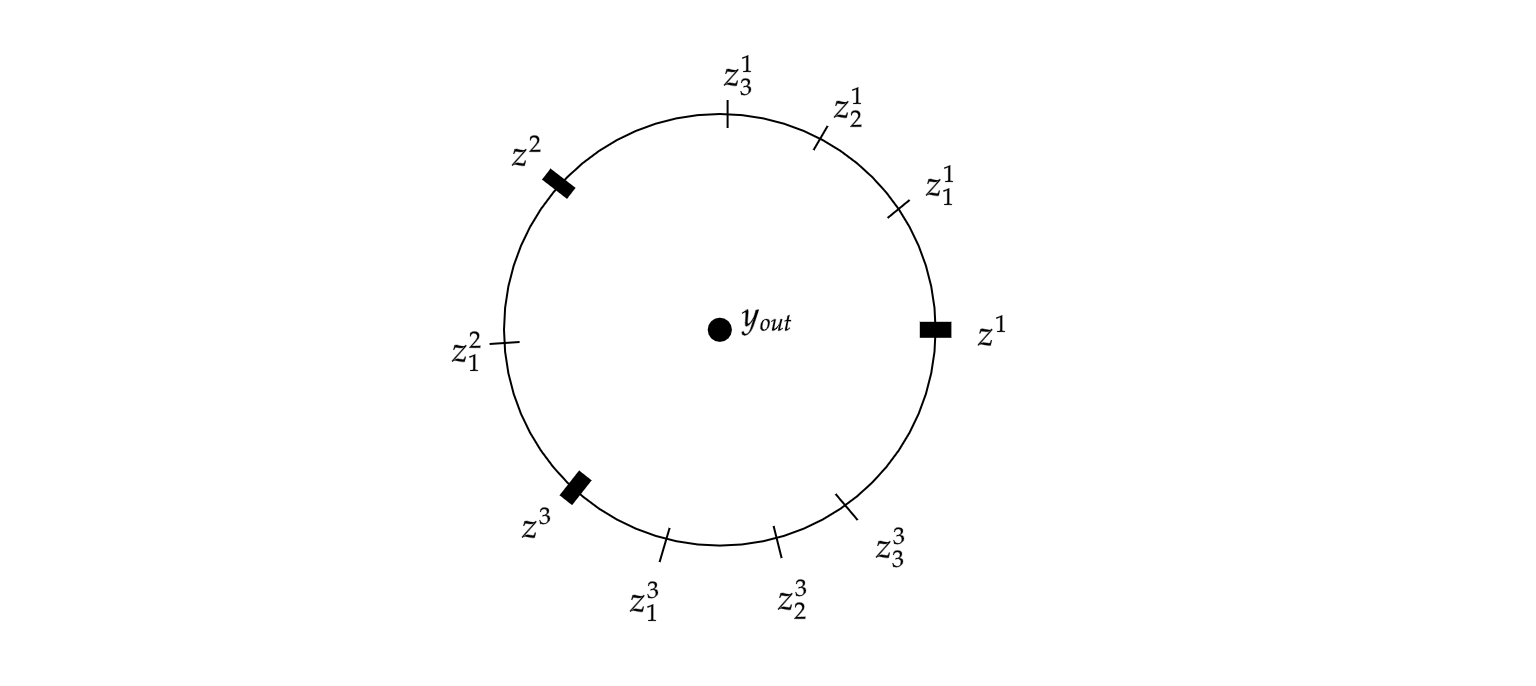}
 \caption{An element of $\mathcal{R}^1_{3,1,3}$}
\end{figure}
$\mathcal{R}^1_{k_1,\cdots,k_p}$ has a Deligne-Mumford compactification given by the following. Let $\mathcal{T}_p$ denote equivalence classes of rooted planar trees with $k_1+\cdots+k_p+p-1$ for some $k_1,\cdots,k_p\geq 0$ leaves, with one internal vertex
marked as \emph{main} and $p-1$ leaves marked as \emph{distinguished}. Those $p-1$ leaves, together with the root, are called the $p$ distinguished semi-infinite edges. Moreover, we require that no two distinguished semi-infinite edges
are adjacent to the same non-main vertex and each non-main vertex is adjacent to $\geq 3$ edges. In particular, one has $|main|\geq p$. 
Then,
\begin{equation}
\overline{\mathcal{R}}^1_{k_1,\cdots,k_p}=\bigsqcup_{T\in\mathcal{T}_p} \mathcal{R}^T, 
\end{equation}
where
\begin{equation}
\mathcal{R}_T:=\mathcal{R}^1_{k_1^T,\cdots,k_p^T}\times \prod_{v\in \mathcal{T}\backslash\{main\}} \mathcal{R}^{|v|}.
\end{equation}
The $k_i^T$'s are defined as follows: among the edges adjacent to the main vertex, there are $p$ special ones that are contained in the respective paths from the main vertex to the $p$ distinguished semi-infinite edges. $k_i^T$ is defined as the number of edges in between the $i$-th and $i+1$-th special edge in counter-clockwise order.\par\indent
In particular, the codimension $1$ boundary of $\overline{\mathcal{R}}^1_{k_1,\cdots,k_p}$ is covered by images of
\begin{equation}
\mathcal{R}^{k'_i+1}\times \mathcal{R}^1_{k_1,\cdots,k_i-k_i'+1,\cdots,k_p},\quad 1\leq i\leq p,
\end{equation}
\begin{equation}
\mathcal{R}^{k_{i-1}'+k_{i}'+2}\times \mathcal{R}^1_{k_1,\cdots,k_{i-1}-k_{i-1}',k_{i}-k_{i}',\cdots,k_p},\quad 1\leq i\leq p.
\end{equation}
under the natural inclusions.
In (4.5), the output of the first disk can be glued to any one of the inputs of the main component labeled by $z^i_1,\cdots,z^i_{k_i-k_i'+1}$, whereas in (4.6), the output of the first disk is glued to the distinguished input $z^i$, and the $k_i+1$-th input of the first disk will become the new distinguished $i$-th input after gluing. \par\indent
For each Lagrangian labeling of the universal families $\mathcal{S}^1_{k_1,\cdots,k_p}\rightarrow \mathcal{R}^1_{k_1,\cdots,k_p}, k_1,\cdots,k_p\geq 0$, we choose a smoothly varying Floer data. Moreover, the Floer data are required to be \emph{consistent} in the sense that over each boundary strata of $\mathcal{R}^1_{k_1,\cdots,k_p}$ that decomposes as a product, it agrees (up to conformal equivalence) with the product of Floer data on the lower dimensional strata (it is understood that we have fixed consistent labeled Floer data for the $\mathcal{R}^d$'s).\par\indent
We briefly discuss the existence of such a consistent choice, even though this is a standard argument, cf. \cite[Section (9g),(9i)]{Sei1}. For each Lagrangian labeling $\mathbf{L}$, there is a fiber bundle
\begin{equation}
F^{\mathbf{L}}(\mathcal{R}^1_{k_1,\cdots,k_p})\rightarrow (\mathcal{R}^1_{k_1,\cdots,k_p})_{\mathbf{L}}
\end{equation}
whose fiber over $r\in (\mathcal{R}^1_{k_1,\cdots,k_p})_{\mathbf{L}}:=\mathcal{R}^1_{k_1,\cdots,k_p}$ is the space $F^{\mathbf{L}}_{S_r}$ of labeled Floer datum on $S_r$. Note that $F^{\mathbf{L}}_{S_r}$ is contractible. A choice of Floer datum for $(\mathcal{R}^1_{k_1,\cdots,k_p})_{\mathbf{L}}$ is just a section of (4.7). The existence argument is done inductively as follows.\par\indent
Suppose we have chosen consistent labeled Floer data for $\mathcal{R}^1_{k_1,\cdots,k_p}$ for all $k_1+\cdots+k_p<N$. Fix some $k_1+\cdots+k_p=N$. Then for each codimension 1 boundary component of $\mathcal{R}^1_{k_1,\cdots,k_p}$ (corresponding to a tree type with one interior node), one defines a section of $F^{\mathbf{L}}(\mathcal{R}^1_{k_1,\cdots,k_p})\rightarrow (\mathcal{R}^1_{k_1,\cdots,k_p})_{\mathbf{L}}$ in a collar neighborhood of that boundary component by gluing at the node. Do this for all codimension 1 boundary components, then by the inductive hypothesis, these local sections agree whenever their domain of definition overlaps (which is a neighborhood of some codimension 2 corner). Therefore, we obtained a section of (4.7) in a collar neighborhood of the entire codimension 1 boundary. Since the fibers of (4.7) are contractible, we can extend this section to all of $(\mathcal{R}^1_{k_1,\cdots,k_p})_{\mathbf{L}}$. By construction, these inductive choices of Floer data are consistent. \par\indent
Fix a consistent choice of Floer data for $\overline{\mathcal{R}}^1_{k_1,\cdots,k_p}$. Let $y_{out}\in CM^*(f)$ be a Morse cochain and
$$\mathbf{x}=\{x^1,x^1_1,\cdots,x^1_{k_1},\cdots,x^p,x^p_1,\cdots,x^p_{k_p}\}$$
be a cyclically composable sequence of morphisms in $\mathcal{F}$. Let
\begin{equation}
\mathcal{M}(\mathcal{R}^1_{k_1,\cdots,k_p}, y_{out}, \mathbf{x})
\end{equation}
be the moduli space of $(r,u)$, where $r\in  \mathcal{R}^1_{k_1,\cdots,k_p}$ and $u: S_r\rightarrow X$ satisfies Floer's equation (3.6) with appropriate Lagrangian boundary conditions and asymptotic conditions specified by $\mathbf{x}$, and that the interior marked point is constrained at $W^u(y_{out})$. For generically chosen Floer data, one can ensure that the above moduli space is regular of dimension
\begin{equation}
k_1+\cdots+k_p+\mathrm{ind}(u)+|y_{out}|-n.
\end{equation}
When $\mathrm{ind}(u)=n-|y_{out}|-k_1-\cdots-k_p$, counting rigid elements in $\mathcal{M}(\mathcal{R}^1_{k_1,\cdots,k_p}, y_{out}, \mathbf{x})$ gives rise to a map
\begin{equation}
{}_p\mathbf{F}_{k_1,\cdots,k_p}: \mathcal{F}_{\Delta}\otimes \mathcal{F}[1]^{\otimes k_1}\otimes\cdots\otimes \mathcal{F}_{\Delta}\otimes \mathcal{F}[1]^{\otimes k_p}\rightarrow CM^*(f)
\end{equation}
of degree $n$. We define ${}_pOC_{k_1,\cdots,k_p}:=(-1)^{\Vec{t}_{{}_pOC}}{}_p\mathbf{F}_{k_1,\cdots,k_p}$, where $\Vec{t}_{{}_pOC}$ denotes the sign twisting datum (compare \cite[(5.66)]{Gan1}, \cite[section 8]{RS}) given by
$
\Vec{t}_{{}_pOC}(\mathbf{x}):=\sum_{i=1}^p((k_1+\cdots+k_{i-1}+1)|x^i|+\sum_{j=1}^{k_i}(k_1+\cdots+k_{i-1}+j+1)|x^i_j|)$.\par\indent
Let ${}_pOC: {}_pCC_*(\mathcal{F})\rightarrow CM^*(f)$ denote the sum $\sum_{k_1,\cdots,k_p}{}_pOC_{k_1,\cdots,k_p}$. 
\begin{prop}
${}_pOC$ is a chain map of degree $n$ (mod 2).
\end{prop}
\emph{Proof Sketch}. To show that ${}_pOC$ is a chain map, consider the case when the moduli space (4.8) is one dimensional. By the consistency condition on Floer data, the boundary of the Gromov compactification of $\mathcal{M}(\mathcal{R}^1_{k_1,\cdots,k_p}, y_{out}, \mathbf{x})$ consists of contributions from domain breaking according to (4.5), (4.6) and semistable strip breaking, which constitute the term ${}_pOC\circ d_{{}_pCC}$; contributions from breaking of Morse trajectory (or equivalently incidence condition constrained to $\partial W^u(y_{out})$), which constitute the term $d_{CM^*(f)}\circ {}_pOC$. Disk and sphere bubbling are excluded since, by monotonicity, a disk or sphere bubbling off decreases the index of the main component by at least $2$, and therefore generically the main component (which is always stable) cannot intersect $W^u(y_{out})$ by our regularity assumption. This implies that, up to sign, $d_{CM*(f)}\circ {}_pOC-{}_pOC\circ d_{{}_pCC}$ is nullhomotopic. The sign verification follows from a tedious computation analogous to \cite[Appendix B]{Gan1}.\qed

\subsection{The $\mathbb{Z}/p\mathbb{Z}$-equivariant open-closed map}
There is $\mathbb{Z}/p$-action on
\begin{equation}
\coprod_{\mathbf{L},k_1,\cdots,k_p} \mathcal{R}^1_{k_1,\cdots,k_p}  
\end{equation}
such that for $r\in\mathcal{R}^1_{k_1,\cdots,k_p}$, the standard generator $\tau\in\mathbb{Z}/p$ acts by rotating the standard representative of $r$ by $e^{2\pi i/p}$, and the Lagrangian labels are rotated accordingly. It is clear that this action uniquely extends to the compactification $\coprod_{\mathbf{L},k_1,\cdots,k_p} \overline{\mathcal{R}}^1_{k_1,\cdots,k_p}$. Let $\sigma_r: S_r\rightarrow S_{\tau(r)}$ be the rotation map, where $S_r, S_{\tau(r)}$ are standard representatives of $r,\tau(r)$, respectively. Then, there is a $\mathbb{Z}/p$-action on the set of consistent choices of Floer data for $\overline{\mathcal{R}}^1_{k_1,\cdots,k_p}$, denoted $F_{{}_pOC}$, where the standard generator $\tau\in\mathbb{Z}/p$ acts by
\begin{equation}
\tau(\epsilon, K,J)_r=(\epsilon_{\tau(r)}\circ \sigma_r, \sigma_r^*K_{\tau(r)}, J_{\tau(r)}\circ \sigma_r).
\end{equation}
Suppose we have an $S^{\infty}$-dependent consistent choices of Floer data for the spaces $\overline{\mathcal{R}}^1_{k_1,\cdots,k_p}$, we say that it is a \emph{$\mathbb{Z}/p$-equivariant Floer data for ${}_pOC$} if 
\begin{equation}
\epsilon_{\tau(w),r}=\epsilon_{w,\tau(r)}\circ \sigma_r,\;\; J_{\tau(w),r}=J_{w,\tau(r)}\circ \sigma_r,\;\;K_{\tau(w),r}=\sigma_r^*K_{w,\tau(r)},
\end{equation}
where $w\in S^{\infty}$. 
\begin{prop}
There exists an $S^{\infty}$-dependent consistent choice of $\mathbb{Z}/p$-equivariant Floer data for ${}_pOC$.    
\end{prop} 
\emph{Proof}. The argument is by induction over the spaces $\overline{\mathcal{R}}^1_{k_1,\cdots,k_p}$ and cells $\Delta_i\subset S^{\infty}$. We first fix a consistent choice of Floer data for the spaces $\overline{\mathcal{R}}^1_{k_1,\cdots,k_p}$ constructed in section 4.1, and designate that to be the choice of Floer data for $\overline{\mathcal{R}}^1_{k_1,\cdots,k_p}\times \Delta_0$. By $\mathbb{Z}/p$-equivariance, we define the Floer data for $\overline{\mathcal{R}}^1_{k_1,\cdots,k_p}\times \tau\Delta_0$ to be pulled back from $\overline{\mathcal{R}}^1_{k_1,\cdots,k_p}\times \Delta_0$ via the action of $\tau$ on $\mathcal{F}_{{}_pOC}$, cf. (4.13).\par\indent
Suppose we have chosen Floer data that is consistent with boundary decomposition and $\mathbb{Z}/p$-equivariant for all $\overline{\mathcal{R}}^1_{k_1',\cdots,k_p'}\times \tau^j\Delta_{i'}$ with $(k_1',\cdots,k_p',i')<(k_1,\cdots,k_p,i)$ and $0\leq j\leq p-1$ if $i'$ is odd and $0\leq j\leq 1$ if $i'$ is even. Now we inductively determine the Floer data for $\overline{\mathcal{R}}^1_{k_1,\cdots,k_p}\times \tau^j\Delta_{i}$ where $0\leq j\leq 1$ if $i$ is even and $0\leq j\leq p-1$ if $i$ is odd. We discuss the case when $i=2k+1$ is odd, as the even case is completely analogous. \par\indent
The boundary of $\overline{\mathcal{R}}^1_{k_1,\cdots,k_p}\times \Delta_{2k+1}$ is covered by
\begin{equation}
(\overline{\mathcal{R}}^{k'_i+1}\times \overline{\mathcal{R}}^1_{k_1,\cdots,k_i-k_i'+1,\cdots,k_p})\times \Delta_{2k+1},\quad 1\leq i\leq p,
\end{equation}
\begin{equation}
(\overline{\mathcal{R}}^{k_{i-1}'+k_{i}'+2}\times \overline{\mathcal{R}}^1_{k_1,\cdots,k_{i-1}-k_{i-1}',k_{i}-k_{i}',\cdots,k_p})\times \Delta_{2k+1},\quad 1\leq i\leq p.
\end{equation}
and
\begin{equation}
\overline{\mathcal{R}}^1_{k_1,\cdots,k_p}\times \Delta_{2k}, \overline{\mathcal{R}}^1_{k_1,\cdots,k_p}\times \tau\Delta_{2k}.  
\end{equation}
By induction hypothesis, we have chosen Floer data for (4.14)-(4.16). Fix a small $\epsilon$ and apply gluing near the codimension $1$ boundary using the prescribed strip-like ends to (4.14) and (4.15), we obtain a choice Floer data for $U_\epsilon(\partial \overline{\mathcal{R}}^1_{k_1,\cdots,k_p})\times \Delta_{2k+1}$, where $U_\epsilon(\partial \overline{\mathcal{R}}^1_{k_1,\cdots,k_p})$ denotes an $\epsilon$-collar neighborhood of the boundary. We then extend the Floer data on 
\begin{equation}
U_\epsilon(\partial \overline{\mathcal{R}}^1_{k_1,\cdots,k_p})\times \Delta_{2k+1}\cup \overline{\mathcal{R}}^1_{k_1,\cdots,k_p}\times \Delta_{2k}\cup\overline{\mathcal{R}}^1_{k_1,\cdots,k_p}\times \tau\Delta_{2k} 
\end{equation} 
smoothly to the entire $\overline{\mathcal{R}}^1_{k_1,\cdots,k_p}\times \Delta_{2k+1}$; this can be done since the space of Floer data on $S_r$ for each $r\in \mathcal{R}^1_{k_1,\cdots,k_p}$ is contractible. Then, we define the Floer data on $\overline{\mathcal{R}}^1_{k_1,\cdots,k_p}\times \tau^j\Delta_{2k+1}, 0\leq j\leq p-1$ to be pulled back from $\overline{\mathcal{R}}^1_{k_1,\cdots,k_p}\times \Delta_{2k+1}$ via the action of $\tau^j$ on $\mathcal{F}_{{}_pOC}$. Note that this does not cause inconsistencies: for $i<j$, the intersection of $\overline{\mathcal{R}}^1_{k_1,\cdots,k_p}\times \tau^i\Delta_{2k+1}$ with $\overline{\mathcal{R}}^1_{k_1,\cdots,k_p}\times \tau^j\Delta_{2k+1}$ is empty unless $i=j-1$, in which case it is $\overline{\mathcal{R}}^1_{k_1,\cdots,k_p}\times \tau^{j-1}\Delta_{2k}$. Moreover, the restriction of the Floer data on $\overline{\mathcal{R}}^1_{k_1,\cdots,k_p}\times \tau^l\Delta_{2k+1}, l=j-1,j$ to $\overline{\mathcal{R}}^1_{k_1,\cdots,k_p}\times \tau^{j-1}\Delta_{2k}$ agree, and are both equivalent to the inductively chosen Floer data on this stratum. This essentially uses the feature that the $\mathbb{Z}/p$-action on $S^{\infty}$ is free. \qed\par\indent
Let $y_{out}\in CM^*(f)$ and
$$\mathbf{x}=\{x^1,x^1_1,\cdots,x^1_{k_1},\cdots,x^p,x^p_1,\cdots,x^p_{k_p}\}\in {}_pCC_*(\mathcal{F}).$$
We define
\begin{equation}
\mathcal{M}(\Delta_i\times \mathcal{R}^1_{k_1,\cdots,k_p},y_{out},\mathbf{x})
\end{equation}
to be the moduli space of $(w,r,u)$, where $w\in \Delta_i, r\in\mathcal{R}^1_{k_1,\cdots,k_p}$ and $u:S_r\rightarrow X$ satisfies
\begin{equation}
(du-Y_{K_{w,r}})^{0,1}_{J_{w,r}}=0
\end{equation}
with appropriate Lagrangian boundary conditions and asymptotics specified by $\mathbf{x}$, and that the interior marked points is constrained on $W^u(y_{out})$. Choosing a generic $\mathbb{Z}/p$-equivariant Floer data for ${}_pOC$, the moduli space (4.18) is regular of dimension
\begin{equation}
k_1+\cdots+k_p+i+\textrm{ind}(u)+|y_{out}|-n.
\end{equation}
When $\textrm{ind}(u)=n-i-|y_{out}|-k_1-\cdots-k_p$, counting rigid elements of (4.18), adjusted by the sign twisting datum $\Vec{t}_{{}_pOC}$, gives a map
\begin{equation}
{}_pOC^i_{k_1,\cdots,k_p}: {}_pCC_*(\mathcal{F})\rightarrow CM^*(f)
\end{equation}
of degree $n-i$ (mod 2). Let
\begin{equation}
{}_pOC^i=\sum_{k_1,\cdots,k_p}{}_pOC^i_{k_1,\cdots,k_p}.
\end{equation}
By definition, ${}_pOC^0={}_pOC$. \par\indent
The codimension one boundary strata of the Gromov compactification $\overline{\mathcal{M}}(\Delta_i\times \mathcal{R}^1_{k_1,\cdots,k_p},y_{out},\mathbf{x})$ consists of contributions coming from domain breaking (4.5) and (4.6), the corresponding moduli space where the parameter $w$ is constrained to $\partial \Delta_i$, semi-stable strip breaking, and Morse trajectory breaking. Note that there is a canonical identification
\begin{equation}
\mathcal{M}(\tau^j(\Delta_i)\times \mathcal{R}^1_{k_1,\cdots,k_p},y_{out},\mathbf{x})\xrightarrow{\cong} \mathcal{M}(\Delta_i\times \mathcal{R}^1_{k_{p-j+1},\cdots,k_p,k_1,\cdots,k_j},y_{out},\tau^j(\mathbf{x}))
\end{equation}
coming from the $\mathbb{Z}/p$-equivariance of Floer data.
Combining all of the above, and using formula (3.13) and (3.14), we deduce that for $i\geq 0$:
\begin{equation}
d_{CM^*(f)}\circ {}_pOC^i- {}_pOC^i\circ d_{{}_pCC_*}=\\
\begin{cases}
{}_pOC^{i-1}\circ (\tau-1),\quad\textrm{when $i$ is odd}\\
{}_pOC^{i-1}\circ (1+\tau+\cdot+\tau^{p-1}),\quad\textrm{when $i$ is even}.
\end{cases}
\end{equation}
We define ${}_pOC^{-1}=0$ so that the above relation is satisfied for $i=0$ (since ${}_pOC^0={}_pOC$ is a chain map by Proposition 4.1). \par\indent
\begin{mydef}
The $\mathbb{Z}/p$\emph{-equivariant open-closed map}
\begin{equation}
OC^{\mathbb{Z}/p}: CC^{\mathbb{Z}/p}_*(\mathcal{F})\rightarrow CM^*(f)[[t,\theta]]
\end{equation}
is defined by
\begin{equation}
\mathbf{x}\mapsto \sum_k \big({}_pOC^{2k}(\mathbf{x})+(-1)^{\|\mathbf{x}\|} {}_pOC^{2k+1}(\mathbf{x})\theta\big)t^k,
\end{equation}
\begin{equation}
\mathbf{x}\theta\mapsto \sum_k \big({}_pOC^{2k}(\mathbf{x})\theta+(-1)^{\|\mathbf{x}\|} \sum_{j=0}^{p-2} {}_pOC^{2k+1} (\mathbf{x}+\tau(\mathbf{x})+\cdot+\tau^j(\mathbf{x}))t\big)t^k,
\end{equation}
and extended $t$-linearly.
\end{mydef}
\begin{prop}
$OC^{\mathbb{Z}/p}$ is a chain map.
\end{prop}
\noindent\emph{Proof}. This is an immediate consequence of (4.24). \qed

\subsection{The closed-open maps}
We briefly review two versions of the closed-open string maps, cf.  \cite[section 5.4-5.6]{Gan1},\cite[section 5.5]{RS}. Let
\begin{equation}
\mathcal{R}^{1,1}_d
\end{equation}
be the moduli space of disks with
\begin{itemize}
    \item $d+1$ boundary marked points $z^-_0,z_1,\cdots,z_d$ labeled counterclockwise such that $z^-_0$ is marked as output and $z_1,\cdots,z_d$ are marked as input,
    \item one interior marked point $y_{in}$ marked as input.
\end{itemize}
After choosing consistent regular Floer data for (4.28), counting rigid solutions to the Floer equation with the usual conditions gives rise to a chain map
\begin{equation}
CO:CM^*(f)\rightarrow CC^*(\mathcal{F})
\end{equation}
of degree $0$. At the level of cohomology, $CO$ is ring map, \cite[Proposition 5.3]{Gan1}. \par\indent
On the other hand, let
\begin{equation}
\mathcal{R}^{1,1}_{r,s}
\end{equation}
be the moduli space of disks with one interior input $y_{in}$, one boundary output $z_{out}$, and $r+s+1$ boundary inputs $z_1,\cdots,z_r,z_{fixed},z_1',\cdots,z_s'$ labeled in clockwise order from $z_{out}$, such that up to automorphism of the disk, $z_{out},z_{fixed}$ and $y_{in}$ lie at $-i,i,0$, respectively. \par\indent
Counting rigid solutions to the moduli problem parametrized by (4.30) gives rise to a chain map
\begin{equation}
{}_2CO:CM^*(f)\rightarrow {}_2CC^*(\mathcal{F}):=\mathrm{hom}_{\mathcal{F}-\mathcal{F}}(\mathcal{F}_{\Delta},\mathcal{F}_{\Delta})
\end{equation}
of degree $0$. Moreover, if $\Psi: CC^*(\mathcal{F})\rightarrow {}_2CC^*(\mathcal{F})$ is the comparison map of (2.18), then
\begin{lemma}\cite[Proposition 5.6]{Gan1}
$\Psi\circ CO$ and ${}_2CO$ are homotopic.\qed
\end{lemma} 

\subsection{Compatibility with quantum Steenrod operations}
Fix a Morse cocycle $b\in CM^*(f)$, then $CO(b)^{\otimes p}$ and ${}_2CO(b)^{\otimes p}$ are cocycles in $(CC^*(\mathcal{F})^{\otimes p})^{\mathbb{Z}/p}$ and $({}_2CC^*(\mathcal{F})^{\otimes p})^{\mathbb{Z}/p}$, respectively. Moreover, by Lemma 2.10, their cohomology class only depends on the cohomology class of $b$. The main result of this paper can be re-stated as follows
\begin{thm}
The diagram
\begin{center}
\begin{tikzcd}[row sep=1.2cm, column sep=0.8cm]
{}_pCC^{\mathbb{Z}/p}_*(\mathcal{F})\arrow[rrrr,"{OC^{\mathbb{Z}/p}}"]\arrow[d,"{\prod^{\mathbb{Z}/p}({}_2CO(b)^{\otimes p},-)}"]& & & &CM^*(f)[[t,\theta]]\arrow[d,"\Sigma_b"]  \\
{}_pCC_*^{\mathbb{Z}/p}(\mathcal{F})\arrow[rrrr,"{OC^{\mathbb{Z}/p}}"]& & & &CM^*(f)[[t,\theta]]
\end{tikzcd}
\end{center}
is homotopy commutative.
\end{thm}
\emph{Proof of Theorem 1.2 given Theorem 4.6}. By definition of $\bigcap^{\mathbb{Z}/p}$ we have
\begin{equation}
\bigcap^{\mathbb{Z}/p}(CO(b),-)=\prod_{eq}(\Psi(CO(b))^{\otimes p},-).
\end{equation}
But Lemma 4.5 implies that $\Psi(CO(b))$ and ${}_2CO(b)$ are homologous in ${}_2CC^*(\mathcal{F})$, and Lemma 2.10 further implies that $\Psi(CO(b))^{\otimes p}$ and ${}_2CO(b)^{\otimes p}$ are homologous in $({}_2CC^*(\mathcal{F})^{\otimes p})^{\mathbb{Z}/p}$. This proves Theorem 1.2. \qed\par\indent
The rest of section 4.4 will be devoted to the proof of Theorem 4.6. We first introduce a new parameter space of disks that will be used in constructing the homotopy in Theorem 4.6. Let
\begin{equation}
\mathcal{Q}^1_{p,k_1,\cdots,k_p}
\end{equation}
be the parameter space of disks with
\begin{itemize}
    \item one interior output $y_{out}$,
    \item $p$ interior inputs $y_1,\cdots,y_p$,
    \item $k_1+\cdots+k_p+p$ boundary inputs $z^1,z^1_1,\cdots,z^1_{k_1},z^2,z^2_1,\cdots,z^2_{k_2},\cdots,z^p,z^p_1,\cdots,z^p_{k_p}$ in counterclockwise order such that up to reparametrization of the disk, $y_{out},y_1,\cdots,y_p,z^1,\cdots,z^p$ lie respectively at $0,r\zeta, \cdots,r\zeta^p, \zeta,\cdots,\zeta^p$ for some $r\in (0,1)$. As before, this particular parametrization of an element will be called the standard representative. 
\end{itemize}
We fix the orientation of $\mathcal{Q}^1_{p,k_1,\cdots,k_p}$ induced by $-dr\wedge dz^1_1\wedge\cdots\wedge dz^1_{k_1}\wedge\cdots\wedge dz^p_1\wedge\cdots\wedge dz^p_{k_p}$. \par\indent
$\mathcal{Q}^1_{p,k_1,\cdots,k_p}$ admits a Deligne-Mumford compactification $\overline{\mathcal{Q}}^1_{p,k_1,\cdots,k_p}$ which is naturally equipped with a submersion $\overline{\mathcal{Q}}^1_{p,k_1,\cdots,k_p}\rightarrow [0,1]$. The boundary strata of $\overline{\mathcal{Q}}^1_{p,k_1,\cdots,k_p}$ is covered by
\begin{equation}
\overline{\mathcal{R}}^1_{k_1,\cdots,k_p}\times \overline{\mathcal{S}}^1_{p+1},\quad (r=0)
\end{equation}
\begin{equation}
\prod_{i=1}^p \overline{\mathcal{R}}^{1,1}_{s_{i-1}',t_i'}\times\overline{\mathcal{R}}^1_{k_1-t_1'-s_1',\cdots,k_p-t_p'-s_p'},\quad (r=1)
\end{equation}
\begin{equation}
\overline{\mathcal{R}}^{k'_i}\times \overline{\mathcal{Q}}^1_{p,k_1,\cdots,k_i-k_i'+1,\cdots,k_p},\quad 1\leq i\leq p,
\end{equation}
\begin{equation}
\overline{\mathcal{R}}^{k_{i-1}'+k_{i}'+1}\times \overline{\mathcal{Q}}^1_{p,k_1,\cdots,k_{i-1}-k_{i-1}',k_{i}-k_{i}',\cdots,k_p},\quad 1\leq i\leq p.
\end{equation}
We explain some of the notations in (4.34)-(4.37). In (4.34), the space $\overline{\mathcal{S}}^1_{p+1}=\mathcal{S}^1_{p+1}$ is just a singleton representing the domain curve $C=(\mathbb{C}P^1, x_0,x_1,\cdots,x_p,x_{\infty})$ (which underlies the definition of quantum Steenrod operations), where the special input $x_0$ lies at $\infty$, the other $p$-inputs  $x_1,\cdots,x_p$ lie at $\zeta,\cdots,\zeta^p$ (the $p$-th roots of unity along the equator) and the output $x_{\infty}$ lies at $\infty$. The output of the main component $y_{out}$ is glued to the input $x_0$ of $C$, and $x_1,\cdots,x_p$ become the new distinguished interior inputs after gluing. In (4.35), the output of a disk in $\overline{\mathcal{R}}^{1,1}_{s_{i-1}',t_i'}$ is glued to the $i$-th distinguished boundary input of the main component, and $z_{fixed}$ will become the new distinguished $i$-th boundary input after gluing. In (4.36), the output of the first disk can be glued to any one of the inputs of the main component labeled by $z^i_1,\cdots,z^i_{k_i-k_i'+1}$. In (4.37), the output of the first disk is glued to the distinguished boundary input $z^i$, and the $k_i+1$-th input of the first disk will become the new distinguished $i$-th boundary input after gluing.
\begin{figure}[H]
 \centering
 \includegraphics[width=1.0\textwidth]{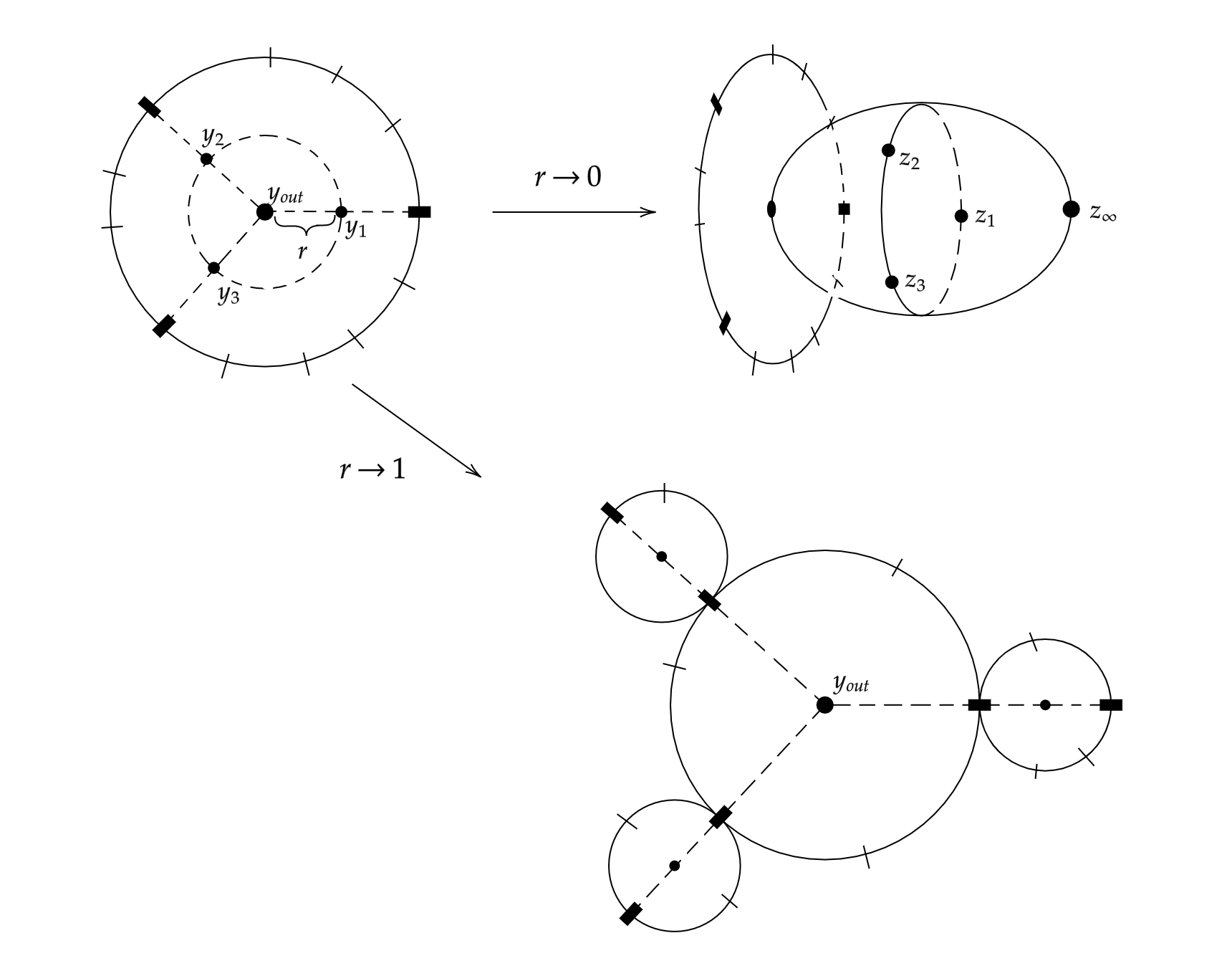}
 \caption{An element of $\mathcal{Q}^1_{p,k_1,\cdots,k_p}$ degenerates as $r\rightarrow 0$ or $r\rightarrow 1$}
\end{figure}
We a priori fix consistent Floer data for the parameter spaces $\mathcal{R}^{d+1}, \mathcal{R}^{1,1}_{r,s}$ as well as consistent $\mathbb{Z}/p$-equivariant ($S^{\infty}$-dependent) Floer data for $\mathcal{R}^1_{k_1,\cdots,k_p}$ and $\mathcal{S}^1_{p+1}$. \par\indent
A choice of $S^{\infty}$-dependent Floer data for the spaces $\mathcal{Q}^1_{p,k_1,\cdots,k_p}, k_1,\cdots,k_p\geq 0$, denoted $(\epsilon_{w,r},K_{w,r},J_{w,r}), w\in S^{\infty}, r\in \overline{\mathcal{Q}}^1_{p,k_1,\cdots,k_p}$, is called \emph{consistent and $\mathbb{Z}/p$-equivariant} if:
\begin{itemize}
    \item It is $\mathbb{Z}/p$-equivariant, i.e.
\begin{equation}
\epsilon_{\tau(w),r}=\epsilon_{w,\tau(r)}\circ \sigma_r,\;\; J_{\tau(w),r}=J_{w,\tau(r)}\circ \sigma_r,\;\;K_{\tau(w),r}=\sigma_r^*K_{w,\tau(r)},
\end{equation}
where $\tau$ denotes the action of the standard generator of $\mathbb{Z}/p$ on $\overline{\mathcal{Q}}^1_{p,k_1,\cdots,k_p}$.
\item When the underlying domain $r$ approaches a boundary strata of type (4.35)-(4.37), the Floer data $(\epsilon_{w,r},K_{w,r},J_{w,r}),w\in S^{\infty}$ decomposes as a product of $(\epsilon_{w,r'},K_{w,r'},J_{w,r'}), w\in S^{\infty}, r'\in \overline{Q}$ (of lower dimension) and the pre-chosen non-equivariant Floer data on the other components. 
\item When the underlying domain $r$ approaches a boundary strata of type (4.34), we require that the Floer data $(\epsilon_{w,r},K_{w,r},J_{w,r}), w\in S^{\infty}$ decomposes as the product of $(\epsilon_{w,r'},K_{w,r'},J_{w,r'}), r'\in \mathcal{R}^1_{p,k_1,\cdots,k_p}$ and $(\epsilon_{w,C},K_{w,C},J_{w,C}), C\in\mathcal{S}^1_{p+1}$ (note the $S^{\infty}$-parameter on both components agree with $w$).
\end{itemize}
The following proposition follows a similar inductive argument as in Proposition 4.2.
\begin{prop}
An $S^{\infty}$-dependent, consistent and $\mathbb{Z}/p$-equivariant Floer data on $\overline{\mathcal{Q}}^1_{k_1,\cdots,k_p}$, $k_1,\cdots,k_p\geq 0$ exists.\qed
\end{prop}

\subsubsection*{Moduli spaces with Floer data parametrized by chains in $S^{\infty}\times S^{\infty}$}
We now define two auxiliary moduli spaces associated with nodal configurations. Fix Morse cochains $y_{\infty},y_1,\cdots,y_p\in CM^*(f)$, a $p$-fold Hochschild chain $\mathbf{x}\in {}_pCC_*(\mathcal{F})$ of type $(k_1,\cdots,k_p)$, and a smooth singular chain $\sigma=(\sigma_1,\sigma_2):\Delta[i]\rightarrow S^{\infty}\times S^{\infty}$, where $\Delta[i]$ is the standard $i$-th simplex in $\mathbb{R}^{i+1}$. We let
\begin{equation}
\mathcal{M}_{\sharp}(\sigma, \mathcal{R}^1_{k_1,\cdots,k_p}\times \mathcal{S}^1_{p+1},y_{\infty},y_1,\cdots,y_p,\mathbf{x})
\end{equation}
be the moduli space of
\begin{equation}
w\in \mathrm{int}(\Delta[i]),\;r\in \mathcal{R}^1_{k_1,\cdots,k_p}, \;u_1:\mathcal{S}_r\rightarrow X, \;u_2:C=\mathbb{C}P^1\rightarrow X    
\end{equation}
satisfying
\begin{equation}
\begin{cases}
(du_1-Y_{K_{\sigma_1(w),\mathcal{S}_r}})^{0,1}_{J_{\sigma_1(w),\mathcal{S}_r}}=0,\quad (du_2-Y_{K_{\sigma_2(w),C}})^{0,1}_{J_{\sigma_2(w),C}}=0,\\
\textrm{appropriate Lagrangian boundary conditions for}\;\mathcal{S}_r,\\
\textrm{boundary marked points of}\;\mathcal{S}_r\;\textrm{are asymptotic to}\;\mathbf{x},\\
\textrm{the points}\;z_{\infty},z_1,\cdots,z_p\;\textrm{on C are constrained to}\;W^u(y_{\infty}),W^s(y_1),\cdots,W^s(y_p),\\
u_1(y_{out})=u_2(z_0).
\end{cases}
\end{equation}
Assuming the moduli space (4.39) is regular in dimension $0$, let
\begin{equation}
(\Sigma\sharp{}_pOC)_{\sigma}:CM^*(f)^{\otimes p}\otimes {}_pCC_*(\mathcal{F})\rightarrow CM^*(f)
\end{equation}
be the corresponding operation defined by counting rigid solutions of (4.39). By linearity, we can define $(\Sigma\sharp{}_pOC)_{\sigma}$ for $\sigma$ any $\mathbf{k}$-coefficient chain in $S^{\infty}\times S^{\infty}$. Similarly, we define
\begin{equation}
\mathcal{M}_{\star}(\sigma,\mathcal{R}^1_{k_1,\cdots,k_p}\times \mathcal{S}^1_{p+1},y_{\infty},y_1,\cdots,y_p,\mathbf{x})
\end{equation}
to be the union over $y\in CM^*(f)$ of
\begin{equation}
w\in \mathrm{int}(\Delta[i]),\; r\in \mathcal{R}^1_{k_1,\cdots,k_p},\; u_1:\mathcal{S}_r\rightarrow X,\; u_2:C=\mathbb{C}P^1\rightarrow X    
\end{equation}
satisfying (all except the last condition are the same as (4.41))
\begin{equation}
\begin{cases}
(du_1-Y_{K_{\sigma_1(w),\mathcal{S}_r}})^{0,1}_{J_{\sigma_1(w),\mathcal{S}_r}}=0,\quad (du_2-Y_{K_{\sigma_2(w),C}})^{0,1}_{J_{\sigma_2(w),C}}=0,\\
\textrm{appropriate Lagrangian boundary conditions for}\;\mathcal{S}_r,\\
\textrm{boundary marked points of}\;\mathcal{S}_r\;\textrm{are asymptotic to}\;\mathbf{x},\\
\textrm{the points}\;z_{\infty},z_1,\cdots,z_p\;\textrm{on C are constrained to}\;W^u(y_{\infty}),W^s(y_1),\cdots,W^s(y_p),\\
u_1(y_{out})\;\textrm{is constrained to}\;W^u(y),\;u_2(z_0)\;\textrm{is constrained to}\;W^s(y).
\end{cases}
\end{equation}
Denote the corresponding operation by counting rigid elements of (4.43) by
\begin{equation}
(\Sigma\star{}_pOC)_{\sigma}: CM^*(f)^{\otimes p}\otimes {}_pCC_*(\mathcal{F})\rightarrow CM^*(f)
\end{equation}
and extend by linearity to all $\mathbf{k}$-coefficient chains of $S^{\infty}\times S^{\infty}$.
\subsubsection*{Proof of Theorem 4.6}
We construct the homotopy required in Theorem 4.6 in three steps.\par\indent
\emph{Step 1}. The first homotopy comes from a parametrized moduli problem associated to $\mathcal{Q}^1_{p,k_1,\cdots,k_p}$. Fix Morse cochains $a_{out},a_1,\cdots,a_p\in CM^*(f)$, $p$-fold Hochschild chain $\mathbf{x}\in {}_pCC_*(\mathcal{F})$ of type $(k_1,\cdots,k_p)$. Let
\begin{equation}
\mathcal{M}(\Delta_i,\mathcal{Q}^1_{p,k_1,\cdots,k_p},a_{\infty},a_1,\cdots,a_p,\mathbf{x})
\end{equation}
be the moduli space of
$$w\in \Delta_i\subset S^{\infty}, r\in \mathcal{Q}^1_{p,k_1,\cdots,k_p}, u:\mathcal{S}_r\rightarrow X$$
satisfying
\begin{equation}
\begin{cases}
(du-Y_{K_{w,\mathcal{S}_r}})^{0,1}_{J_{w,\mathcal{S}_r}}=0,\\
\textrm{appropriate Lagrangian boundary conditions for}\;\mathcal{S}_r,\\
\textrm{boundary marked points of}\;\mathcal{S}_r\;\textrm{are asymptotic to}\;\mathbf{x},\\
\textrm{interior marked points}\;y_{out},y_1,\cdots,y_p\;\textrm{are constrained to}\;W^u(a_{out}),W^s(a_1),\cdots,W^s(a_p).
\end{cases}
\end{equation}
Fixing a generic, $\mathbb{Z}/p$-equivariant and consistent choice of Floer data, (4.47) is regular of dimension
\begin{equation}
k_1+\cdots+k_p+\mathrm{ind}(u)+|y_{out}|-|y_1|-\cdots-|y_p|-n+i+1.
\end{equation}
Let
\begin{equation}
\mathcal{H}^i:CM^*(f)^{\otimes p}\otimes {}_pCC_*(\mathcal{F})\rightarrow {}_pCC_*(\mathcal{F})
\end{equation}
be the operation of degree $n-i-1$ (mod 2) obtained by counting rigid elements of (4.47). 
If $b\in CM^*(f)$ is a Morse cocycle, we also write
\begin{equation}
\mathcal{H}^i_b:=\mathcal{H}^i(\overbrace{b,\cdots,b}^{p\;\mathrm{times}},-): {}_pCC_*(\mathcal{F})\rightarrow {}_pCC_*(\mathcal{F}),
\end{equation}
which is a map of degree $n+p|b|-i-1$ (mod 2). \par\indent
By considering the boundary strata of the Gromov compactification of the $1$-dimensional moduli space of (4.47), we obtain that
$$
\mathcal{H}^i_b\circ d_{{}_pCC_*}\pm d_{CM^*(f)}\circ \mathcal{H}^i_b=(\Sigma\sharp{}_pOC)_{\delta(\Delta_i),b}\pm {}_pOC^i\circ\prod({}_2CO(b)^{\otimes p}\otimes-)$$
\begin{equation}
\pm
\begin{cases}
H^{i-1}_b\circ (\tau-1),\quad\textrm{when i is odd},\\
H^{i-1}_b\circ(1+\tau+\cdots+\tau^{p-1}),\quad\textrm{when i is even}.
\end{cases}
\end{equation}
In (4.52), $\mathcal{H}^i_b\circ d_{{}_pCC_*}$ comes from domain breaking in (4.36),(4.37) and semi-stable strip breaking. $d_{CM^*(f)}\circ \mathcal{H}^i_b$ comes from Morse trajectory breaking. $(\Sigma\sharp{}_pOC)_{\delta(\Delta_i),b}$ (where $\delta:S^{\infty}\rightarrow S^{\infty}\times S^{\infty}$ is the diagonal map) and ${}_pOC^i\circ\prod({}_2CO(b)^{\otimes p}\otimes-)$ comes from the limiting behavior $r\rightarrow 0$ and $r\rightarrow 1$, which corresponds to domain breaking of (4.34) and (4.35), respectively. Finally, the last term comes from the boundary strata where the parameter $w\in S^{\infty}$ is constrained to $\partial\Delta_i$, together with canonical identifications
$$\mathcal{M}(\tau^j(\Delta_i)\times \mathcal{Q}^1_{p,k_1,\cdots,k_p},y_{out},b,\cdots,b,\mathbf{x})\cong$$
\begin{equation}
\mathcal{M}(\Delta_i\times \mathcal{Q}^1_{p,k_{p-j+1},\cdots,k_p,k_1,\cdots,k_j},y_{out},b,\cdots,b,\tau^j(\mathbf{x})).
\end{equation}
As usual, we can put the $\mathcal{H}^i_b$ together into a $t$-linear map $H^{\mathbb{Z}/p}_b:CC^{\mathbb{Z}/p}_*(\mathcal{F})\rightarrow CM^*(f)[[t,\theta]]$ by
\begin{equation}
\begin{cases}
\mathbf{x}\mapsto \sum_k\big(\mathcal{H}^{2k}_b(\mathbf{x})+(-1)^{\|x\|} \mathcal{H}^{2k+1}_b(\mathbf{x})\theta\big)t^k\\
\mathbf{x}\theta\mapsto \sum_k\big(\mathcal{H}^{2k}_b(\mathbf{x})\theta+(-1)^{\|x\|} \mathcal{H}^{2k+1}_b((1-\tau)^{p-2}\mathbf{x})t\big)t^k.
\end{cases}
\end{equation}
Note that since we are in characteristic $p$, $1+\tau+\cdots+\tau^{p-1}=(\tau-1)^{p-1}$. (4.52) implies that $\mathcal{H}^{\mathbb{Z}/p}_b$ provides a homotopy between
\begin{equation}
OC^{\mathbb{Z}/p}\circ \prod^{\mathbb{Z}/p}({}_2CO(b)^{\otimes p}\otimes-)    
\end{equation}
and the $t$-linear map $(\Sigma\sharp{}_pOC)^{\mathbb{Z}/p}_{\delta(\Delta),b}: CC^{\mathbb{Z}/p}_*(\mathcal{F})\rightarrow CM^*(f)[[t,\theta]]$ defined by
\begin{equation}
\begin{cases}
\mathbf{x}\mapsto \sum_k\big((\Sigma\sharp{}_pOC)_{\delta(\Delta_{2k}),b}(\mathbf{x})+(-1)^{\|x\|} (\Sigma\sharp{}_pOC)_{\delta(\Delta_{2k+1}),b}(\mathbf{x})\theta\big)t^k\\
\mathbf{x}\theta\mapsto \sum_k\big((\Sigma\sharp{}_pOC)_{\delta(\Delta_{2k}),b}(\mathbf{x})\theta +(-1)^{\|x\|}(\Sigma\sharp{}_pOC)_{\delta(\Delta_{2k+1}),b}((1-\tau)^{p-2}\mathbf{x})t\big)t^k.
\end{cases}
\end{equation}

\emph{Step 2}. The second homotopy involves a standard trick of inserting a finite length Morse trajectory (and letting its length go to infinity), cf. \cite{PSS}.
\begin{figure}[H]
 \centering
 \includegraphics[width=1.0\textwidth]{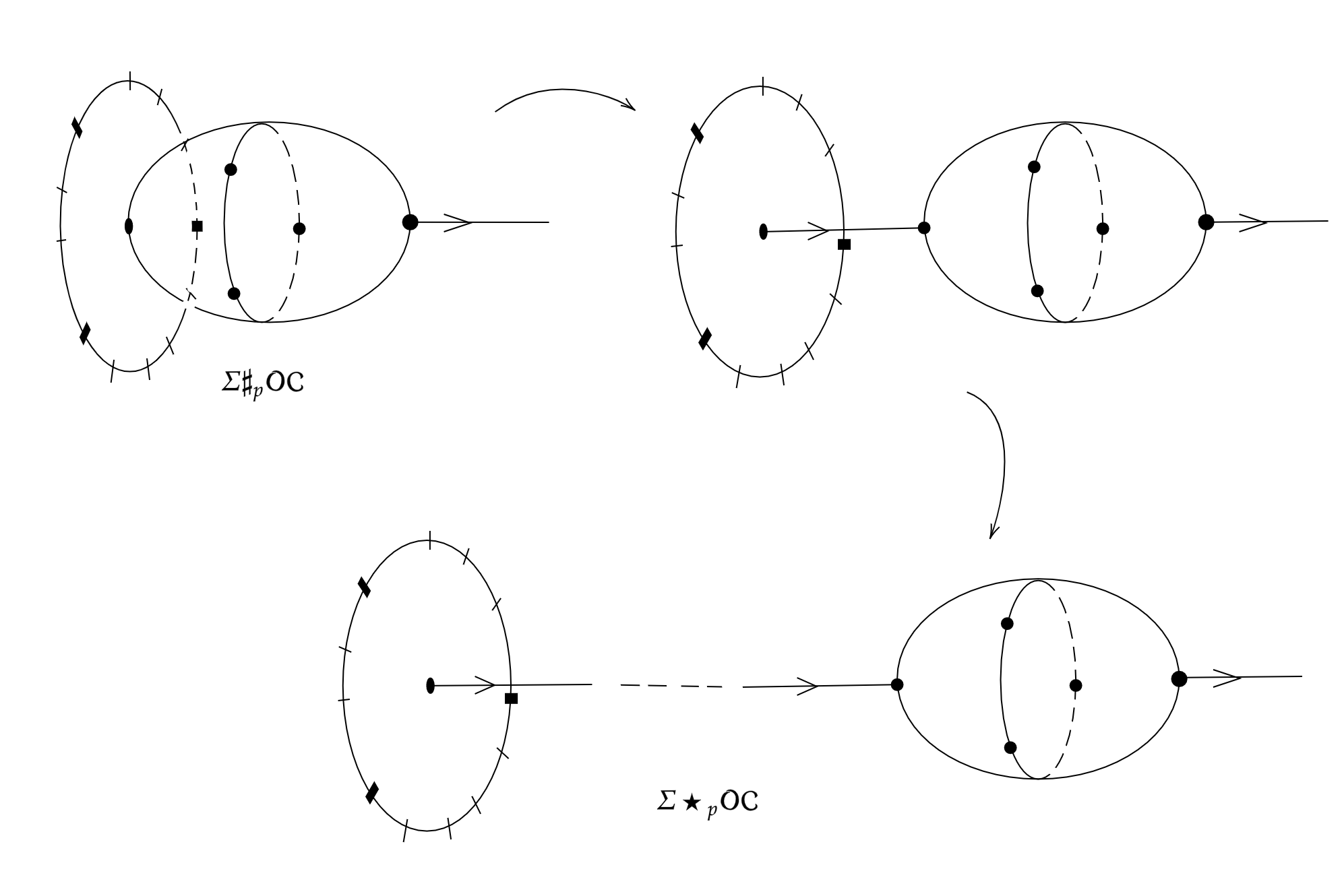}
 \caption{Inserting finite Morse trajectory with varying length}
\end{figure}
The goal is to construct a homotopy between $(\Sigma\sharp{}_pOC)^{\mathbb{Z}/p}_{\delta(\Delta),b}$ and $(\Sigma\star{}_pOC)^{\mathbb{Z}/p}_{\delta(\Delta),b}$. \par\indent
Fix Morse cocycle $b$, Morse cochain $y_{\infty}$ and $p$-fold Hochschild chain $\mathbf{x}$ of type $(k_1,\cdots,k_p)$ as before. Let
\begin{equation}
\mathcal{M}(\delta(\Delta_i), T,\mathcal{R}^1_{k_1,\cdots,k_p}\times\mathcal{S}^1_{p+1},y_{\infty},b,\cdots,b,\mathbf{x})
\end{equation}
be the moduli space of
$$w\in \Delta_i\subset S^{\infty}, T\in(0,\infty), \gamma:[0,T]\rightarrow X, r\in \mathcal{R}^1_{k_1,\cdots,k_p}, u_1:\mathcal{S}_r\rightarrow X, u_2:C=\mathbb{C}P^1\rightarrow X$$
satisfying
\begin{equation}
\begin{cases}
\dot{\gamma}=-\nabla f,\\
(du_1-Y_{K_{w,\mathcal{S}_r}})^{0,1}_{J_{w,\mathcal{S}_r}}=0,\quad (du_2-Y_{K_{w,C}})^{0,1}_{J_{w,C}}=0,\\
\textrm{appropriate Lagrangian boundary conditions for}\;\mathcal{S}_r,\\
\textrm{boundary marked points of}\;\mathcal{S}_r\;\textrm{are asymptotic to}\;\mathbf{x},\\
\textrm{the points}\;z_{\infty},z_1,\cdots,z_p\;\textrm{on C are constrained to}\;W^u(y_{\infty}),W^s(b),\cdots,W^s(b),\\
u_1(y_{out})=\gamma(0), u_2(z_0)=\gamma(T).
\end{cases}
\end{equation}
Then, as $T\rightarrow 0$, we obtain the nodal moduli space defining $(\Sigma\sharp{}_pOC)_{\delta(\Delta_i),b}$ and as $T\rightarrow \infty$, the finite Morse trajectory breaks into two semi-infinite Morse trajectories meeting at some critical point, giving rise to $(\Sigma\star{}_pOC)_{\delta(\Delta_i),b}$. Putting these together into an equivariant operation (cf. (4.54)), we obtain a homotopy between $(\Sigma\sharp{}_pOC)^{\mathbb{Z}/p}_{\delta(\Delta),b}$ and $(\Sigma\star{}_pOC)^{\mathbb{Z}/p}_{\delta(\Delta),b}$.\par\indent
\emph{Step 3}. Finally we show that $(\Sigma\star{}_pOC)^{\mathbb{Z}/p}_{\delta(\Delta),b}$ is homotopic to $\Sigma_b\circ OC^{\mathbb{Z}/p}$, hence completing the proof of Theorem 4.6. The basic idea is elementary topology: one would like to decompose the diagonal cell $\delta(\Delta_i)\subset S^{\infty}\times S^{\infty}$, which parametrizes the Floer data for the operation $(\Sigma\star{}_pOC)^{\mathbb{Z}/p}_{\delta(\Delta),b}$, into product cells, which parametrize the composition $\Sigma_b\circ OC^{\mathbb{Z}/p}$. On the homology level, this is
\begin{lemma}\cite[Lemma 2.1.]{SW}
Let $\delta:S^{\infty}/\mathbb{Z}/p\rightarrow S^{\infty}/\mathbb{Z}/p\times S^{\infty}/\mathbb{Z}/p$ be the diagonal map and $\delta_*: H_*(S^{\infty}/\mathbb{Z}/p;\mathbb{F}_p)\rightarrow H_*(S^{\infty}/\mathbb{Z}/p;\mathbb{F}_p)^{\otimes 2}$ the induced map on homology, then
\begin{equation}\delta_*\Delta_i=
\begin{cases}
\displaystyle\sum_{i_1+i_2=i} \Delta_{i_1}\otimes \Delta_{i_2},\quad\textrm{if i is odd or}\;p=2,\\
\displaystyle\sum_{\substack{i_1+i_2=i\\ i_k\;\mathrm{even}}} \Delta_{i_1}\otimes \Delta_{i_2},\quad\textrm{if i is even and}\;p>2.\\
\end{cases}
\end{equation}\qed
\end{lemma}
For our purpose, a chain level enhancement of Lemma 4.8 is needed, for the reason that a general element $\mathbf{x}\in {}_pCC_*(\mathcal{F})$ is not $\mathbb{Z}/p$-invariant and thus a cohomology class in $S^{\infty}/\mathbb{Z}/p\times S^{\infty}/\mathbb{Z}/p$ does not give rise to a well-defined operation. The chain level refinement of Lemma 4.8 is the following. 
\begin{lemma}
Fix a ground field $\mathbf{k}$ of characteristic $p$. There exists a sequence of $\mathbf{k}$-coefficient chains $C_i$ of dimension $i$ in $S^{\infty}\times S^{\infty}$, $i\geq 0$, satisfying the recursive property that $\partial C_{i+1}=$
\begin{equation}
\begin{cases}
\displaystyle \delta(\Delta_i)-\sum_{i_1+i_2=i} \Delta_{i_1}\times \Delta_{i_2}-\sum_{\substack{i_1+i_2=i\\ i_1\;\mathrm{odd}}}\Delta_{i_1}\times (\tau-1)\Delta_{i_2}-(\tau\times \tau-1)C_i,\quad\textrm{if i is odd},\\
\displaystyle \delta(\Delta_i)-\sum_{\substack{i_1+i_2=i\\ i_k\;\mathrm{even}}} \Delta_{i_1}\times \Delta_{i_2}-\sum_{\substack{i_1+i_2=i\\ i_1\;\mathrm{odd}}}\sum_{0\leq k<j\leq p-1}\tau^k\Delta_{i_1}\times \tau^j\Delta_{i_2}-(1+\tau\times \tau+\cdots+(\tau\times \tau)^{p-1})C_i,\quad\textrm{if i is even},\\
\end{cases}
\end{equation}
where $\tau\times \tau$ denotes the diagonal action $\tau\times \tau(w_1,w_2):=(\tau(w_1),\tau(w_2))$. 
\end{lemma}
The Proof of Lemma 4.9 will be given in Appendix B. \par\indent
\emph{Proof of Theorem 4.6 given Lemma 4.9}. Let
\begin{equation}
(\Sigma\star{}_pOC)_{C_i,b}:{}_pCC_*(\mathcal{F})\rightarrow CM^*(f)
\end{equation}
be the operation obtained by counting rigid elements of the moduli space (cf. (4.39))
\begin{equation}
\mathcal{M}_{\star}(C_i,\mathcal{R}^1_{k_1,\cdots,k_p}\times \mathcal{S}^1_{p+1},y_{\infty},b,\cdots,b,\mathbf{x}).
\end{equation}
Then, by considering the boundary strata of the $1$-dimensional component of (4.62), and using (4.60), we obtain that
$$
(\Sigma\star{}_pOC)_{C_{i+1},b}\circ d_{{}_pCC^*}\pm d_{CM^*(f)}\circ (\Sigma\star{}_pOC)_{C_{i+1},b}=(\Sigma\star{}_pOC)_{\delta(\Delta_i),b}\pm
$$
\begin{equation}
\begin{cases}
\displaystyle\sum_{i_1+i_2=i}(\Sigma\star{}_pOC)_{\Delta_{i_1}\times\Delta_{i_2},b}\pm \sum_{\substack{i_1+i_2=i\\i_1\;\textrm{odd}}}(\Sigma\star{}_pOC)_{\Delta_{i_1}\times (\tau-1)\Delta_{i_2},b}\pm(\Sigma\star{}_pOC)_{(\tau\times \tau-1)C_i,b},\;\textrm{i odd}\\
\displaystyle \sum_{\substack{i_1+i_2=i\\i_k\;\textrm{even}}}(\Sigma\star{}_pOC)_{\Delta_{i_1}\times\Delta_{i_2},b}\pm \sum_{\substack{i_1+i_2=i\\i_1\;\textrm{odd}}}\sum_{0\leq k<j\leq p-1}(\Sigma\star{}_pOC)_{\tau^k\Delta_{i_1}\times \tau^j\Delta_{i_2},b}\pm(\Sigma\star{}_pOC)_{(1+\tau\times \tau+\cdots+(\tau\times \tau)^{p-1})C_i,b},\;\textrm{i even}.
\end{cases}
\end{equation}
Since
\begin{equation}
(\Sigma\star{}_pOC)_{\Delta_{i_1}\times \Delta_{i_2},b}=\Sigma^{i_2}_b\circ {}_pOC^{i_1}
\end{equation}
and the Floer data are $\mathbb{Z}/p$-equivariant, we can simplify (4.63) as
$$
(\Sigma\star{}_pOC)_{C_{i+1},b}\circ d_{{}_pCC^*}\pm d_{CM^*(f)}\circ (\Sigma\star{}_pOC)_{C_{i+1},b}=(\Sigma\star{}_pOC)_{\delta(\Delta_i),b}\pm
$$
\begin{equation}
\begin{cases}
\displaystyle\sum_{i_1+i_2=i}\Sigma^{i_2}_b\circ{}_pOC^{i_1} \pm(\Sigma\star{}_pOC)_{C_i,b}\circ(\tau-1),\;\textrm{i odd}\\
\displaystyle \sum_{\substack{i_1+i_2=i\\i_k\;\textrm{even}}}\Sigma_b^{i_2}\circ{}_pOC^{i_1}\pm \sum_{\substack{i_1+i_2=i\\i_1\;\textrm{odd}}}\Sigma_b^{i_2}\circ {}_pOC^{i_1}\circ(\tau-1)^{p-2}\pm(\Sigma\star{}_pOC)_{C_i,b}\circ(1+\tau+\cdots+\tau^{p-1}),\;\textrm{i even}.
\end{cases}
\end{equation}
We let $\tilde{\mathcal{H}}^i_b:=(\Sigma\star{}_pOC)_{C_{i+1},b}$ (notice the shift in $i$) and define $\tilde{\mathcal{H}}^{\mathbb{Z}/p}_b:{}_pCC^{\mathbb{Z}/p}_*(\mathcal{F})\rightarrow  CM^*(f)[[t,\theta]]$ by
\begin{equation}
\begin{cases}
\mathbf{x}\mapsto \sum_k\big(\tilde{\mathcal{H}}^{2k}_b(\mathbf{x})+(-1)^{\|x\|} \tilde{\mathcal{H}}^{2k+1}_b(\mathbf{x})\theta\big)t^k\\
\mathbf{x}\theta\mapsto \sum_k\big(\tilde{\mathcal{H}}^{2k}_b(\mathbf{x})\theta+(-1)^{\|x\|} \tilde{\mathcal{H}}^{2k+1}_b((\tau-1)^{p-2}\mathbf{x})t\big)t^k.
\end{cases}
\end{equation}
By (4.65), 
\begin{equation}
\tilde{\mathcal{H}}_b^{\mathbb{Z}/p}\circ d_{{}_pCC_*^{\mathbb{Z}/p}}\pm d_{CM^*(f)}\circ \tilde{\mathcal{H}}_b^{\mathbb{Z}/p}\pm (\Sigma\star{}_pOC)_{\delta(\Delta),b}^{\mathbb{Z}/p}\pm \Sigma_b\circ {}_pOC^{\mathbb{Z}/p}
\end{equation}
is the $t$-linear map
\begin{equation}
\begin{cases}
\displaystyle\mathbf{x}\mapsto\sum_k\sum_{\substack{i_1+i_2=2k\\i_2\;\textrm{odd}}}\Sigma^{i_2}_b\Big({}_pOC^{i_1}((\tau-1)^{p-2}x)\Big)t^k\\
\displaystyle\mathbf{x}\theta\mapsto \sum_k\sum_{\substack{i_1+i_2=2k+1\\ i_2\;\textrm{odd}}}\Sigma^{i_2}_b\Big({}_pOC^{i_1}((\tau-1)^{p-2})x\Big)t^{k+1},
\end{cases}
\end{equation}
which is the composition of the chain map
\begin{equation}
\begin{cases}
\displaystyle\mathbf{x}\mapsto\sum_k\sum_{\substack{i_1+i_2=2k\\i_2\;\textrm{odd}}}\Sigma^{i_2}_b\Big({}_pOC^{i_1}(x)\Big)t^k\\
\displaystyle\mathbf{x}\theta\mapsto \sum_k\sum_{\substack{i_1+i_2=2k+1\\ i_2\;\textrm{odd}}}\Sigma^{i_2}_b\Big({}_pOC^{i_1}(x)\Big)t^{k+1}.
\end{cases}
\end{equation}
with $x\mapsto (\tau-1)^{p-2}x$. Since $p>2$, Lemma 2.9 implies that $(\tau-1)^{p-2}$ is nullhomotopic as an endomorphism of $CC_*^{\mathbb{Z}/p}(\mathcal{F})$. Therefore, (4.68) is nullhomotopic and the proof is complete. \qed
\begin{rmk}
When $p=2$, the analogue of (4.60) is
$$\partial C_{i+1}=\delta(\Delta_i)+\sum_{i_1+i_2=i}\Delta_{i_1}\times\Delta_{i_2}+\sum_{\substack{i_1+i_2=i\\i_1\;\mathrm{odd}}}\Delta_{i_1}\times(1+\tau)\Delta_{i_2}+(1+\tau\times \tau)C_i.
$$
However, we do not discuss this case further in the current paper.
\end{rmk} 

\renewcommand{\theequation}{5.\arabic{equation}}
\setcounter{equation}{0}
\section{Case study: the intersection of two quadrics in $\mathbb{C}P^5$}
Let $X$ be the intersection of two quadrics in $\mathbb{C}P^5$, which is a monotone symplectic manifold of real dimension $6$. We know that $H^{2i}(X,\mathbb{Z})=\mathbb{Z}$ and $H^3(X,\mathbb{Z})=\mathbb{Z}^4$. Let $1,h=h_2, h_4,h_6$ be the even generators. In \cite{Don}, it was shown that there is a canonical isomorphism $\eta: H_1(\Sigma_2)\rightarrow H^3(X)$, where $\Sigma_2$ is a genus two surface. In \cite[Example 6.8]{SW}, the quantum Steenrod powers involving even cohomology classes have been computed using covariant constancy. The goal of this section is to compute $QSt$ for an odd cohomology class, conditional on certain assumptions on the arithmetic Fukaya category of $X$. \par\indent
Recall we have a decomposition of the Fukaya category
\begin{equation}
\mathrm{Fuk}(X,\mathbb{C})=\mathrm{Fuk}(X,\mathbb{C})_{-8}\oplus \mathrm{Fuk}(X,\mathbb{C})_0\oplus \mathrm{Fuk}(X,\mathbb{C})_{8}.
\end{equation}
A celebrated result of \cite{Smi} (cf. Theorem 1.1 loc.cit.) shows that there is a quasi-equivalence $F: Tw^{\pi}\mathrm{Fuk}(X,\mathbb{C})_0\simeq Tw^{\pi}\mathrm{Fuk}(\Sigma_2,\mathbb{C})$, where $Tw^{\pi}$ denotes first taking twisted complexes and then idempotent completion. Since applying $Tw^{\pi}$ does not change the Hochschild invariants, by abuse of notation, we also denote by $F$ the induced isomorphism on Hochschild (co)homology/cyclic homology. \par\indent
In this example, we take $R=\mathbb{Z}[\frac{1}{2}]$. Condition A1) and A2) follow from the computations in Lemma 5.4 below. Moreover, condition A3) follows from \cite[Corollary 4.16]{Smi}. We remark that even though loc.cit. uses complex coefficients, the results therein regarding generation and the open-closed map being an isomorphism (specifically Corollary 3.11, Lemma 4.15 and Corollary 4.16 in loc.cit.) hold more generally over $R$. An interesting question is whether Smith's equivalence holds over any field of characteristic not $2$. 
\begin{conj}
Over a field $K$ of characteristic not $2$, there is an $A_{\infty}$-quasi-equivalence
\begin{equation}
F: Tw^{\pi}\mathrm{Fuk}(X)_0\simeq Tw^{\pi}\mathrm{Fuk}(\Sigma_2).  
\end{equation}
\end{conj}
Unfortunately, Conjecture 5.1 cannot hold at this level of generality. In fact, we will show in section 5.1 that Smith's equivalence imposes the following simple arithmetic property on $K$.
\begin{prop}
If Conjecture 5.1 holds, then $K$ must contain a square root of $-1$.    
\end{prop}
Now we go to positive characteristics. Fix an odd prime $p$ and a field $\mathbf{k}$ of characteristic $p$. Given Proposition 5.2, it is natural to ask whether Smith's quasi-equivalence holds over any $\mathbf{k}$ that contains a $\sqrt{-1}$. While we do not obtain such a result in this paper, we show in section 5.2 that Smith's equivalence over $\mathbf{k}$ determines the quantum Steenrod powers of all odd cohomology classes of $X$.
\begin{prop}
Let $\mathbf{k}$ be a field of odd characteristic $p$ that contains a $\sqrt{-1}$. Assuming Conjecture 5.1 holds over $\mathbf{k}$, then
\begin{equation}
QSt^X(\eta(\gamma_i))=(-1)^{\frac{p-1}{2}}(\frac{p-1}{2}!)t^{\frac{p-1}{2}}(1+\frac{1}{256}t^2+\frac{81}{262144}t^4+\cdots)\eta(\gamma_i).
\end{equation}
\end{prop}
We remark that the $t^{\frac{p-1}{2}}$-coefficient of $QSt^X(\eta(\gamma_i))$ in (5.3) is $(-1)^{\frac{p-1}{2}}(\frac{p-1}{2}!)$, which implies that $Q\Xi_{X}=(-1)^{\frac{p-1}{2}}\mathrm{id}$ when $\mathbf{k}=\mathbb{F}_p$, $p=1$ mod $4$, cf. \cite[Definition 1.8]{Sei3}. This recovers the computation of $Q\Xi_{X}$ in \cite[Example 9.8]{Sei3}.

\subsection{Quantum cohomology and open-closed maps}
\begin{lemma} We record some computations from \cite{Don}.
\begin{enumerate}[label=\arabic*)]
    \item The ordinary cup products on $H^*(X,\mathbb{Z})$ is given by
\begin{equation}
h\cup h=4h_4,\;\; h\cup h_4=h_6,\;\;  \eta(\gamma_1)\cup\eta(\gamma_2)=(\gamma_1\cdot\gamma_2)h_6,  
\end{equation}
where $(\;\cdot\;)$ denotes the intersection pairing on $H_1(\Sigma_2)$.
\item The nontrivial quantum products on $QH^*(X,\mathbb{Z})$ is given by
\begin{equation}
h\star h=4h_4+4,\;\; h\star h_4=h_6+2h,\;\;h\star h_6=4h_4+4,\;\;h_4\star h_4=2h_4+3    
\end{equation}
and
\begin{equation}
\eta(\gamma_i)\star\eta(\gamma_j)=(\gamma_i\cdot\gamma_j)\frac{h^3-16h}{4},
\end{equation}
where we use powers, e.g. $h^3$, to denote powers of quantum products. 
\item The generalized eigenvalues of $c_1\star=2h\star$ are $\{-8,0,8\}$. Over $R=\mathbb{Z}[\frac{1}{2}]$, there is a generalized eigenspace decomposition
\begin{align}
QH^*(X)&=QH^*(X)_{-8}\oplus QH^*(X)_0\oplus QH^*(X)_8\nonumber\\
&=\langle \frac{-h_6+4h_4-3h_2+4}{4}\rangle\oplus\langle 1-\frac{h^2}{16}, \eta(H_1(\Sigma_2)),h^3-16h\rangle \oplus \langle \frac{h^3+4h^2}{16}\rangle.
\end{align}
\item Let $v=1-\frac{h^2}{16}$ and $w=h^3-16h$. Then, $v$ is the $0$-idempotent, i.e. the unit element of $QH^*(X)_0$. Moreover, if we let $\langle\cdot,\cdot\rangle_X=\int_X(\cdot\cup\cdot)$ denote the Poincare pairing on $X$, then
\begin{equation}
\langle v,v\rangle_X=0,\;\;\langle v,w\rangle_X=4.   
\end{equation}
\end{enumerate}\qed
\end{lemma}

\begin{lemma}
Suppose Conjecture 5.1 holds over $K$. We denote $\Xi=CO_{\Sigma_2}^{-1}\circ F\circ CO_X$ and $\Theta=OC_{\Sigma_2}\circ F\circ OC_X^{-1}$, where $OC_{X}$ and $CO_{X}$ are understood to be restricted to $QH^*(X,K)_0$ and $HH^*(\mathrm{Fuk}(X,K)_0)$, respectively. Then, we have
\begin{enumerate}[label=\arabic*)]
    \item $K$ must contain a square root of $-1$, denoted $\sqrt{-1}$, and
    \begin{equation}
      \Theta(v)=\frac{\sqrt{-1}\epsilon}{2},\;\;\Theta(w)=8\sqrt{-1}\epsilon H  
    \end{equation}
    where $H$ is the standard generator of $H^2(\Sigma_2,\mathbb{Z})$ (viewed as an element of $H^2(\Sigma_2,K)$), and $\epsilon\in\{-1,1\}$ is a sign we leave undetermined for now. 
    \item \begin{equation}
      \Theta(\eta(\gamma))=\frac{\sqrt{-1}\epsilon}{2}\Xi(\eta(\gamma)),\;\;\gamma\in H_1(\Sigma_2).  
    \end{equation}
\end{enumerate}
\end{lemma}
\emph{Proof.} Recall that $v=1-\frac{h^2}{16}$ and $w=h^3-16h$. Since $\Xi$ is a unital ring map, we have $\Xi(v)=1$. We now determine $\Xi(w)$. The crucial property is that $\Xi(c_1(X)|_{QH^*(X)_0})=c_1(\Sigma_2)$ (this follows from the fact that the cyclic open-closed map intertwines connections). As a result, $\Xi(w)=-8\Xi(c_1(X)\star v)=-8c_1(\Sigma_2)\cup \Xi(v)=16H$ (note cup product and quantum product on $\Sigma_2$ agree).    \par\indent
With the constant $s$ yet to be determined, we now proceed to compute $\Theta$. Let $\Theta(v)=a+bH$ and $\Theta(w)=c+dH$. Since $\Theta$ is a map of $QH^*(X)_0$-modules via $\Xi$, $0=\Theta(w\star w)=\Xi(w)\cup\Theta(w)=16H\cup (c+dH)=16c H$. Thus we must have $c=0$, whence $d\neq 0$. Similarly, $dH=\Theta(w)=\Theta(w\star v)=\Xi(w)\cup\Theta(v)=16aH$, and thus $a\neq 0$ and $d=16a$. Since $\Theta$ preserves Poincare pairing up to a minus sign, by Lemma 5.4 4), one has $0=-\langle v,v\rangle_X=\langle \Theta(v),\Theta(v)\rangle_{\Sigma_2}=2ab$, which implies that $b=0$. Similarly, $-4=-\langle v,w\rangle_X=\langle \Theta(v),\Theta(w)\rangle_{\Sigma_2}=ad$.\par\indent
Now, note that $d=16a,ad=-4$ is solvable only when the coefficient field contains a root of $-1$, in which case $a=\frac{\sqrt{-1}\epsilon}{2}, d=8\sqrt{-1}\epsilon$, for some $\epsilon\in\{-1,1\}$. The ambiguity of signs is canceled out in the computations below. For the moment, we do not determine $\Xi$ or $\Theta$ restricted to the odd dimensional cohomology, but simply note that $\Theta(\eta(\gamma_i))=\Xi(\eta(\gamma_i))\cup \Theta(v)=\frac{\sqrt{-1}\epsilon}{2}\Xi(\eta(\gamma_i))$, again by the fact that $\Theta$ is a map of $QH^*(X)_0$-modules.\qed\par\indent
In particular, Lemma 5.5 1) implies Proposition 5.2, and in particular it implies that Smith's equivalence cannot hold over $\mathbf{k}=\mathbb{F}_p$ for some prime $p\equiv 3$ mod $4$. This agrees with an observation of Seidel in \cite[Example 9.8]{Sei3} concerning formal groups. \par\indent 
An important computation that will be used later is the determination of $\Theta$ and $\Xi$ restricted $H^3(X,K)$, up to the effect of an autoequivalence of $\mathrm{Fuk}(\Sigma_2)$. The precise statements are given in Lemma 5.6. We first fix some notations.\par\indent
Recall there is a \emph{Chern character} map, cf. \cite[3.1]{Shk},
\begin{equation}
ch: K_0(\mathrm{Fuk}(X,K)_0)\rightarrow    HH_0(\mathrm{Fuk}(X,K)_0)
\end{equation}
defined by $ch(L):=[\mathrm{id}_L]$. Where $[\mathrm{id}_L]$ denotes the cohomological unit of $L$ viewed as a length $0$ Hochschild homology class. Shklyarov showed that (cf. Proposition 4.4 loc.cit.)
\begin{equation}
\langle ch(L), ch(L')\rangle_{Sh} =\chi(L,L'),    
\end{equation}
where $\chi(L,L')=\sum_i (-1)^i \dim HF^i(L,L')$ denotes the Euler pairing and $\langle-,-\rangle_{Sh}$ denotes the Shklyarov pairing on Hochschild homology. We say that an object $L$ has \emph{integral Chern character} if $OC(ch(L))\in \mathrm{im}(H^3(X,\mathbb{Z})\rightarrow H^3(X,K))$.
\begin{lemma}
\begin{enumerate}[label=\arabic*)]
    \item Any object $L\in \mathrm{Fuk}(X,K)_0$ has integral Chern character.
    \item Assume we are in the situation of Lemma 5.5. Let $\{\gamma_1,\gamma_2,\gamma_3,\gamma_4\}$ be a standard symplectic basis of $H_1(\Sigma_2,\mathbb{Z})$ with respect to the intersection pairing.  Then with respect to the basis $\{\eta(\gamma_i)\}_{i=1}^4$ of $H^3(X,K)$ and $\{\mathrm{PD}(\gamma_i)\}_{i=1}^4$ of $H^1(\Sigma_2,K)$, 
    \begin{equation}\Theta|_{H^3(X,K)}: H^3(X,K)\rightarrow H^1(\Sigma_2,K)
    \end{equation}
    has integral coefficients.  
    \item Assume we are in the situation of Lemma 5.5. There exists an autoequivalence $G\in \mathrm{Auteq}(\mathrm{Fuk}(\Sigma_2,K))$ such that if we replace $F$ by $G\circ F$, then
    \begin{equation}
        \Theta=\begin{bmatrix}
          I_2&\mathbf{0}\\
          \mathbf{0} &-I_2
        \end{bmatrix}
    \end{equation} 
    with respect to the basis $\{\eta(\gamma_i)\}_{i=1}^4$ and $\{\mathrm{PD}(\gamma_i)\}_{i=1}^4$.
\end{enumerate}
\end{lemma}
\emph{Proof}. 1) Take an integral basis $\{\gamma_i\}_{i=1}^4$ of $H_1(\Sigma_2,\mathbb{Z})$. These give rise to vanishing cycles $V_{\gamma_i}\subset X, 1\leq i\leq 4$, cf. \cite[Section 4.2]{Smi}. As $V_{\gamma_i}$ bounds no Maslov index $2$ disk, $OC(ch(V_{\gamma_i}))$ is the image of $\mathrm{PD}([V_{\gamma_i}])$ under $H^3(X,\mathbb{Z})\rightarrow H^3(X,K)$. By (5.12), and the fact that $OC$ intertwines the Shklyarov pairing with the Poincare pairing up to a sign, we conclude that for any $L\in \mathrm{Fuk}(X,K)_0$, the Poincare pairing of $OC(ch(L))$ with $\mathrm{PD}([V_{\gamma_i}])$ is integral, for any $1\leq i\leq 4$. Since $[V_{\gamma_i}]$ form a basis for $H_3(X,\mathbb{Z})$, $OC(ch(L))$ must be integral. \par\indent
2) By \cite[Lemma 2.19]{AS}, any object in $\mathrm{Fuk}(\Sigma_2,K)$ has integral Chern character. In particular, for $1\leq i\leq 4$, $OC_{\Sigma_2}(HH_*(F)(ch_{V_{\gamma_i}}))=OC_{\Sigma_2}(ch_{F(V_{\gamma_i})})$ has integral coefficients when expanded with respect to the basis $\{\mathrm{PD}(\gamma_i)\}_{i=1}^4$. Since $OC_X(ch(V_{\gamma_i}))=\mathrm{PD}([V_{\gamma_i}]), 1\leq i\leq 4$ form an integral basis of $H^3(X)$, 2) follows. \par\indent
3) Since $OC$ intertwines pairings up to $(-1)^{\frac{n(n+1)}{2}}$ (where $n$ is the complex dimension), we conclude that with respect to the basis $\{\eta(\gamma_i)\}_{i=1}^4$ and $\{\mathrm{PD}(\gamma_i)\}_{i=1}^4$, the matrix form of $\Theta$ is in $\mathrm{Sp}_4^-(\mathbb{Z})$, the set of anti-symplectic matrices. Therefore, to prove 3), it suffices to show that the natural map
\begin{equation}
\mathrm{Auteq}(\mathrm{Fuk}(\Sigma_2,K))\rightarrow \mathrm{Aut}(HH_0(\mathrm{Fuk}(\Sigma_2,K)),\langle-,-\rangle_{Sh})\xrightarrow[OC]{\simeq}   \mathrm{Aut}(H^1(\Sigma_2,K),\langle-,-\rangle_{\Sigma_2})=\mathrm{Sp}_4(K)
\end{equation}
surjects onto the image of $\mathrm{Sp}_4(\mathbb{Z})\rightarrow \mathrm{Sp}_4(K)$. This can be seen by the following argument. \par\indent
Let $\Gamma(\Sigma_2)$ denote the symplectic mapping class group of $\Sigma_2$. Then, there is a natural map \begin{equation}\Gamma(\Sigma_2)\rightarrow \mathrm{Aut}(H^1(\Sigma_2,K),\langle-,-\rangle_{\Sigma_2})=\mathrm{Sp}_4(K)
\end{equation}
whose image is $\mathrm{Sp}_4(\mathbb{Z})$; this is because $\mathrm{Sp}_4(\mathbb{Z})$ is generated by algebraic Dehn twists, which can be lifted to symplectic Dehn twists along simple closed curves. Moreover, (5.16) factors as a homomorphism (depending on a choice of \emph{balancing}, cf. \cite[Section 2.6]{AS}) $\Gamma(\Sigma_2)\rightarrow \mathrm{Auteq}(\mathrm{Fuk}(\Sigma_2,K))$ composed with (5.15). Hence the image of (5.15) must also contain $\mathrm{Sp}_4(\mathbb{Z})$.\qed 
\subsection{A computation of $QSt^X(\eta(\gamma_i))$}
For the rest of this subsection, fix a field $\mathbf{k}$ of odd characteristic $p$ containing a $\sqrt{-1}$. Under the same assumptions as Lemma 5.5 (with $K=\mathbf{k}$), we now outline a way to compute the quantum Steenrod operations on $H^3(X,\mathbf{k})$ using the categorical interpretation of $QSt$ proved in Theorem 1.2. First,
\begin{align}
&QSt^X(\eta(\gamma_i))=Q\Sigma^X_{\eta(\gamma_i)}(1)\\
&=Q\Sigma^X_{\eta(\gamma_i)}(pr_0(1))\qquad (pr_0(1)\;\;\textrm{is the projection of $1$ onto $QH^*(X)[[t,\theta]]_0$})\\
&=OC_X^{\mathbb{Z}/p}\Big(CO_X(\eta(\gamma_i))\cap^{\mathbb{Z}/p} OC_X^{\mathbb{Z}/p}(pr_0(1))\Big)\qquad \textrm{(by Theorem 1.1)}\\
&=OC_X^{\mathbb{Z}/p}\circ F^{-1}\circ (OC^{\mathbb{Z}/p}_{\Sigma_2})^{-1}\Big(Q\Sigma^{\Sigma_2}_{\Xi(\eta(\gamma_i))} (OC^{\mathbb{Z}/p}_{\Sigma_2}\circ F\circ (OC^{\mathbb{Z}/p}_X)^{-1}pr_0(1))\Big).
\end{align}
At this point, we make a few observations. Firstly, since $\Sigma_2$ has no genus $0$ Gromov-Witten invariants, the quantum Steenrod operations is classical: 
\begin{equation}
Q\Sigma^{\Sigma_2}_{y}(c)=St^{\Sigma_2}(y)\cup c=(\frac{p-1}{2}!)t^{\frac{p-1}{2}}y^{(1)}\cup c,     
\end{equation}
for $y\in H^1(\Sigma_2,\mathbf{k})$. Here, $(-)^{(1)}$ denotes the relative Frobenius along $H^1(\Sigma_2,\mathbb{F}_p)\rightarrow H^1(\Sigma_2,\mathbf{k})$. Concretely, expanding $y$ in terms of a $\mathbf{k}$-basis coming from $H^1(\Sigma_2,\mathbb{F}_p)$, then $y^{(1)}$ is obtained by raising each coefficient to the $p$-th power. Secondly, by Theorem 1.3, we have
\begin{equation}
OC_X^{\mathbb{Z}/p}\circ F^{-1}\circ (OC^{\mathbb{Z}/p}_{\Sigma_2})^{-1}=OC_X^{S^1}\circ F^{-1}\circ (OC^{S^1}_{\Sigma_2})^{-1}(\otimes_{\mathbf{k}[[t]]}\mathbf{k}[[t,\theta]]).
\end{equation}
Hence we can write $(OC^{\mathbb{Z}/p}_{\Sigma_2}\circ F\circ (OC^{\mathbb{Z}/p}_X)^{-1}pr_0(1))=\sum_{n\geq 0}(c_n+d_nH)t^n$ for some coefficients $c_n,d_n$. Then,
\begin{equation}
Q\Sigma^{\Sigma_2}_{\Xi(\eta(\gamma_i))} (OC^{\mathbb{Z}/p}_{\Sigma_2}\circ F\circ (OC^{\mathbb{Z}/p}_X)^{-1}pr_0(1))=(\frac{p-1}{2}!)t^{\frac{p-1}{2}}\sum_{n\geq 0} c_n\Xi(\eta(\gamma_i))^{(1)}t^n.
\end{equation}
Finally, we observe that the $t$-connection restricted to the summand $\eta(H_1(\Sigma_2))[[t]]$ is trivial (as both $c_1\star$ and $\mu$ are trivial when applied to the image of $\eta$). As a consequence, $OC_X^{S^1}\circ F^{-1}\circ (OC^{S^1}_{\Sigma_2})^{-1}|_{H^1(\Sigma_2)[[t]]}=\Theta^{-1}[[t]]$ because a regular connection over $\mathbb{Z}$ has no nontrivial automorphism whose constant term is the identity. By Lemma 5.6 2), 
$\Theta(\eta(\gamma_i))\in \mathrm{im}(H^1(\Sigma_2,\mathbb{F}_p)\rightarrow H^1(\Sigma_2,\mathbf{k}))$. In particular, $\Theta^{(1)}=\Theta$ since the $p$-th power of an element in $\mathbb{F}_p$ is itself. Therefore, $QSt^X(\eta(\gamma_i))$ is equal to
\begin{align}
(\frac{p-1}{2}!)t^{\frac{p-1}{2}}\sum_{n\geq 0} c_n\Theta^{-1}\Xi(\eta(\gamma_i))^{(1)}t^n&\overset{(5.10)}{=}(\frac{p-1}{2}!)t^{\frac{p-1}{2}}\sum_{n\geq 0} c_n\Theta^{-1}(-2\sqrt{-1}\epsilon\Theta(\eta(\gamma_i)))^{(1)}t^n\nonumber\\
&=(-2\sqrt{-1}\epsilon)^p(\frac{p-1}{2}!)t^{\frac{p-1}{2}}\Theta^{-1}\Theta^{(1)}\eta(\gamma_i)\sum_{n\geq 0} c_n t^n\nonumber\\
&=(-2\sqrt{-1}\epsilon)(-1)^{\frac{p-1}{2}}(\frac{p-1}{2}!)t^{\frac{p-1}{2}} \eta(\gamma_i)\sum_{n\geq 0} c_n t^n.
\end{align}
Therefore, to finish the computation, it suffices to determine $c_n$. These coefficients can be extracted from the so-called \emph{$R$-matrix}, cf. Lemma 5.7; it turns out the $R$-matrix will also depend on the sign $\epsilon$, which cancels out the $\epsilon$ in (5.24).  
\begin{lemma} \cite[Lemma B.1]{Hug}
Any isomorphism $\varphi: QH^*(X,R)\cong R\oplus QH^*(\Sigma_2,R)\oplus R$ that intertwines the operations $c_1(X)\star$ and $(8,c_1(\Sigma_2)\star,-8)$ can be uniquely extended (meaning is the $t=0$ term of a unique $k[[t]]$-linear map) to an isomorphism 
\begin{equation}
\tilde{\varphi}: QH^*(X,R)[[t]]\cong \mathcal{E}^{\frac{8}{t}}\oplus QH^*(\Sigma_2,R)[[t]]\oplus \mathcal{E}^{\frac{-8}{t}}
\end{equation}
that intertwines the $t$-connections. If one chooses basis for both sides of (5.25), the matrix form of $\tilde{\varphi}$ is called the \emph{$R$-matrix} with constant term $\varphi$. \qed
\end{lemma}
We remark that 1) \cite[Lemma B.1]{Hug} was originally proved over $\mathbb{C}$, but the computations involved make clear they work over any coefficient ring with $2$ inverted and 2) it suffices to determine $\tilde{\varphi}$ on the even degree part of $QH^*(X,R)[[t]]$, this is because on the odd degree part $\eta(H_1(\Sigma_2))[[t]]$, the $t$-connection is trivial. \par\indent
Consider the isomorphism 
\begin{equation}
(\frac{-h_6+4h_4-3h_2+4}{4}, \Theta^{-1}, \frac{h^3+4h^2}{16}): R\oplus QH^*(\Sigma_2,R)\oplus R\rightarrow QH^*(X,R).
\end{equation}
When restricted to the even degree part, with respect to the basis 
\begin{equation}
\{\frac{-h_6+4h_4-3h_2+4}{4}, w, v, \frac{h^3+4h^2}{16}\}
\end{equation}
of $QH^{even}(X,R)$ and the basis
\begin{equation}
\{(\sqrt{-1}\epsilon,0,0),(0,-2H,0),(0,1,0),(0,0,\sqrt{-1}\epsilon)\}
\end{equation}
of $R\oplus QH^{even}(\Sigma_2,R)\oplus R$, the inverse of (5.26) is given by the matrix
\begin{equation}
 \begin{bmatrix}
 1 &0&0&0\\
 0&-4&0&0\\
 0&0&\frac{1}{2}&0\\
 0&0&0&1
 \end{bmatrix}\sqrt{-1}\epsilon.
\end{equation}
Since the cyclic open-closed map preserves connection, by the uniqueness property of $R$-matrix in Lemma 5.7, one concludes that when restricted to the even degree part, $OC_X^{S^1}\circ F^{-1}\circ (OC_{\Sigma_2}^{S^1})^{-1}$ agrees with the $R$-matrix whose constant term is (5.29). The next lemma follows from using the computation in \cite[Appendix B.1]{Hug}. 
\begin{lemma}
With notations as above, the $R$ matrix with constant term (5.29) is given by
\begin{equation}
 \begin{bmatrix}
 1 &0&0&0\\
 0&-4&0&0\\
 0&0&\frac{1}{2}&0\\
 0&0&0&1
 \end{bmatrix}\sqrt{-1}\epsilon-
 \begin{bmatrix}
 -\frac{7}{64}& \frac{1}{2}&   \frac{1}{32}&    0\\
 \frac{1}{16}&   0& -\frac{5}{512}& \frac{1}{16}\\
 -\frac{1}{8}&  -2&      0&  \frac{1}{8}\\
  0& \frac{1}{2}&  -\frac{1}{32}& \frac{7}{64}
 \end{bmatrix}\sqrt{-1}\epsilon t+
 \begin{bmatrix}
 -\frac{15}{8192}&-\frac{7}{128}& \frac{9}{2048}&   \frac{1}{1024}\\
 -\frac{27}{2048}&  \frac{7}{128}&      0&  \frac{27}{2048}\\
  \frac{1}{32}&     0& -\frac{3}{512}&     \frac{1}{32}\\
 \frac{1}{1024}&  \frac{7}{128}& \frac{9}{2048}& -\frac{15}{8192}
 \end{bmatrix}\sqrt{-1}\epsilon t^2\nonumber
\end{equation}
\begin{equation}
-\begin{bmatrix}
-\frac{389}{524288}& -\frac{135}{16384}&  \frac{393}{262144}&   \frac{65}{65536}\\
     \frac{33}{8192}&         0& -\frac{435}{1048576}&    \frac{33}{8192}\\
 -\frac{9}{2048}&      \frac{3}{128}&            0&     \frac{9}{2048}\\
 -\frac{65}{65536}& -\frac{135}{16384}&  -\frac{393}{262144}& \frac{389}{524288}
\end{bmatrix} \sqrt{-1}\epsilon t^3+
\begin{bmatrix}
 -\frac{38421}{134217728}& -\frac{3069}{1048576}& \frac{11907}{16777216}&     \frac{3537}{8388608}\\
    -\frac{7533}{4194304}&    \frac{999}{262144}&              0&     \frac{7533}{4194304}\\
           \frac{9}{8192}&             0&     -\frac{63}{524288}&           \frac{9}{8192}\\
     \frac{3537}{8388608}&  \frac{3069}{1048576}& \frac{11907}{16777216}& -\frac{38421}{134217728}  
\end{bmatrix}\sqrt{-1}\epsilon t^4\pm\cdots
\end{equation}\qed
\end{lemma}
\emph{Proof of Proposition 5.3}. By (5.24), it suffices to compute the coefficients $c_n$, which we obtain via the following steps. First, express $1$ in terms of our chosen basis (5.27) (viewed as a basis of $QH^{even}(X,R)[[t]]$ over $k[[t]]$)
\begin{equation}
1=\frac{1}{8}\cdot\frac{-h_6+4h_4-3h_2+4}{4}+0\cdot w+1\cdot v+  \frac{1}{8}\cdot \frac{h^3+4h^2}{16}.
\end{equation}
Then, apply the $R$-matrix of Lemma 5.6 to the resulting vector $(\frac{1}{8},0,1,\frac{1}{8})$. Next, apply the matrix 
\begin{equation}
\begin{bmatrix}
0&1&0&0\\
0&0&1&0
\end{bmatrix},  
\end{equation}
which under the basis (5.28) corresponds to the projection 
\begin{equation}\mathcal{E}^{\frac{8}{u}}\oplus QH^{even}(\Sigma_2,R)[[t]]\oplus \mathcal{E}^{-\frac{8}{u}}\rightarrow QH^{even}(\Sigma_2,R)[[t]].
\end{equation}
These steps combined have the effect of computing 
\begin{equation}
OC^{\mathbb{Z}/p}_{\Sigma_2}\circ F\circ (OC^{\mathbb{Z}/p}_X)^{-1}(pr_0(1))
\end{equation}
in terms of the basis $\{-2H,1\}$ of $QH^{even}(\Sigma_2,R)[[t]]$, and $\sum_{n\geq 0} c_nt^n$ is simply the coefficient in front of $1$. Unfolding the computations we obtain Proposition 5.3. \qed\par\indent
Things to note about the coefficients in (5.3):
\begin{enumerate}[label=\arabic*)]
    \item All the denominators are powers of 2, hence the expression in fact makes sense over $\mathbb{Z}[\frac{1}{2}]$. 
    \item For fixed $p$, all high enough coefficients have numerators divisible by $p$; in principle, one can compute exactly how large that is. Therefore when reduced mod $p$, this gives a polynomial expression in $t$.
    \item If we add back the Novikov variable $q$ (which we have set to $1$ throughout the paper) of degree $2$, then for degree reasons (since $QSt^X(\eta(\gamma_i))$ has total degree $3p$) (5.3) involve no powers of $q$ greater than $p$. Geometrically, this says that $p$-fold cover curves do not contribute to $QSt^X(\eta(\gamma_i))$.
\end{enumerate}

\renewcommand{\theequation}{6.\arabic{equation}}
\setcounter{equation}{0}
\section{A conjectural $B$-side formula}
We expect that Theorem 1.2 has applications to the study of mirror symmetry over integers or a number field, such as studied in the work of \cite{LPe}, \cite{LPo} and \cite{GHHPS}. The role of quantum Steenrod operations in arithmetic mirror symmetry was envisioned by \cite{Sei3}, which explored a parallelism between the formal group law on the quantum cohomology of a closed monotone symplectic manifold and the formal group law of its Fukaya category. The former is closely related to a certain coefficient of the quantum Steenrod powers, cf. Theorem 1.9 and Lemma 2.10 in loc.cit., and gives rise to potential mirror constructions that have interesting arithmetic meanings, cf. Conjecture 9.12 in loc.cit.\par\indent
Theorem 1.2 of the current paper gives the first explicit identification of quantum Steenrod operation as a Fukaya categorical invariant. Given this, it is natural to ask that in a hypothetical arithmetic mirror situation, what are the operations on the B-side that mirror the quantum Steenrod operations? In particular, these are operations of (a suitable version of) polyvector fields on the de-Rham cohomology. We expect that the mirror operations are related to `Frobenius $p$-th power maps', as well as the `$p$-th power contraction' $i^{[p]}$ studied for instance in \cite[Lemma 2.1]{BK}. The latter is closely related to the restricted Lie structure on polyvector fields of a characteristic $p$ scheme and play an important role in the study of deformation quantization in positive characteristics. More precisely, we expect the following correspondence. \par\indent
Let $(A,f)$ be a (formal) smooth commutative algebra of dimension $d<p$ over $\mathbf{k}$ equipped with a superpotential $f\in A$. Loosely speaking, one expects Hochschild-Kostant-Rosenberg type quasi-isomorphisms
\begin{equation}
(\bigwedge^*TA, \iota_{df})\simeq CC^*(\mathrm{MF}(A,f))    
\end{equation}
and
\begin{equation}
(\Omega^*_{A}[[t]], td-df\wedge)\simeq CC^{S^1}_*(\mathrm{MF}(A,f)),
\end{equation}
where $\mathrm{MF}$ denotes the $\mathbb{Z}/2$-graded dg category of matrix factorizations. 
See \cite[Theorem 1.8, Corollary 1.13]{CT} for flavors of these results when $f$ has isolated singularities. 
\begin{conj}
Under the identification $CC^{S^1}_*\oplus CC^{S^1}_*\theta\simeq CC^{\mathbb{Z}/p}_*$ of Theorem 1.3 1), the HKR-type quasi-isomorphisms (6.1) and (6.2) intertwine the following two Frobenius $p$-linear graded multiplicative actions:
\begin{enumerate}[label=\arabic*)]
    \item The $\mathbb{Z}/p$-equivariant cap product $\bigcap^{\mathbb{Z}/p}: HH^*(\mathrm{MF}(A,f))\otimes HH^{\mathbb{Z}/p}_*(\mathrm{MF}(A,f))\rightarrow HH^{\mathbb{Z}/p}_*(\mathrm{MF}(A,f))$ 
    \item An action (after applying $\otimes_{\mathbf{k}[[t]]}\mathbf{k}[[t,\theta]]$)
\begin{equation}
H^*((\bigwedge^*TA, \iota_{df}))\otimes H^*(\Omega^*_{A}[[t]], td-df\wedge)\rightarrow H^*(\Omega^*_{A}[[t]], td+df\wedge)    
\end{equation}
where
\begin{itemize}
    \item `twisted functions' $f\in A/\iota_{df}(TA)$ act as multiplication by $f^p$;
    \item `twisted vector fields' $D\in \ker(\iota_{df}:TA\rightarrow A)/\iota_{df}(\bigwedge^2TA)$ act as 
    \begin{equation}
        i^{[p]}_{D}:=(\iota_{D^p}-\mathcal{L}_{D}^{p-1}\iota_D)t^{\frac{p-1}{2}}.
    \end{equation}
\end{itemize}
\end{enumerate}
\end{conj}
\subsection{Example: $A_N$ singularities}
Let us look at Conjecture 6.1 in the example of $A_N$-singularities for the prime $p=3$. Namely, we set $A=\mathbf{k}[[z]], f=z^N$ for some $3\nmid N$. This has an isolated singularity at $0$, and thus $H^*(\bigwedge^*TA,\iota_{df})$ is concentrated in degree $0$, where it is just the twisted functions $\mathbf{k}[[z]]/Nz^{n-1}$. As an algebra, this is generated by the element $z$. On the other hand, $H^*(\Omega^*_{A}[[t]], td-df\wedge)$ is concentrated in degree $1$ and is freely generated as a $\mathbf{k}[[t]]$-module by the elements $dz,zdz,\cdots,z^{N-2}dz$. One can easily calculate the action of $z$ on these generators (recall that $z$ acts as multiplication by $z^3$), and the result is that at the level of cohomology,
\begin{itemize}
    \item when $0\leq k<N-4$, 
    \begin{equation}[z^3\cdot z^kdz]=[z^{k+3}dz];
    \end{equation}
    \item for the remaining cases, we have 
    \begin{equation}
        [z^3\cdot z^{N-4}dz]=0,\;\; [z^3\cdot z^{N-3}dz]=-\frac{1}{N}[tdz],\;\; \textrm{and}\;\;[z^3\cdot z^{N-2}dz]=\frac{1}{N}[tzdz].
        \end{equation}
\end{itemize}
\subsubsection*{Computations in a minimal categorical model}
In this subsection, we study a purely algebraic computation of the $\mathbb{Z}/p$-equivariant cap product (cf. Definition 2.12) for a minimal categorical model of an $A_N$ singularity. We show that it agrees with the computation in (6.5) and (6.6). \par\indent
Now we describe the setup. Fix a coefficient field $\mathbf{k}$ of characteristic $p$. Consider the formal affine line $\mathbf{k}[[x]]$ equipped with the formal superpotential $W=r_2x^2+r_3x^3+\cdots$. \cite[Theorem 5.8]{Dyc} proved that the category of matrix factorization of $W$ has a compact generator whose endomorphism dg algebra is $A_{\infty}$-quasi-equivalent to the minimal unital $A_{\infty}$-algebra on the exterior algebra $\mathbf{k}\langle 1,\partial\rangle=\mathbf{k}1\oplus\mathbf{k}\partial$ (where $\partial$ has odd degree) given by
\begin{equation}
m_i(\partial,\cdots,\partial)=r_i, \;\;i\geq 2.    
\end{equation}
We denote this minimal $A_{\infty}$-algebra by $A_W$. Our goal is to study the $\mathbb{Z}/p$-equivariant cap product action for $\mathrm{MF}(W)$---and since Hochschild invariants are Morita invariant---it suffices to consider this action for the $A_{\infty}$-algebra $A_W$. \par\indent
Let $N\geq 2$ be the minimal integer such that $r_N\neq 0$, and further assume that $p\nmid N$. Then, $A_W$ is in fact $A_{\infty}$-quasi-equivalent to $A_{r_Nz^N}$. It suffices to show that one can make an invertible formal change of variables that transforms $W$ into $z^N$. As a first step, we make the change of variables
\begin{equation}
x\mapsto y-\frac{1}{N}\frac{r_{N+1}}{r_N}y^2.    
\end{equation}
Then the formal superpotential $W$ becomes
\begin{equation}
W=r_Ny^N+r_{N+2}'y^{N+2}+\cdots   
\end{equation}
Then, we can make a change of variables 
\begin{equation}
y\mapsto z-\frac{1}{N}\frac{r'_{N+2}}{r_N}z^3    
\end{equation}
to kill the $y^{N+2}$ terms, and so on. The infinite composition of these changes of variables is formally well-defined, and gives the desired quasi-equivalence of $A_W$ and $A_{r_Nz^N}$. Without loss of generality, we set $r_N=1$ (and all other $r_i$'s zero), and call $A_{z^N}$ \emph{the minimal $A_{\infty}$-algebra of an $A_N$-singularity}. \par\indent
The main result of this section is the following.
\begin{thm}
Let $p=3$. With respect to the generators ${}_3\tilde{R}^0_{\partial},\cdots,{}_3\tilde{R}^{N-2}_{\partial}$ of $HH^{\mathbb{Z}/3}_*(A_{z^N})$, the action $\bigcap^{\mathbb{Z}/3}(\varphi^1_1,-), \varphi^1_1\in HH^*(A_{z^N})$ (the definition of ${}_3\tilde{R}^0_{\partial},\cdots,{}_3\tilde{R}^{N-2}_{\partial}$ and $\varphi^1_1$ is given in (6.21) and (6.28)) is given by the matrix
\begin{equation}
\begin{pmatrix}
\mathbf{0}&\mathbf{0}&\mathbf{I}_{N-4}\\
\mathbf{0}&\mathbf{0}&\mathbf{0}\\
\begin{bmatrix}
-\frac{t}{N}&0\\
0&\frac{t}{N}
\end{bmatrix}&\mathbf{0}&\mathbf{0}
\end{pmatrix},
\end{equation}
where $\mathbf{I}$ denotes the identity matrix.
\end{thm}
In particular, under the identifications $\varphi^1_1\leftrightarrow [z]$ and ${}_3\tilde{R}^k_{\partial}\leftrightarrow [z^{N-2-k}dz]$, $0\leq k\leq N-2$, this agrees with the computation in (6.5) and (6.6).\par\indent
\textbf{Hochschild homology}. Let $\overline{CC}_*(A_{z^N})$ denote the normalized Hochschild chain complex of $A_{z^N}$, i.e. the quotient of $CC_*(A_{z^N})$ by all chains of the form $x_0|\cdots| 1|\cdots|x_n$, where the unit $1$ is not in the $0$-th position. Here, we have used the shorthand notation $a_0|a_1|\cdots|a_n$ for $a_0\otimes a_1\otimes\cdots\otimes a_n$. Since $A_{z^N}$ is strictly unital, the homology of $\overline{CC}_*(A_{z^N})$ computes $HH_*(A_{z^N})$. As a $\mathbb{Z}/2$-graded vector space, 
\begin{align}
\overline{CC}_*(A_{z^N})&= \overline{CC}_{even}(A_{z^N})\oplus \overline{CC}_{odd} (A_{z^N})\nonumber\\
&=\bigoplus_{k\geq 0} \mathbf{k}\langle 1|\overbrace{\partial|\cdots|\partial}^k\rangle\oplus\bigoplus_{k\geq 0} \mathbf{k}\langle\partial|\overbrace{\partial|\cdots|\partial}^k\rangle,
\end{align}
where $\mathbf{k}\langle\cdots\rangle$ denotes the free $\mathbf{k}$-vector space generated by element(s) in the bracket. 
The Hochschild differential is given by
\begin{equation}
b(1|\overbrace{\partial|\cdots|\partial}^k)=0\;,\quad b(\partial|\overbrace{\partial|\cdots|\partial}^k)=N\cdot 1|\overbrace{\partial|\cdots|\partial}^{k-N+1}.
\end{equation}
From (6.13) it is straightforward to conclude that 
\begin{equation}
HH_{even}(A_{z^N})=0\;,\quad HH_{odd}(A_{z^N})=\mathbf{k}^{\oplus N-1},
\end{equation}
and moreover, the Hochschild cycles 
\begin{equation}
\partial|\overbrace{\partial|\cdots|\partial}^k\;,\quad 0\leq k\leq N-2.    
\end{equation}
descend to a set of generators for $HH_{odd}(A_{z^N})$.\par\indent
\textbf{Hochschild cohomology}. Again, since $A_{z^N}$ is unital, we can use the normalized Hochschild cochain complex $\overline{CC}^*(A_{z^N})\subset CC^*(A_{z^N})$ to compute Hochschild cohomology. Recall that this is the subcomplex consisting of all multilinear maps $\phi$ which vanishes on elements of the form $x_1|\cdots|1|\cdots|x_n$. In particular,  $\phi\in \overline{CC}^*(A_{z^N})$ is uniquely determined by its value on the elements \begin{equation}
\overbrace{\partial|\cdots|\partial}^k\;,\quad k\geq 0.   
\end{equation} 
Therefore,
\begin{align}
\overline{CC}^*(A_{z^N})&=\overline{CC}^{even}(A_{z^N})\oplus\overline{CC}^{odd}(A_{z^N})\\
&=\prod_{k\geq 0}\mathrm{Hom}(\mathbf{k}\langle \overbrace{\partial|\cdots|\partial}^k\rangle, \mathbf{k}\langle 1\rangle)\oplus\prod_{k\geq 0}\mathrm{Hom}(\mathbf{k}\langle \overbrace{\partial|\cdots|\partial}^k\rangle, \mathbf{k}\langle \partial\rangle).
\end{align}
Let $\varphi_{1}^k$ denote the even Hochschild cochain that sends $\overbrace{\partial|\cdots|\partial}^k$ to $1$ and all other generators to $0$; similarly, let $\varphi_{\partial}^k$ denote the odd Hochschild cochain that sends $\overbrace{\partial|\cdots|\partial}^k$ to $\partial$ and all other generators to $0$. Then the Hochschild cochain differential is given by
\begin{equation}
\sum_{k\geq 0} a_k\varphi^k_1\mapsto0\;,\quad \sum_{k\geq 0} a_k\varphi^k_{\partial}\mapsto\sum_{n\geq 0}(N a_{n+1-N}) \varphi^n_1. 
\end{equation}
From (6.19) we conclude that
\begin{equation}
HH^{odd}(A_{z^N})=0\;,\quad HH^{even}(A_{z^N})=\mathbf{k}^{\oplus N-1},    
\end{equation}
and moreover the Hochschild cocycles 
\begin{equation}
\varphi^k_{1}\;,\quad 0\leq k\leq N-2    
\end{equation}
descend to a set of generators for $HH^{even}(A_{z^N})$. Finally, we note that under the cup (Yoneda) product, $(\varphi^1_{1})^{\cup k}=\varphi^k_{1}$.Thus, $HH^{even}(A_{z^N})$ is generated by $\varphi^1_{1}$ as an algebra.   \par\indent
\textbf{$\mathbb{Z}/p$-equivariant Hochschild homology}. We demonstrate the computation for $p=3$; in this subsection, $\mathbf{k}$ will be a field of characteristic $3$. First, we find explicit generators for the $3$-fold Hochschild homology of $A_{z^N}$. Since $A_{z^N}$ is unital, there exists an isomorphism 
\begin{equation}
 {}_3HH_*(A_{z^N})\cong HH_*(A_{z^N})=\begin{cases}
 0\;,\quad\textrm{if}\;*=\mathrm{even}\\
 \mathbf{k}^{\oplus N-1}\;,\quad\mathrm{if}\;*=\mathrm{odd}
 \end{cases}  
\end{equation}
induced by the quasi-isomorphism $\Phi^0_3$ of (2.34). Therefore, it suffices to find lifts of the chain level generators (6.15) under $\Phi^0_3$. For $0\leq k<N-1$, we define 
\begin{equation}
{}_3R^k_{\partial}:=\sum_{k_1+k_2+k_3=k} 1|\overbrace{\partial|\cdots|\partial}^{k_1}|1|\overbrace{\partial|\cdots|\partial}^{k_2}|\partial|\overbrace{\partial|\cdots|\partial}^{k_3}\in {}_3\overline{CC}_*(A_{z^N}).
\end{equation}
In the expression (6.23), the three entries $1,1,\partial$ not under the overbrace are the three bimodule entries. Since $k<N-1$, if we apply the $3$-fold Hochschild differential to ${}_3R^k_{\partial}$, the only terms involve applying an $m_2$ to an expression that contains a unit $1$. Hence, it is easy to see that ${}_3R^k_{\partial}, 0\leq k<N-1$ are cycles. Moreover, if $(k_1,k_2)\neq (0,0)$, then 
\begin{equation}
\Phi^0_3(1|\overbrace{\partial|\cdots|\partial}^{k_1}|1|\overbrace{\partial|\cdots|\partial}^{k_2}|\partial|\overbrace{\partial|\cdots|\partial}^{k_3})=0.  
\end{equation}
Thus, 
\begin{equation}
\Phi^0_3({}_3R^k_{\partial})=\Phi^0_3(1|1|\partial|\overbrace{\partial|\cdots|\partial}^{k})= \partial|\overbrace{\partial|\cdots|\partial}^{k},  
\end{equation}
and the cycles $\{{}_3R^k_{\partial}\}_{0\leq k<N-1}$ descend to a set of generators for ${}_3HH_*(A_{z^N})$.\par\indent
Now we proceed to find generators for $HH^{\mathbb{Z}/3}_*(A_{z^N})$. By Theorem 1.3 1), we know that abstractly $HH^{\mathbb{Z}/3}_{even}(A_{z^N})=\mathbf{k}[[t,\theta]]^{\oplus N-1}\theta$ and $HH^{\mathbb{Z}/3}_{odd}(A_{z^N})=\mathbf{k}[[t,\theta]]^{\oplus N-1}$, so it suffices to find lifts of the generators ${}_3R^k_{\partial}, 0\leq k<N-1$ under $HH^{\mathbb{Z}/3}_{odd}(A_{z^N})\xrightarrow{t,\theta=0} HH_{odd}(A_{z^N})$. \par\indent
Let $\tau\in\mathbb{Z}/3$ be the standard generator, which acts on $\overline{CC}^{\mathbb{Z}/3}_*(A_{z^N})$ via (2.35). Then, 
\begin{align}
(\tau-1){}_3R^{k}_{\partial}&=\sum_{k_1+k_2+k_3=k}(\partial|\overbrace{\partial|\cdots|\partial}^{k_1}|1|\overbrace{\partial|\cdots|\partial}^{k_2}|1|\overbrace{\partial|\cdots|\partial}^{k_3})-(1|\overbrace{\partial|\cdots|\partial}^{k_1}|1|\overbrace{\partial|\cdots|\partial}^{k_2}|\partial|\overbrace{\partial|\cdots|\partial}^{k_3})  \nonumber\\
&=b^3\Big(\frac{1}{N}\sum_{k_1+k_2+k_3=k+N-1}\partial|\overbrace{\partial|\cdots|\partial}^{k_1}|1|\overbrace{\partial|\cdots|\partial}^{k_2}|\partial|\overbrace{\partial|\cdots|\partial}^{k_3}\Big)
\end{align}
where $b^3$ denotes $3$-fold Hochschild differential. Next, we observe that
$$(1+\tau+\tau^2)\Big(\frac{1}{N}\sum_{k_1+k_2+k_3=k+N-1}\partial|\overbrace{\partial|\cdots|\partial}^{k_1}|1|\overbrace{\partial|\cdots|\partial}^{k_2}|\partial|\overbrace{\partial|\cdots|\partial}^{k_3}\Big)=$$
\begin{equation}b^3\Big(-\frac{1}{N^2}\sum_{k_1+k_2+k_3=k+2(N-1)}\partial|\overbrace{\partial|\cdots|\partial}^{k_1}|\partial|\overbrace{\partial|\cdots|\partial}^{k_2}|\partial|\overbrace{\partial|\cdots|\partial}^{k_3}\Big).
\end{equation}
Moreover, in (6.27), the term inside the bracket of $b^3(\cdots)$ is $\mathbb{Z}/3$-invariant. In particular, for $0\leq k<N-1$ 
\begin{align}
{}_3\tilde{R}^k_{\partial}:&={}_3R^k_{\partial}+\frac{1}{N}\sum_{k_1+k_2+k_3=k+N-1}\partial|\overbrace{\partial|\cdots|\partial}^{k_1}|1|\overbrace{\partial|\cdots|\partial}^{k_2}|\partial|\overbrace{\partial|\cdots|\partial}^{k_3}\theta \nonumber\\
&+\frac{1}{N^2}\sum_{k_1+k_2+k_3=k+2(N-1)}\partial|\overbrace{\partial|\cdots|\partial}^{k_1}|\partial|\overbrace{\partial|\cdots|\partial}^{k_2}|\partial|\overbrace{\partial|\cdots|\partial}^{k_3}t\quad\in \overline{CC}^{\mathbb{Z}/3}_{odd}(A_{z^N})
\end{align}
is a cycle in $\overline{CC}^{\mathbb{Z}/3}_*(A_{z^N})$ whose constant term is ${}_3R^k_{\partial}$. As a result,
$\{{}_3\tilde{R}^k_{\partial}\}_{0\leq k<N-1}$ descends to a set of generators for $HH^{\mathbb{Z}/3}_{*}(A_{z^N})$ as a free $\mathbf{k}[[t,\theta]]$-module.\par\indent
\textbf{The $\mathbb{Z}/p$-equivariant cap product}. Continuing the above discussion, we would like to determine the action $\bigcap^{\mathbb{Z}/3}$ of Definition 2.12 for the $A_{\infty}$-algebra $A_{z^N}, 3\nmid N$. Since $\bigcap^{\mathbb{Z}/3}$ is a multiplicative and Frobenius-linear action, and $HH^*(A_{z^N})$ is generated by $\varphi_{1}^1$ as an algebra, it suffices to compute $\bigcap^{\mathbb{Z}/3}(\varphi^1_{1},{}_3\tilde{R}^k_{\partial})$, for $0\leq k<N-1$. \par\indent
\emph{Proof of Theorem 6.2}. 
First of all, observe that by the $A_{\infty}$-structure equations (6.7) and $\varphi^1_{1}(\partial)=1$, on the chain level:
\begin{equation}
\bigcap^{\mathbb{Z}/3}(\varphi^1_{1}, x|\overbrace{\partial|\cdots|\partial}^{k_1}|y|\overbrace{\partial|\cdots|\partial}^{k_2}|z|\overbrace{\partial|\cdots|\partial}^{k_3})=x|\overbrace{\partial|\cdots|\partial}^{k_1-1}|y|\overbrace{\partial|\cdots|\partial}^{k_2-1}|z|\overbrace{\partial|\cdots|\partial}^{k_3-1},   
\end{equation}
where we set the term to be $0$ if $k_i-1<0$ for some $i\in\{1,2,3\}$. Now we discuss a few cases.\par\indent
\textbf{Case I: $3\leq k<N-1$}. As an immediate consequence of (6.29) and the definition of ${}_3\tilde{R}^k_{\partial}$ (6.28), we have
\begin{equation}
\bigcap^{\mathbb{Z}/3}(\varphi^1_{1}, {}_3\tilde{R}^k_{\partial})={}_3\tilde{R}^{k-3}_{\partial}, 3\leq k<N-1.     
\end{equation}
\textbf{Case II: $k=0$}. By (6.28) and (6.29), we compute that (note $\bigcap^{\mathbb{Z}/3}(\varphi^1_{1},{}_3R^0_{\partial})=0$)
\begin{align}
    \bigcap^{\mathbb{Z}/3}(\varphi^1_{1}, {}_3\tilde{R}^0_{\partial})&=\frac{1}{N}\sum_{k_1+k_2+k_3=N-4}\partial|\overbrace{\partial|\cdots|\partial}^{k_1}|1|\overbrace{\partial|\cdots|\partial}^{k_2}|\partial|\overbrace{\partial|\cdots|\partial}^{k_3}\theta \nonumber\\
&+\frac{1}{N^2}\sum_{k_1+k_2+k_3=2N-5}\partial|\overbrace{\partial|\cdots|\partial}^{k_1}|\partial|\overbrace{\partial|\cdots|\partial}^{k_2}|\partial|\overbrace{\partial|\cdots|\partial}^{k_3}t\quad\in \overline{CC}^{\mathbb{Z}/3}_{*}(A_{z^N}).
\end{align}
The claim is that after passing to cohomology,
\begin{equation}
 \bigcap^{\mathbb{Z}/3}(\varphi^1_{1}, {}_3\tilde{R}^0_{\partial})=-\frac{1}{N}{}_3\tilde{R}^{N-3}_{\partial}t\quad\in HH^{\mathbb{Z}/3}_{*}(A_{z^N}).
\end{equation}
We now discuss the proof of (6.32) in the case $N\equiv 1$ mod $3$; the case $N\equiv 2$ mod $3$ will be completely parallel. \par\indent
First, observe that 
\begin{equation}
\sum_{k_1+k_2+k_3=N-4}\partial|\overbrace{\partial|\cdots|\partial}^{k_1}|1|\overbrace{\partial|\cdots|\partial}^{k_2}|\partial|\overbrace{\partial|\cdots|\partial}^{k_3}\theta=b^3(Q_1\theta),
\end{equation}
where
\begin{equation}Q_1=-\sum_{k_1+k_2+k_3=N-3}k_3\cdot1|\overbrace{\partial|\cdots|\partial}^{k_1}|1|\overbrace{\partial|\cdots|\partial}^{k_2}|\partial|\overbrace{\partial|\cdots|\partial}^{k_3}.
\end{equation}
Next, we compute that
\begin{align}
(1+\tau+\tau^2)Q_1t&=-\sum_{k_1+k_2+k_3=N-3}\Big(k_3\cdot 1|\overbrace{\partial|\cdots|\partial}^{k_1}|1|\overbrace{\partial|\cdots|\partial}^{k_2}|\partial|\overbrace{\partial|\cdots|\partial}^{k_3}+k_1\cdot\partial|\overbrace{\partial|\cdots|\partial}^{k_1}|1|\overbrace{\partial|\cdots|\partial}^{k_2}|1|\overbrace{\partial|\cdots|\partial}^{k_3}  \nonumber\\
&\quad+k_2\cdot 1|\overbrace{\partial|\cdots|\partial}^{k_1}|\partial|\overbrace{\partial|\cdots|\partial}^{k_2}|1|\overbrace{\partial|\cdots|\partial}^{k_3}\Big).
\end{align}
Using the fact that $N=1\in\mathbb{F}_3$, we compute that
\begin{align}
b^3(Q_2)&=\sum_{k_1+k_2+k_3=N-3}\Big(k_1\cdot1|\overbrace{\partial|\cdots|\partial}^{k_1}|1|\overbrace{\partial|\cdots|\partial}^{k_2}|\partial|\overbrace{\partial|\cdots|\partial}^{k_3}-k_1\cdot\partial|\overbrace{\partial|\cdots|\partial}^{k_1}|1|\overbrace{\partial|\cdots|\partial}^{k_2}|1|\overbrace{\partial|\cdots|\partial}^{k_3}\Big)\nonumber\\
&\quad+\sum_{k_1+k_2+k_3=2N-5}\partial|\overbrace{\partial|\cdots|\partial}^{k_1}|\partial|\overbrace{\partial|\cdots|\partial}^{k_2}|\partial|\overbrace{\partial|\cdots|\partial}^{k_3},
\end{align}
where 
\begin{equation}
Q_2=-\sum_{k_1+k_2+k_3=2N-4}k_1\cdot\partial|\overbrace{\partial|\cdots|\partial}^{k_1}|1|\overbrace{\partial|\cdots|\partial}^{k_2}|\partial|\overbrace{\partial|\cdots|\partial}^{k_3}.    
\end{equation}
On the other hand,
\begin{align}
b^3(Q_3)&=\sum_{k_1+k_2+k_3=N-3}\Big(k_2\cdot1|\overbrace{\partial|\cdots|\partial}^{k_1}|1|\overbrace{\partial|\cdots|\partial}^{k_2}|\partial|\overbrace{\partial|\cdots|\partial}^{k_3}-k_2\cdot1|\overbrace{\partial|\cdots|\partial}^{k_1}|\partial|\overbrace{\partial|\cdots|\partial}^{k_2}|1|\overbrace{\partial|\cdots|\partial}^{k_3}\Big),   
\end{align}
where 
\begin{equation}
Q_3=-\sum_{k_1+k_2+k_3=2N-4}k_2\cdot1|\overbrace{\partial|\cdots|\partial}^{k_1}|\partial|\overbrace{\partial|\cdots|\partial}^{k_2}|\partial|\overbrace{\partial|\cdots|\partial}^{k_3}.
\end{equation}
Combining (6.36) and (6.38), and using that $k_1+k_2=1-k_3\in\mathbb{F}_3$, we obtain
\begin{align}
b^3(Q_2+Q_3)&=(1+\tau+\tau^2)Q_1+\sum_{k_1+k_2+k_3=2N-5}\partial|\overbrace{\partial|\cdots|\partial}^{k_1}|\partial|\overbrace{\partial|\cdots|\partial}^{k_2}|\partial|\overbrace{\partial|\cdots|\partial}^{k_3}\nonumber\\
&\quad+\sum_{k_1+k_2+k_3=N-3}1|\overbrace{\partial|\cdots|\partial}^{k_1}|1|\overbrace{\partial|\cdots|\partial}^{k_2}|\partial|\overbrace{\partial|\cdots|\partial}^{k_3}. 
\end{align}
By doing similar calculations one can show that
\begin{equation}
(\tau-1)(Q_2+Q_3)+\sum_{k_1+k_2+k_3=2N-4}\partial|\overbrace{\partial|\cdots|\partial}^{k_1}|1|\overbrace{\partial|\cdots|\partial}^{k_2}|\partial|\overbrace{\partial|\cdots|\partial}^{k_3}=b^3(Q_4),   
\end{equation}
where
\begin{equation}
Q_4=\sum_{k_1+k_2+k_3=3N-5}k_2\cdot\partial|\overbrace{\partial|\cdots|\partial}^{k_1}|\partial|\overbrace{\partial|\cdots|\partial}^{k_2}|\partial|\overbrace{\partial|\cdots|\partial}^{k_3},
\end{equation}
and 
\begin{equation}
(1+\tau+\tau^2)Q_4=\sum_{k_1+k_2+k_3=3N-5}\partial|\overbrace{\partial|\cdots|\partial}^{k_1}|\partial|\overbrace{\partial|\cdots|\partial}^{k_2}|\partial|\overbrace{\partial|\cdots|\partial}^{k_3}.
\end{equation}
Combining all of the above, we obtain that on the chain level,
\begin{equation}
\bigcap^{\mathbb{Z}/3}(\varphi^1_1,{}_3\tilde{R}^0_{\partial})+{}_3\tilde{R}^{N-3}t=d_{\overline{CC}^{\mathbb{Z}/3}}(Q_1\theta+(Q_2+Q_3)t-Q_4t\theta),
\end{equation}
which shows (6.32) for $N\equiv1$ mod $3$.\par\indent
\textbf{Case III: $k=1$}. Using similar calculations as above, one shows that at the level of cohomology,
\begin{equation}
\bigcap^{\mathbb{Z}/3}(\varphi^1_1,{}_3\tilde{R}^1_{\partial})=\frac{1}{N}{}_3\tilde{R}^{N-2}t.
\end{equation}
\textbf{Case IV: $k=2$}. Again, using similar calculations, one shows that at the level of cohomology,
\begin{equation}
\bigcap^{\mathbb{Z}/3}(\varphi^1_1,{}_3\tilde{R}^2_{\partial})=0.
\end{equation}
Putting these computations together gives Theorem 6.2.\qed 

\begin{appendices}
\renewcommand{\theequation}{A.\arabic{equation}}
\setcounter{equation}{0}
\section{Grading}
\subsection{$\mathbb{Z}$-grading on the monotone Fukaya category}
An $R$-linear \emph{pre-graded} $A_{\infty}$-category $\mathcal{A}$ is the data of a collection of objects $\mathrm{ob}\mathcal{A}$; for each pair of objects $X_0,X_1$, an $R$-module $\mathrm{hom}_{\mathcal{A}}(X_0,X_1)$; multilinear structure maps
\begin{equation}
\mu_{\mathcal{A}}^d: \mathrm{hom}_{\mathcal{A}}(X_0,X_1)\otimes\cdots\otimes \mathrm{hom}_{\mathcal{A}}(X_{d-1},X_d)\rightarrow \mathrm{hom}_{\mathcal{A}}(X_0,X_d),\;\;d\geq 1, X_0,\cdots,X_d\in\mathrm{ob}\mathcal{A},
\end{equation}
satisfying the $A_{\infty}$-relations (2.2). A \emph{$\mathbb{Z}$-grading} on a pre-graded $A_{\infty}$-category $\mathcal{A}$ is a map, for each pair of objects $X_0, X_1$, 
\begin{equation}
Gr: \mathrm{hom}_{\mathcal{A}}(X_0,X_1)\rightarrow \mathrm{hom}_{\mathcal{A}}(X_0,X_1)
\end{equation}
such that
\begin{equation}
[Gr, \mu^d_{\mathcal{A}}]=(2-d)\mu^d_{\mathcal{A}}.
\end{equation}
Let $X$ be a closed monotone symplectic manifold, and $R$ a base coefficient ring. Recall from section 3.2 that for each $\lambda\in R$ there is a $\mathbb{Z}/2$-graded $A_{\infty}$-category $\mathrm{Fuk}(X,R)_{\lambda}$ over $R$. In this section, following \cite[Appendix A.1]{Hug}, we upgrade $Fuk(X,R)_{\lambda}$ to a $\mathbb{Z}$-graded $A_{\infty}$-category $\mathrm{Fuk}(X,R[q,q^{-1}])_{\lambda}$ over $R[q,q^{-1}]$, where $q$ is a formal variable of degree $2$.   \par\indent
The objects of $\mathrm{Fuk}(X,R[q,q^{-1}])_{\lambda}$ are oriented spin Lagrangian submanifolds $L$ equipped with $R^*$-local systems, such that the Maslov index $2$ disk potential (3.3) is $\lambda$. \emph{For simplicity, for the rest of this section, we will ignore Hamiltonian perturbations when defining Lagrangian intersections and moduli spaces} (and assume they are already transverse), as they are standard and can be done in the same way as in section 3.2. \par\indent
Let $\mathcal{L}X\rightarrow X$ denote the oriented Lagrangian Grassmannian bundle over $X$. Any oriented Lagrangian submanifold $L\subset X$ gives rise to a section $s_L: L\rightarrow \mathcal{L}X|_L$ defined by $x\mapsto T_xL$. Let $L_0, L_1$ be two objects of $\mathrm{Fuk}(X,R[q,q^{-1}])_{\lambda}$. For simplicity, assume they are equipped with trivial local systems. For each $p\in L_0\cap L_1$, let $\mathcal{L}(L_0,L_1,p)$ denote the homotopy classes of paths $\tilde{p}:[0,1]\rightarrow \mathcal{L}_pX$ with $\tilde{p}(0)=T_pL_0$ and $\tilde{p}(1)=T_pL_2$. Define
\begin{equation}
CF^*(L_0,L_1)=\bigoplus_{p\in L_0\cap L_1} R\langle \mathcal{L}(L_0,L_1,0)\rangle.    
\end{equation}
Given a generator $\tilde{p}\in\mathcal{L}(L_0,L_1,p)$, we define its $\mathbb{Z}$-grading to be the Maslov index of $\tilde{p}$ composed with the canonical short path from $T_pL_1$ to $T_pL_0$, cf. \cite[Definition 1.8]{Aur}, which we denote by $\mathrm{ind}(\tilde{p})$. In other words, $Gr(\tilde{p}):=\mathrm{ind}(\tilde{p})\tilde{p}$. Note that the homotopy class of $\tilde{p}$ is uniquely determined by $\mathrm{ind}(\tilde{p})$. The $R[q,q^{-1}]$-module structure on $CF^*(L_0,L_1)$ is defined by letting $q$ send $\tilde{p}$ to the unique homotopy class of paths with Maslov index $\mathrm{ind}(\tilde{p})+2$.  \par\indent
We now define the $A_{\infty}$-structure maps. Let $L_0,\cdots,L_d$ be objects of $\mathrm{Fuk}(X,R[q,q^{-1}])_{\lambda}$ with trivial local systems, and let $\tilde{p}_i\in\mathcal{L}(L_{i-1},L_{i},p_i)$ for some $p_i\in L_{i-1}\cap L_i$, $i=1,2,\cdots,d$. For $p_0\in L_d\cap L_1$ and a homotopy class $\alpha\in \pi_2(X,L_0\cup\cdots\cup L_d)$, let 
\begin{equation}
\mathcal{M}(p_0,\cdots,p_d;\alpha)
\end{equation} 
denote the moduli space of pairs $(r,u)$, where $r\in \mathcal{R}^{d+1}, u: \mathcal{S}_r\rightarrow X$ a map of homotopy class $\alpha$ satisfying $du^{0,1}_J=0$ and the following asymptotic and Lagrangian boundary conditions: let $D^2\backslash\{z_0,\cdots,z_d\}$ be a representative of $\mathcal{S}_r$, then we require that $t$ sends $z_i$ to $p_i$ and the arc $\partial_{i-1}D^2=(z_{i-1},z_i)$ to $L_{i-1}$.
Then, the structure maps $\mu^d_q$ of $\mathrm{Fuk}(X,R[q,q^{-1}])_{\lambda}$ is defined as
\begin{equation}
\mu^d_q(\tilde{p}_1,\cdots,\tilde{p}_d):=\sum_{\substack{p_0\in L_d\cap L_0\\ u|\mathrm{ind}([u])=2-d}}\#\mathcal{M}(p_0,\cdots,p_d;[u])\tilde{p}_0(\tilde{p}_1,\cdots,\tilde{p}_d;[u]),
\end{equation}
where $\tilde{p}_0(\tilde{p}_1,\cdots,\tilde{p}_d;[u])\in \mathcal{L}(L_k,L_0,p_0)$ is given as follows. Fixing a trivialization of $u^*\mathcal{L}X$, we may view the concatenation $\gamma=(s_{L_0})|_{\partial_0D^2}\circ \tilde{p}_1\circ \cdots\circ \tilde{p}_d\circ (s_{L_d})|_{\partial_d D^2}$ as a path from $T_{p_0}L_0$ to $T_{p_0}L_d$ in $\mathcal{L}_{p_0}X$. We define $\tilde{p}_0(\tilde{p}_1,\cdots,\tilde{p}_d;[u])$ to be the unique homotopy class of paths from $T_{p_0}L_d$ to $T_{p_0}L_0$ such that $\tilde{p}_0(\tilde{p}_1,\cdots,\tilde{p}_d;[u])\circ \gamma$ has Maslov index $0$. By definition of the Maslov index, it is immediate that (A.6) defines an $R[q,q^{-1}]$-multilinear operation of degree $2-d$ (i.e. satisfies equation (A.3)). It is also clear from definition that after restricting to $q=1$ (which has the additional effect of collapsing the $\mathbb{Z}$-grading to a $\mathbb{Z}/2$-grading), we recover the $\mathbb{Z}/2$-graded $A_{\infty}$-category $\mathrm{Fuk}(X,R)_{\lambda}$ from section 3.2. \par\indent
There is also a $\mathbb{Z}$-graded $R[q,q^{-1}]$-linear version of quantum cohomology $(QH^*(X,R[q,q^{-1}]),\star_q)$, cf. \cite[section 3a]{SW} for a definition using the Morse chain model $CM^*(f,R[q,q^{-1}])$, which is the model we use in this paper. Recall that the underlying graded vector space of $CM^*(f,R[q,q^{-1}])$ is freely generated by the critical points of $f$ over $R[q,q^{-1}]$. There is a \emph{quantum $q$-connection} on $QH^*(X,R[q,q^{-1}])[[t]]$ given by
\begin{equation}
\nabla^{QH}_{tq\frac{d}{dq}}:=tq\frac{d}{dq}+c_1\star_q.
\end{equation}
We now describe a (chain level) $\mathbb{Z}$-graded $R[q,q^{-1}]$-linear open-closed map
\begin{equation}
OC_q: CC_*(\mathrm{Fuk}(X,R[q,q^{-1}])_{\lambda}\rightarrow CM^{*+n}(f,R[q,q^{-1}])    
\end{equation}
that recovers the usual open-closed map when restricted to $q=1$. Let $L_0,\cdots,L_d$ be objects of $\mathrm{Fuk}(X,R[q,q^{-1}])_{\lambda}$ and $\tilde{p}_i\in \mathcal{L}(L_{i-1},L_i,p_i)$ for some $p_i\in L_{i-1}\cap L_i$. Let $y_{out}$ be a critical point of $f$. Let $\mathcal{R}^1_{d+1}$ be the parameter space of disks with $d+1$ boundary marked points and $1$ interior marked point, up to automorphism; let $\mathcal{S}$ be the associated universal curve. For a homotopy class $\alpha \in \pi_2(X,L_0\cup\cdots\cup L_d)$, we define 
\begin{equation}
\mathcal{M}(p_0,\cdots,p_d,y_{out};\alpha)    
\end{equation}
to be the moduli space of pairs $r\in \mathcal{R}^1_{d+1}, u:\mathcal{S}_r\rightarrow X$ of homotopy class $\alpha$ satisfying $(du)^{0,1}_J=0$ and the following asymptotic, incidence and Lagrangian boundary conditions: let $D^2\backslash\{z_0,\cdots,z_d\}, y_0\in D^2$ be a representative of $\mathcal{S}_r$, then we require that $u(z_i)=p_i, u(\partial_i D^2)\subset L_i$ and that $u(y_0)\in W^u(y_{out})$. Define
\begin{equation}
OC_q(\tilde{p}_0\otimes \cdots\otimes\tilde{p}_d):=\sum_{\substack{y_{out}\in\mathrm{crit(f)}\\ [u]|\mathrm{ind}([u])=n-d-|y_{out}|}} \#\mathcal{M}(p_0,\cdots,p_d,y_{out};[u])q^{\frac{\mathrm{ind}(\tilde{p_0},\cdots,\tilde{p}_d)}{2}} y_{out},  
\end{equation}
where $\mathrm{ind}(\tilde{p_0},\cdots,\tilde{p}_d)$ is defined as the Maslov index of $\tilde{p}_0\circ(s_{L_0})|_{\partial_0D^2}\circ \tilde{p}_1\circ \cdots\circ \tilde{p}_d\circ (s_{L_d})|_{\partial_d D^2}$, viewed as a loop in $\mathcal{L}_{p_0}X$ after fixing a trivialization of $u^*\mathcal{L}X$. Using essentially the same argument as in \cite[Proposition 5.1]{Gan1}, taking into account the Maslov index of $\tilde{p}_i$, one can show that (A.9) is a $\mathbb{Z}$-graded $R[q,q^{-1}]$-linear chain map. One can furthermore enhance this into an $S^1$-equivariant version, following \cite{Gan2}, and obtain a $\mathbb{Z}$-graded $R[q,q^{-1}]$-linear chain map
\begin{equation}
OC_q^{S^1}:  CC^{S^1}_*(\mathrm{Fuk}(X,R[q,q^{-1}])_{\lambda}\rightarrow CM^{*+n}(f,R[q,q^{-1}])[[t]].  
\end{equation}
\begin{thm}(\cite[Theorem 6.3.5]{PS})
On the level of cohomology, $OC_q^{S^1}$ intertwines the Getzler-Gauss-Manin $q$-connection with the quantum $q$-connection (A.7).\qed    
\end{thm}
We refer the readers to \cite[Definition 3.16]{Hug} for a definition of the Getzler-Gauss-Manin $q$-connection on $HH^{S^1}_*(\mathrm{Fuk}(X,R[q,q^{-1}])_{\lambda}$.

\subsection{From the $q$-connection to the $t$-connection}
Recall that in section 1.3, we introduced the quantum $t$-connection (1.10)
\begin{equation}
\nabla^{QH}_{t^2\frac{d}{dt}}=t^2\frac{d}{dt}+\frac{\mu}{t}-\frac{c_1\star}{t^2},
\end{equation}
which has a second order pole at $t=0$. In fact, one can obtain the quantum $t$-connection from the quantum $q$-connection, and the process is crucially related to the issue of $\mathbb{Z}$-grading. Define the \emph{total degree operator} $\mathrm{Deg}: QH^*(X,R[q,q^{-1}])[[t]]\rightarrow QH^*(X,R[q,q^{-1}])[[t]]$ by
\begin{equation}
\mathrm{Deg}:=2(q\frac{d}{dq}+t\frac{d}{dt}+\mu).  
\end{equation}
As its name suggests, the effect of $\mathrm{Deg}$ applied to an element $\beta\in QH^*(X,R[q,q^{-1}])[[t]]$ is to multiply $\beta$ by its total $\mathbb{Z}$-grading (i.e. combining the grading from $H^*(X)$, from $q$ and from $t$), shifted by $n=\dim_{\mathbb{C}}X$ (cf. the definition of $\mu$ in (1.10)). Then, by the formulae (A.7) and (A.12), the $t$-connection is obtained as
\begin{equation}
\nabla^{QH}_{t^2\frac{d}{dt}}=\frac{1}{2}t\,\mathrm{Deg}-\nabla^{QH}_{tq\frac{d}{dq}}.    
\end{equation}
To be more precise, (1.10) is taken to be the restriction of (A.14) to $q=1$, but we will abuse terminology and call both the quantum $t$-connection.\par\indent
Similarly, the Getzler-Gauss-Manin $t$-connection can be obtained from its $q$-version via
\begin{equation}
\nabla^{GGM}_{t^2\frac{d}{dt}}:=\frac{1}{2}tGr^--\nabla^{GGM}_{tq\frac{d}{dq}}.    
\end{equation}
In (A.15), the \emph{total grading operator} $Gr^-$ is defined as
\begin{equation}
Gr^-:=\mathcal{L}_{Gr}+\Gamma+2t\frac{d}{dt},    
\end{equation}
where $\mathcal{L}_{Gr}$ denotes the Lie action of $Gr$ (A.2), viewed as a length $1$ Hochschild cochain, on Hochschild chains, and $\Gamma(x_0\otimes x_1\otimes\cdots\otimes x_k)=-kx_0\otimes x_1\otimes\cdots\otimes x_k$ is the length operator on Hochschild chains. The effect of $Gr^-$ on an element $\alpha=x_0\otimes x_1\otimes \cdots x_k t^l\in CC^{S^1}_*(\mathrm{Fuk}(X,R[q,q^{-1}])_{\lambda})$ is to multiply it by
\begin{equation}
|x_0|+|x_1|+\cdots+|x_k|+k+2l,  
\end{equation}
which is the total $\mathbb{Z}$-grading of $\alpha$ viewed as a cyclic chain (note that $k$ is absorbed into the reduced grading convention (2.21) on Hochschild chains). As a result of the above discussion, we obtain the following.
\begin{cor}
On the level of cohomology, $OC^{S^1}$ intertwines the quantum $t$-connection with the Getzler-Gauss-Manin $t$-connection.    
\end{cor}
\emph{Proof.} From the fact that $OC^{S^1}_q$ of (A.11) is $\mathbb{Z}$-graded, i.e. it intertwines the total grading operators $Gr^-$ with $\mathrm{Deg}$, Corollary A.2 immediately follows from Theorem A.1 and formulae (A.14), (A.15) for the $t$-connections (after setting $q=1$). \qed

\renewcommand{\theequation}{B.\arabic{equation}}
\setcounter{equation}{0}
\section{Proof of Lemma 4.9}
The proof is by induction and uses the fact that $S^{\infty}\times S^{\infty}$ is contractible. We set $C_{-1}=0$ and $C_0=0$. Suppose we have found $C_{i'}$ satisfying (4.60) for all $i'\leq i$. If $i$ is odd, we have
\begin{align}
&\partial\big(\delta(\Delta_i)-\sum_{i_1+i_2=i} \Delta_{i_1}\times \Delta_{i_2}-\sum_{\substack{i_1+i_2=i\\ i_1\;\mathrm{odd}}}\Delta_{i_1}\times (\tau-1)\Delta_{i_2}-(\tau\times \tau-1)C_i\big)\notag\\
=&(\tau\times \tau-1)\delta(\Delta_{i-1})-\sum_{\substack{i_1+i_2=i\\ i_1\;\mathrm{odd}}}\big((\tau-1)\Delta_{i_1-1}\times \Delta_{i_2}+(-1)^{i_1}\Delta_{i_1}\times (1+\tau+\cdots+\tau^{p-1})\Delta_{i_2-1}\big)\notag\\
&-\sum_{\substack{i_1+i_2=i\\ i_1\;\mathrm{even}}}\big((1+\tau+\cdots+\tau^{p-1})\Delta_{i_1-1}\times \Delta_{i_2}+(-1)^{i_1}\Delta_{i_1}\times (\tau-1)\Delta_{i_2-1}\big)\notag\\
&-\sum_{\substack{i_1+i_2=i\\ i_1\;\mathrm{odd}}}(\tau-1)\Delta_{i_1-1}\times (\tau-1)\Delta_{i_2}-(\tau\times \tau-1)\partial C_i.
\end{align}
By induction hypothesis, we have
\begin{equation}
(\tau\times \tau-1)\partial C_i=(\tau\times \tau-1)\Big(\delta(\Delta_{i-1})-\sum_{\substack{i_1+i_2=i-1\\ i_k\;\mathrm{even}}} \Delta_{i_1}\times \Delta_{i_2}-\sum_{\substack{i_1+i_2=i-1\\ i_1\;\mathrm{odd}}}\sum_{0\leq k<j\leq p-1}\tau^k\Delta_{i_1}\times \tau^j\Delta_{i_2}\Big).
\end{equation}
Combining (B.1) and (B.2), we have
\begin{align}
&\partial\big(\delta(\Delta_i)-\sum_{i_1+i_2=i} \Delta_{i_1}\times \Delta_{i_2}-\sum_{\substack{i_1+i_2=i\\ i_1\;\mathrm{odd}}}\Delta_{i_1}\times (\tau-1)\Delta_{i_2}-(\tau\times \tau-1)C_i\big)\notag\\
=&-\sum_{\substack{i_1+i_2=i\\ i_1\;\mathrm{odd}}}\big((\tau-1)\Delta_{i_1-1}\times \Delta_{i_2}+(-1)^{i_1}\Delta_{i_1}\times (1+\tau+\cdots+\tau^{p-1})\Delta_{i_2-1}\big)\notag\\
&-\sum_{\substack{i_1+i_2=i\\ i_1\;\mathrm{even}}}\big((1+\tau+\cdots+\tau^{p-1})\Delta_{i_1-1}\times \Delta_{i_2}+(-1)^{i_1}\Delta_{i_1}\times (\tau-1)\Delta_{i_2-1}\big)\notag\\
&-\sum_{\substack{i_1+i_2=i\\ i_1\;\mathrm{odd}}}(\tau-1)\Delta_{i_1-1}\times (\tau-1)\Delta_{i_2}+(\tau\times \tau-1)\Big(\sum_{\substack{i_1+i_2=i-1\\ i_k\;\mathrm{even}}} \Delta_{i_1}\times \Delta_{i_2}+\sum_{\substack{i_1+i_2=i-1\\ i_1\;\mathrm{odd}}}\sum_{0\leq k<j\leq p-1}\tau^k\Delta_{i_1}\times \tau^j\Delta_{i_2}\Big)\notag\\
=&-\sum_{\substack{i_1+i_2=i\\ i_1\;\mathrm{odd}}}(-1)^{i_1}\Delta_{i_1}\times (1+\tau+\cdots+\tau^{p-1})\Delta_{i_2-1}-\sum_{\substack{i_1+i_2=i\\ i_1\;\mathrm{even}}}(1+\tau+\cdots+\tau^{p-1})\Delta_{i_1-1}\times \Delta_{i_2}\notag\\
&+(\tau\times \tau-1)\sum_{\substack{i_1+i_2=i-1\\ i_1\;\mathrm{odd}}}\sum_{0\leq k<j\leq p-1}\tau^k\Delta_{i_1}\times \tau^j\Delta_{i_2}\notag\\
=&\;0.
\end{align}
Since $S^{\infty}\times S^{\infty}$ is contractible, there exists $C_{i+1}$ such that
\begin{equation}
\partial C_{i+1}=\delta(\Delta_i)-\sum_{i_1+i_2=i} \Delta_{i_1}\times \Delta_{i_2}-\sum_{\substack{i_1+i_2=i\\ i_1\;\mathrm{odd}}}\Delta_{i_1}\times (\tau-1)\Delta_{i_2}-(\tau\times \tau-1)C_i.
\end{equation}
If $i$ is even, we have
\begin{align}
&\partial\big(\delta(\Delta_i)-\sum_{\substack{i_1+i_2=i\\ i_k\;\mathrm{even}}} \Delta_{i_1}\times \Delta_{i_2}-\sum_{\substack{i_1+i_2=i\\ i_1\;\mathrm{odd}}}\sum_{0\leq k<j\leq p-1}\tau^k\Delta_{i_1}\times \tau^j\Delta_{i_2}-(1+\tau\times \tau+\cdots+\tau\times \tau^{p-1})C_i\big)\notag\\
=&(1+\tau\times \tau+\cdots+\tau\times \tau^{p-1})\delta(\Delta_{i-1})-(1+\tau\times \tau+\cdots+\tau\times \tau^{p-1})\partial C_i\notag\\
&-\sum_{\substack{i_1+i_2=i\\ i_k\;\mathrm{even}}}\big((1+\tau+\cdots+\tau^{p-1})\Delta_{i_1-1}\times \Delta_{i_2}+(-1)^{i_1}\Delta_{i_1}\times (1+\tau+\cdots+\tau^{p-1})\Delta_{i_2-1}\big)\notag\\
&-\sum_{\substack{i_1+i_2=i\\ i_1\;\mathrm{odd}}}\sum_{0\leq k<j\leq p-1}\big(\tau^k(\tau-1)\Delta_{i_1-1}\times \tau^j\Delta_{i_2}+(-1)^{i_1}\tau^k\Delta_{i_1}\times \tau^j(\tau-1)\Delta_{i_2-1}\big).
\end{align}
By induction hypothesis, we have
\begin{equation}
(1+\tau\times \tau+\cdots+\tau\times \tau^{p-1})\partial C_i=(1+\tau\times \tau+\cdots+\tau\times \tau^{p-1})\big(\delta(\Delta_{i-1})-\sum_{i_1+i_2=i-1} \Delta_{i_1}\times \Delta_{i_2}-\sum_{\substack{i_1+i_2=i-1\\ i_1\;\mathrm{odd}}}\Delta_{i_1}\times (\tau-1)\Delta_{i_2}\big).
\end{equation}
Combining (B.5) and (B.6), we obtain
\begin{align}
&\partial\big(\delta(\Delta_i)-\sum_{\substack{i_1+i_2=i\\ i_k\;\mathrm{even}}} \Delta_{i_1}\times \Delta_{i_2}-\sum_{\substack{i_1+i_2=i\\ i_1\;\mathrm{odd}}}\sum_{0\leq k<j\leq p-1}\tau^k\Delta_{i_1}\times \tau^j\Delta_{i_2}-(1+\tau\times \tau+\cdots+\tau\times \tau^{p-1})C_i\big)\notag\\
=&(1+\tau\times \tau+\cdots+\tau\times \tau^{p-1})\big(\sum_{i_1+i_2=i-1} \Delta_{i_1}\times \Delta_{i_2}+\sum_{\substack{i_1+i_2=i-1\\ i_1\;\mathrm{odd}}}\Delta_{i_1}\times (\tau-1)\Delta_{i_2}\big)\notag\\
&-\sum_{\substack{i_1+i_2=i\\ i_k\;\mathrm{even}}}\big((1+\tau+\cdots+\tau^{p-1})\Delta_{i_1-1}\times \Delta_{i_2}+(-1)^{i_1}\Delta_{i_1}\times (1+\tau+\cdots+\tau^{p-1})\Delta_{i_2-1}\big)\notag\\
&-\sum_{\substack{i_1+i_2=i\\ i_1\;\mathrm{odd}}}\sum_{0\leq k<j\leq p-1}\big(\tau^k(\tau-1)\Delta_{i_1-1}\times \tau^j\Delta_{i_2}+(-1)^{i_1}\tau^k\Delta_{i_1}\times \tau^j(\tau-1)\Delta_{i_2-1}\big).
\end{align}
Next, observe that
\begin{equation}
\sum_{0\leq k<j\leq p-1}\tau^k(\tau-1)\Delta_{i_1}\times \tau^j\Delta_{i_2}=(1+\tau\times \tau+\cdots+\tau\times \tau^{p-1})(\Delta_{i_1}\times \Delta_{i_2})-\Delta_{i_1}\times (1+\tau+\cdots+\tau^{p-1})\Delta_{i_2}
\end{equation}
and
\begin{equation}
\sum_{0\leq k<j\leq p-1}\tau^k\Delta_{i_1}\times \tau^j(\tau-1)\Delta_{i_2}=(1+\tau+\cdots+\tau^{p-1})\Delta_{i_1}\times \Delta_{i_2}-(1\times \tau)(1+\tau\times \tau+\cdots+\tau\times \tau^{p-1})(\Delta_{i_1}\times \Delta_{i_2}).
\end{equation}
Plugging (B.8) and (B.9) into (B.7) we obtain
\begin{equation}
\partial\big(\delta(\Delta_i)-\sum_{\substack{i_1+i_2=i\\ i_k\;\mathrm{even}}} \Delta_{i_1}\times \Delta_{i_2}-\sum_{\substack{i_1+i_2=i\\ i_1\;\mathrm{odd}}}\sum_{0\leq k<j\leq p-1}\tau^k\Delta_{i_1}\times \tau^j\Delta_{i_2}-(1+\tau\times \tau+\cdots+\tau\times \tau^{p-1})C_i\big)=0
\end{equation}
and by contractibility there exists some $C_{i+1}$ with
\begin{equation}
\partial C_{i+1}=\delta(\Delta_i)-\sum_{\substack{i_1+i_2=i\\ i_k\;\mathrm{even}}} \Delta_{i_1}\times \Delta_{i_2}-\sum_{\substack{i_1+i_2=i\\ i_1\;\mathrm{odd}}}\sum_{0\leq k<j\leq p-1}\tau^k\Delta_{i_1}\times \tau^j\Delta_{i_2}-(1+\tau\times \tau+\cdots+\tau\times \tau^{p-1})C_i.
\end{equation}
This concludes the induction.\qed

\end{appendices}


\end{document}